\documentclass[twocolumn]{svjour3nojournal}          
\smartqed  
\usepackage{microtype}
\usepackage{graphicx}
\usepackage{comment}
\graphicspath{ {./figures/} {./} }
\usepackage{newtxtext,newtxmath}
\usepackage[T1]{fontenc}
\usepackage[utf8]{inputenc}
\usepackage[english]{babel}
\usepackage{bm}
\usepackage[toc,title]{appendix}
\usepackage{url}
\usepackage{hyperref}
\usepackage{float}
\usepackage[labelfont=bf,font=small,labelsep=space]{caption}
\usepackage{subcaption}
\usepackage{xfrac}
\usepackage[dvipsnames]{xcolor}
\usepackage{cite,comment}
\usepackage{tikz}
    \usetikzlibrary{arrows}
    \usetikzlibrary{calc}
    \usetikzlibrary{cd}
    \tikzstyle{every picture}+=[remember picture]
\usepackage{subfiles}
\usepackage{commands}
\usepackage{multirow}
\usepackage{standalone}
\usepackage{pgfplots}
\usetikzlibrary{
        pgfplots.fillbetween
    }
\pgfplotsset{compat=1.16}

\makeatletter
\newcommand*{\bigcdot}{}
\DeclareRobustCommand*{\bigcdot}{%
  \mathbin{\mathpalette\bigcdot@{}}%
}
\newcommand*{\bigcdot@scalefactor}{.7}
\newcommand*{\bigcdot@widthfactor}{1.15}
\newcommand*{\bigcdot@}[2]{%
  \sbox0{$#1\vcenter{}$}
  \sbox2{$#1\cdot\m@th$}%
  \hbox to \bigcdot@widthfactor\wd2{%
    \hfil
    \raise\ht0\hbox{%
      \scalebox{\bigcdot@scalefactor}{%
        \lower\ht0\hbox{$#1\bullet\m@th$}%
      }%
    }%
    \hfil
  }%
}
\makeatother

\newcommand\numberthis{\addtocounter{equation}{1}\tag{\theequation}}
\providecommand{\numberTblEq}[1]{\refstepcounter{tblEqCounter}\label{#1}\thetag{\thetblEqCounter}}

\newcommand{\dbar}{d\hspace*{-0.08em}\bar{}\hspace*{0.1em}}
\newcommand{\ul}{\mathbf}
\newcommand{\R}{\mathbb{R}}
\DeclareMathAlphabet\gothic{U}{euf}{m}{n}
\DeclareMathOperator\Span{Span}

\newcommand\cA{{\omega}}
\newcommand\cF{{\mathcal F}}
\newcommand\ve{\varepsilon}
\newcommand\bM{{\mathbb M}}

\journalname{Journal of Mathematical Imaging and Vision}
\begin{document}

\title{Geodesic Tracking via New Data-driven Connections of Cartan Type for Vascular Tree Tracking
}


\author{Nicky~van den Berg \and Bart~Smets 
\and Gautam Pai \and \mbox{Jean-Marie Mirebeau} \and \mbox{Remco Duits}
}


\institute{Nicky J. van den Berg \and Bart M.N. Smets \and Gautam Pai \and Remco Duits \at
              Department of Mathematics and Computer Science, Eindhoven University of Technology, The Netherlands \\
              \email{n.j.v.d.berg@tue.nl, b.m.n.smets@tue.nl, g.pai@tue.nl, 
              r.duits@tue.nl}           
           \and
           Jean-Marie Mirebeau \at
           Department of Mathematics, Centre Borelli, ENS Paris-Saclay, CNRS, University Paris-Saclay, 91190, Gif-sur-Yvette, France, \email{jm.mirebeau@gmail.com}
}

\date{November 9, 2023}

\maketitle

\begin{abstract}
We introduce a data-driven version of the plus Cartan connection on the homogeneous space $\mathbb{M}_2$ of 2D positions and orientations. We formulate a theorem that describes all shortest and straight curves (parallel velocity and parallel momentum, respectively) with respect to this new data-driven connection and corresponding Riemannian manifold. Then we use these shortest curves for geodesic tracking of complex vasculature in multi-orientation image representations defined on $\mathbb{M}_{2}$. The data-driven Cartan connection characterizes the Hamiltonian flow of all geodesics. It also allows for improved adaptation to curvature and misalignment of the (lifted) vessel structure that we track via globally optimal geodesics. We compute these geodesics numerically via steepest descent on distance maps on $\mathbb{M}_2$ that we compute by a new modified anisotropic fast-marching method. 

Our experiments range from tracking single blood vessels with fixed endpoints to tracking complete vascular trees in retinal images. Single vessel tracking is performed in a single run in the multi-orientation image representation, where we project the resulting geodesics back onto the underlying image. The complete vascular tree tracking requires only two runs and avoids prior segmentation, placement of extra anchor points, and dynamic switching between geodesic models. 

Altogether we provide a geodesic tracking method using a single, flexible, transparent, data-driven geodesic model providing globally optimal curves which correctly follow highly complex vascular structures in retinal images. 

All~experiments~in this article~can~be~reproduced~via~documented  \emph{Mathematica}~notebooks~available~
at~\cite{githubNicky}.


\end{abstract}


\section{Introduction}\label{sec:Introduction}

Retinal images are often used to examine the vascular system with optical scanning devices that image the vasculature in the retina noninvasively. 
The vasculature in the eye is known to be typically representative of the vasculature throughout the body. This allows doctors to monitor the circulatory system and aids in the diagnosis of different kinds of diseases like diabetes, hypertension \cite{bekkers2017phd,weiler2015arteriole,sasongko2016retinal} and Alzheimer's disease \cite{colligris2018ocular}. Typically, high levels of tortuosity in the vasculature are a biomarker for such diseases \cite{bekkers2015curvature,Kalitzeos2013retinal,Cheung2012retinal,Sasongko2011retinal}. Successful automatic vessel tracking detects complex vasculature and aids the effective diagnosis of such diseases. 
Here,  geometric models come into play via geodesic tracking methods where geodesics are the shortest paths that follow the biological blood vessels. They help in tracking and subsequent analysis of the vascular tree in the retina originating from the optic nerve \cite{chen2018dynamic,Liu2021tubular,Liu2019vessel,bekkers2015pde}.

Geodesic tracking has been extensively studied where many prevalent approaches perform the tracking in the standard 2D image domain \cite{cootes1995active,Kimmel1997GeodesicActiveControus,cohen2001multiple,sapiro2006geometric,peyre2010geodesic,Chen2021GeneralizedAssymmetric}. For many methods in this category, calculating the geodesics in $\mathbb{R}^2$ leads to certain difficulties in accurately following the 
blood vessel. For example, one common difficulty is the inaccurate tracking of crossing structures and bifurcations. This has motivated methods that aim to lift the image function to higher dimensional spaces. For example, the space of positions and orientations \cite{pechaud2009brain,Portegies2019Brain} or radius-lifted spaces \cite{PhDDaChen} where the lifting yields the benefit of disentangling seemingly complex crossing structures in the retinal images. 

\begin{figure}
    \centering
    \begin{subfigure}[b]{0.22\textwidth}
        \includegraphics[width=\textwidth]{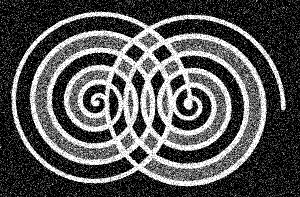}
        \caption{Image $f$.}
    \end{subfigure}
    \hfill
    \begin{subfigure}[b]{0.22\textwidth}
        \includegraphics[width=\textwidth]{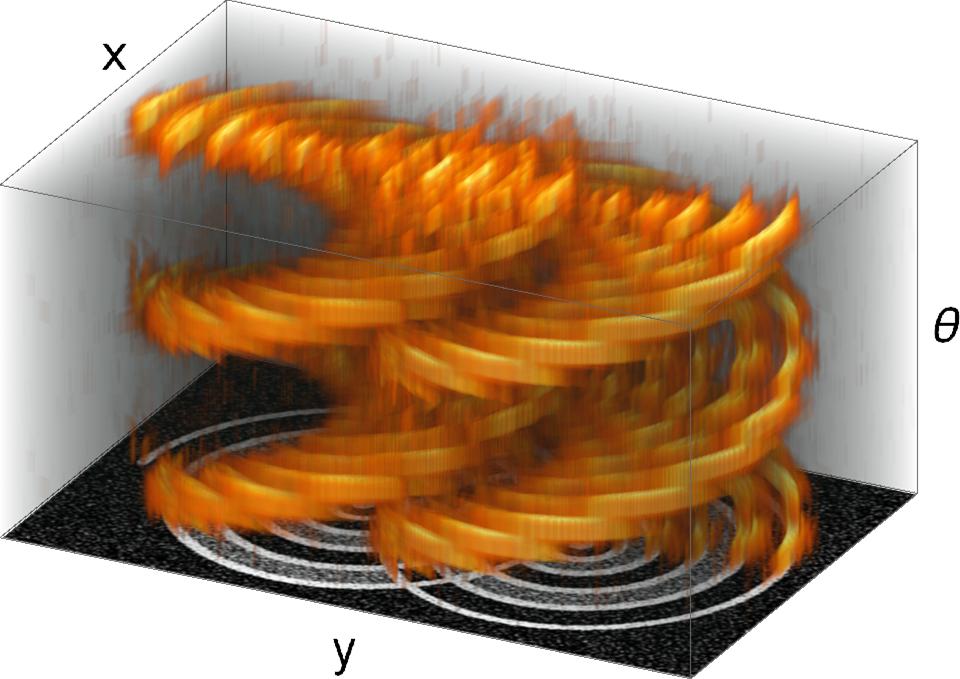}
        \caption{Orientation Score $U$ belonging to image $f$.}
        \label{fig:LiftedOrientationScores}
    \end{subfigure}
    \hfill
    \begin{subfigure}[b]{0.48\textwidth}
    \centering{
        \includegraphics[width=0.25\textwidth]{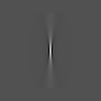}
        \includegraphics[width=0.25\textwidth]{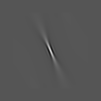}
            }
        \caption{Visualization of the real part of the cake-wavelet, with orientations $\theta=0$ (left) and $\theta=3\pi/16$ (right), used to lift image $f$ to the orientation score.}
        \label{fig:cakewavelets}
    \end{subfigure}
    \caption{(a) Visualization of a grayscale image $f:\R^2 \to \R$ and (b) its corresponding  orientation score $U: \mathbb{M}_{2} \to \R$ in the space of positions and orientations $\mathbb{M}_2$ given by \eqref{score}, using a standard cakewavelet as depicted in \ref{fig:cakewavelets}. We use a volume-rendering where the orange spirals indicate data-points $\ul{p}=(x,y,\theta)\in \mathbb{M}_2$ with high amplitudes $|U(\ul{p})|$. }
    \label{fig:OrientationScores}
\end{figure}

In this article, we focus on the methods
\cite{benmansour2011tubular,bekkers2014multi,PhDDaChen,bekkers2015pde} that perform the geodesic tracking in the 3D-space of positions and orientations $\mathbb{M}_{2}$. This is based on lifted images, or so-called \textit{orientation scores}. The well-known benefit of this lifted approach is that lines involved in crossings are manifestly disentangled in $\mathbb{M}_{2}$. As visualized in Figure~\ref{fig:OrientationScores}, the crossing circles in the image become disjoint spirals (cf. Fig. \ref{fig:LiftedOrientationScores}) in the homogeneous space of positions and orientations.

\begin{figure}
    \centering
    \begin{subfigure}[t]{0.225\textwidth}
    \centering
        \includegraphics[width=0.75\textwidth]{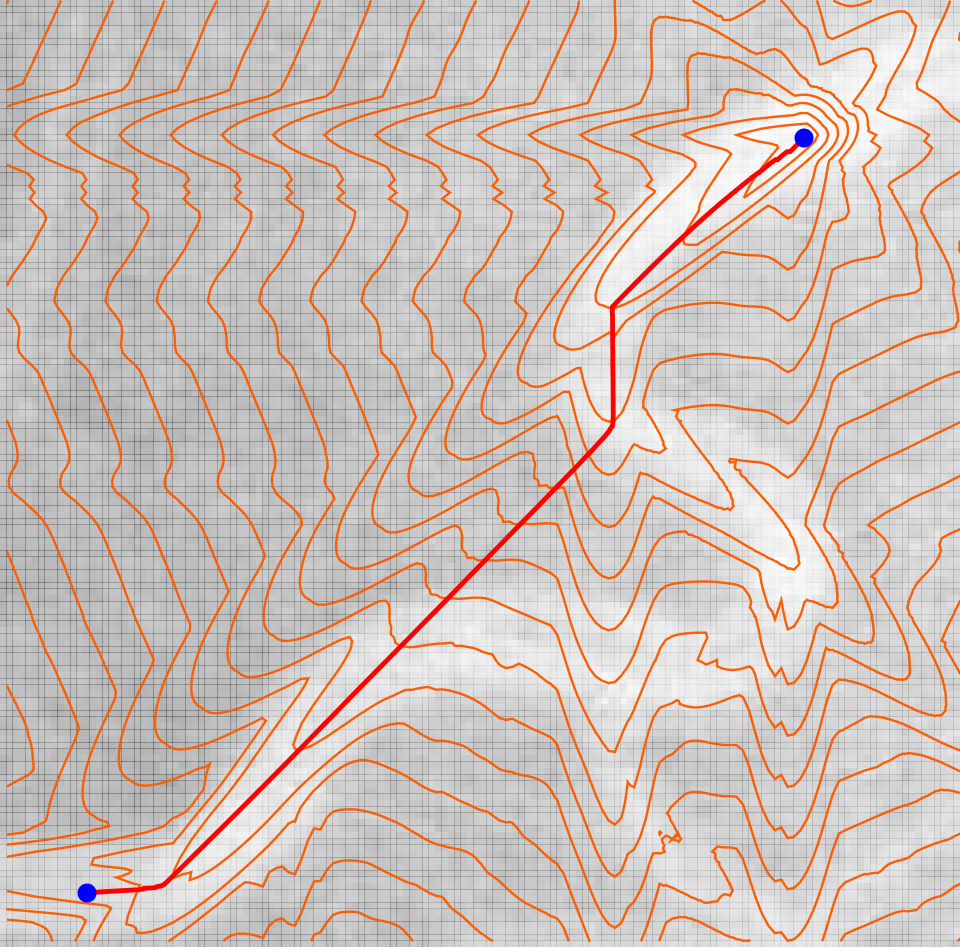}
        \caption{Previous geodesic tracking.}
    \end{subfigure}
    \hfill
    \begin{subfigure}[t]{0.225\textwidth}
    \centering
        \includegraphics[width=0.75\textwidth]{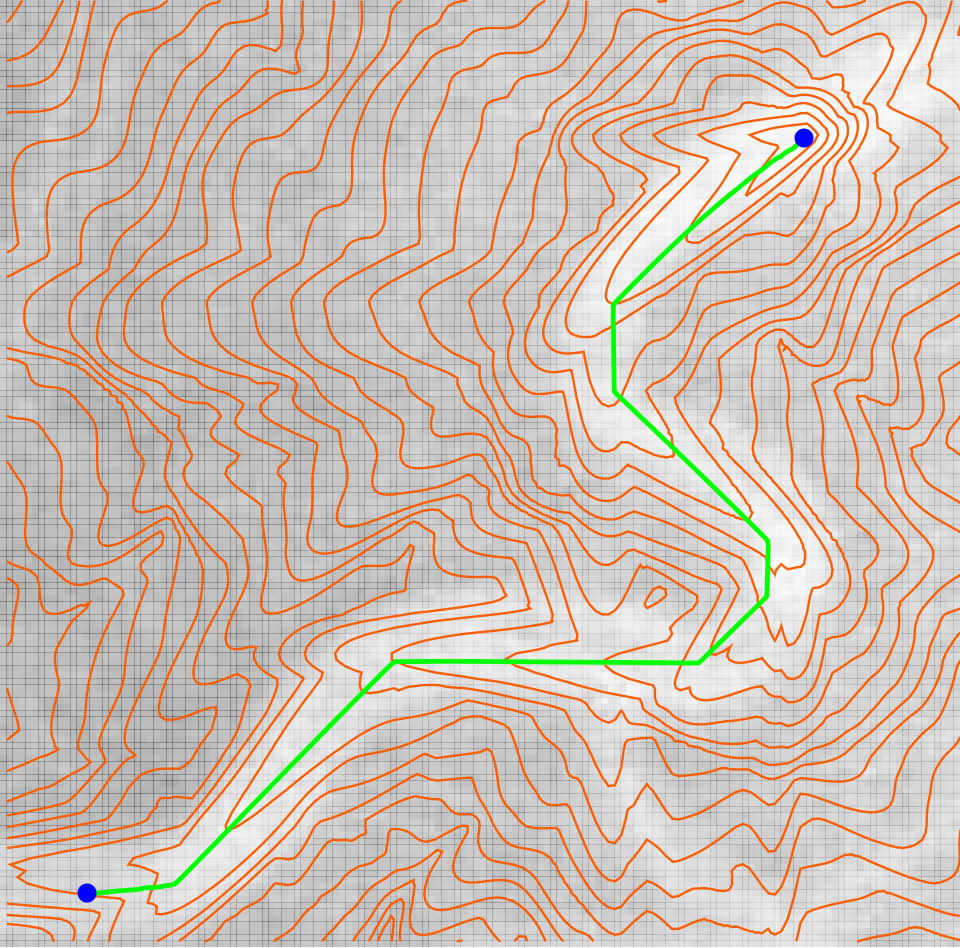}
        \caption{Proposed geodesic tracking.}
    \end{subfigure}
    \caption{Orientations sampling bias in geodesic tracking. Sampling bias can lead to wrong tracking results, and our new model will overcome this as we will show later in more detail (Fig. \ref{fig:S-curve}).}
    \label{fig:deviationFromHorizontality}
\end{figure}

However, practical considerations of working in (the domain $\mathbb{M}_2$ of) orientation scores, like memory reduction and enabling low computation times, result in some undesirable effects. For example, using a limited number of orientations leads to imperfections in the computation of the orientation scores. Hence, some vessels can be assigned a near angular coordinate that may not reflect their true orientation, and therefore does not align with the vessel data correctly. We denote this problem as 
`misalignment' (also referred to as \textit{deviation from horizontality} \cite{franken2009crossing}).
Moreover, considering a limited number of orientations results in a sampling bias on orientations, and thereby the possibility of missing high curvature regions yielding 
poor curvature adaptation (cf.~Fig.~\ref{fig:deviationFromHorizontality}).


In this article, we provide a novel, data-driven tracking model that improves upon existing geodesic tracking methods. Our model demonstrates an improved curvature adaptation, reduces misalignment, and exhibits a high degree of geometric interpretability.

We will aim for a \textit{single} geometric Finslerian model to deal with complex vasculature without requiring heavy pre-processing (e.g.~placement of anchor points, pre-skeletonization) and associated extra parameters, and without suffering from the `cusp problems' reported in \cite{charlot,duits2014association,bekkers2015pde}. 

The cusp problem is tackled by creating an asymmetric Finslerian model\footnote{In Finsler geometry \cite{bao2000Introduction,duits2018optimal}, the norm of tangent vectors may not be induced by an inner product. Recall that in the Riemannian setting, one does have $\mathcal{F}(\p,\dot{\p})=\sqrt{\mathcal{G}_\p(\dot{\p},\dot{\p})}$.} $(\mathbb{M}_2,\mathcal{F}^U)$ extension of the data-driven Riemannian manifold similar to the much less data-driven techniques in \cite{duits2018optimal}. For a quick impression of such a `cusp' in a spatially projected sub-Riemannian geodesic, see Fig.~\ref{fig:cusp} in Appendix~\ref{app:CostFunction}. Clearly, cusps are undesirable for vascular tracking, and an asymmetric Finslerian version of the Riemannian manifold tackles this problem. Intuitively, cusps in spatial projections of sub-Riemannian geodesics arise sometimes as optimal paths of a `Reed-Shepp' car (imagine a car driving along the geodesic track) \cite{duits2018optimal,reeds1990optimal} where the car was required to use its reverse gear to follow the optimal path. In the asymmetric Finslerian model we turn off the reverse gear of the car, while allowing for `in-place rotations' see Fig.~\ref{fig:InPlaceRotation}.

Pre-processing techniques for geodesic tracking  such as pre-skeletonization and iterative placement of anchor/key points are typically used in conjunction with Bézier curves \cite{absil2009optimization} or splines on Lie-groups \cite{absil2009optimization,bekkers2018nilpotent}, but often require additional parameters and fine-tuning. 
Specifically, extensive use of anchor points implies that anchors get relatively close to each other, and then the choice of geometric model in between becomes increasingly less relevant (even non-data-driven sub-Riemannian distance approximations suffice as shown in the work by Bekkers et al. \cite{bekkers2018nilpotent}). 
As a result, this reduces the geometric interpretability of the overall model. 
In this work, we therefore aim for a single geometrical model. Therefore, we will not use pre-processing, pre-skeletonization \cite{PhDDaChen}, multiple anchor points \cite{chen2016vessel}, and connectivity by perceptional grouping \cite{cohen2001multiple,abbasi2017retrieving,bekkers2018nilpotent}, even though these techniques are 
theoretically interesting and applicable.

\begin{figure}
    \centering
    \includegraphics[width=0.45\textwidth]{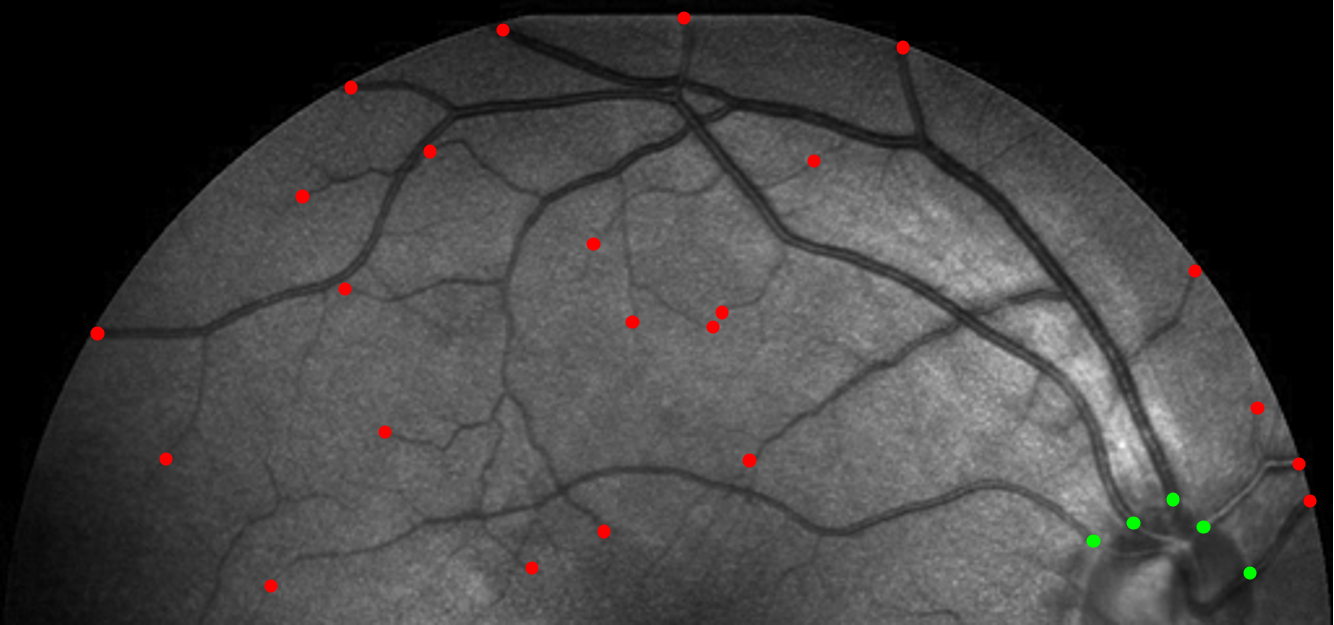}
    
    \includegraphics[width=0.45\textwidth]{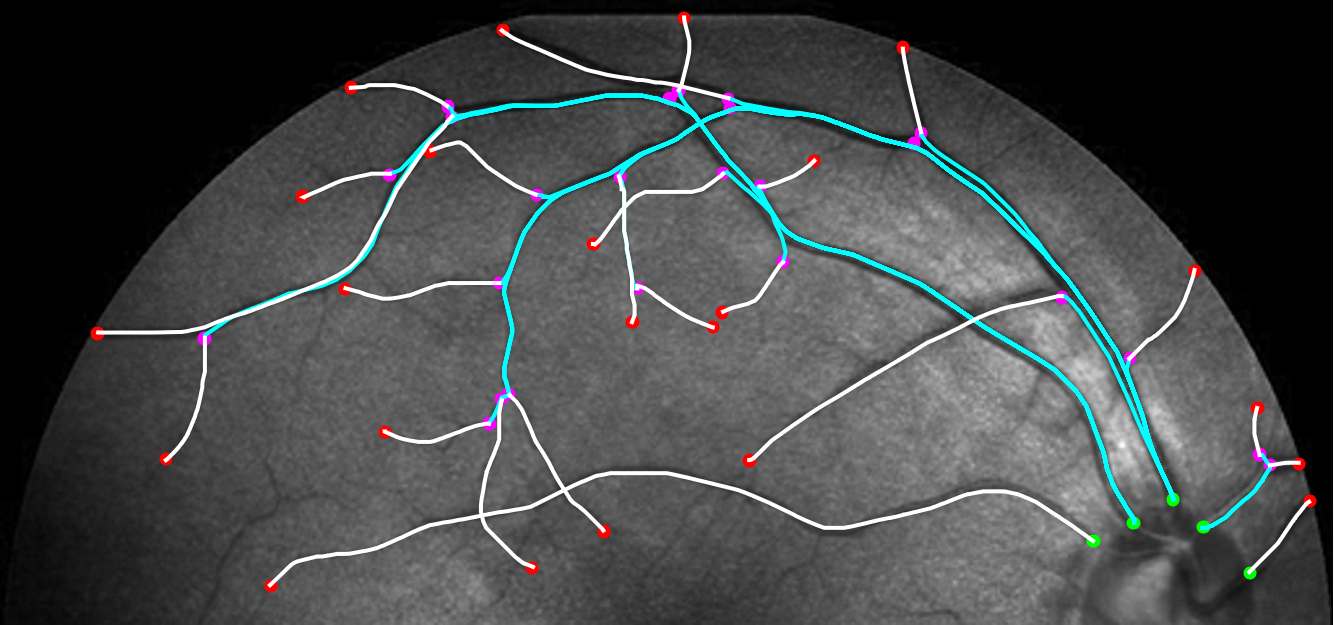}
    \caption{This tracking result (bottom) of a vascular tree in an optical image of the eye (top) is calculated with only two runs of the anisotropic fast marching algorithm. In the images, seeds, bifurcations, and tips are indicated by green, purple, and red points respectively. The white and cyan lines denote the tracking results obtained in the first and second run respectively. Details follow in the experimental section \ref{sec:experiments}.}
    \label{fig:1StepTrackingCompleteVesselTree}
\end{figure}

\begin{figure}
    \centering
    \includegraphics[width=0.4\textwidth]{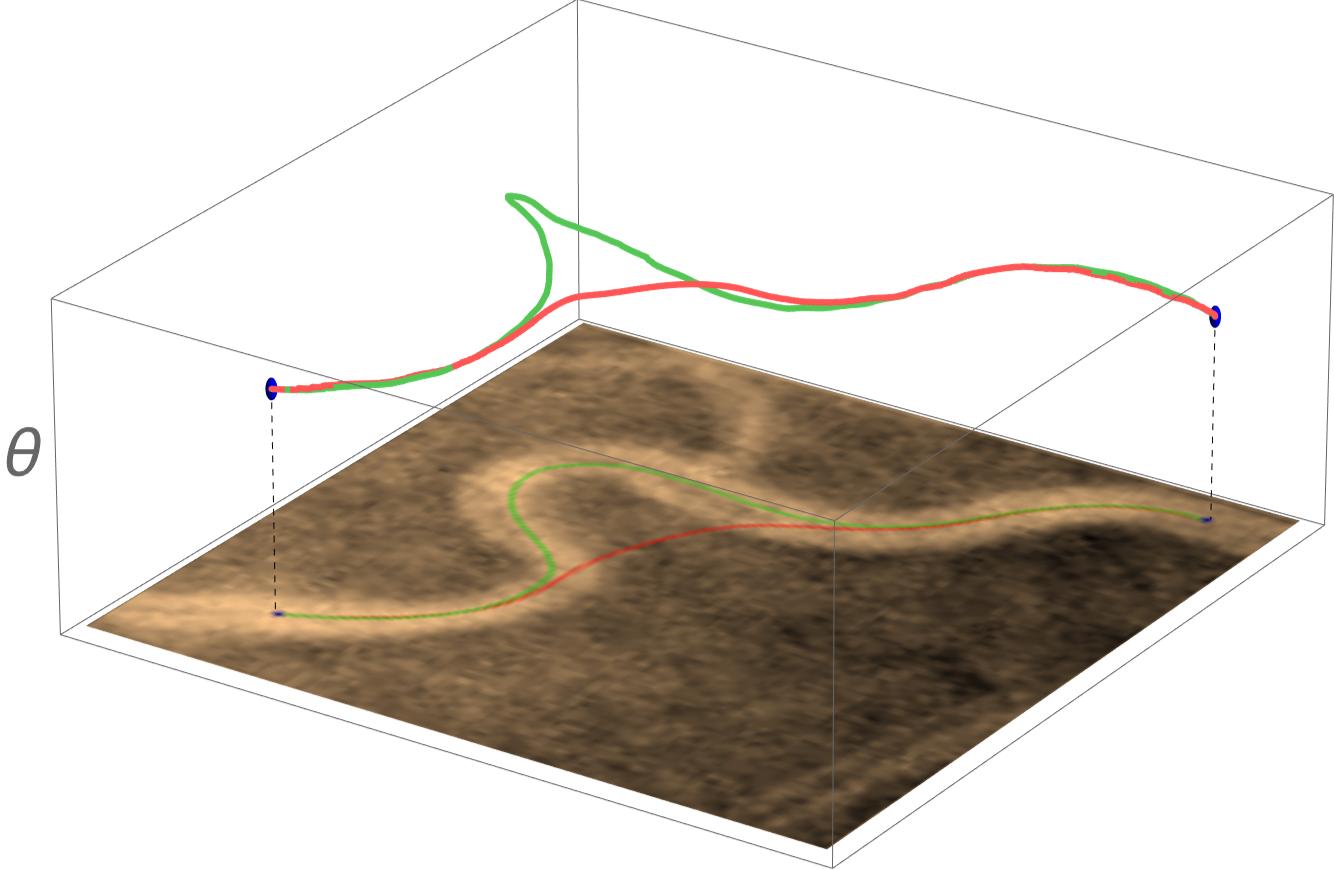}
    \caption{Tracking result of the previous left-invariant model \cite{duits2018optimal} (red), and the new data-driven left-invariant model (green). The tracking performed in the lifted space of positions and orientations is projected back onto the input 2D image. Our proposed model (in green) demonstrates a significantly improved accuracy in adapting to curvature of blood vessels in the optical image.}
    \label{fig:1vesselTracking}
\end{figure}

In tracking an entire vascular tree, we limit the number of anchor points to at most one (which is computed without explicit manual supervision) and only use the boundary conditions for each vessel (see Fig.~\ref{fig:1StepTrackingCompleteVesselTree}). 
Thanks to our new modified version of the anisotropic fast marching algorithm \cite{mirebeau2018fast}, we can now better adapt to curvature and spatial misalignment efficiently (see Fig.~\ref{fig:1vesselTracking}). We also address common pitfalls at complex overlapping structures, where one must impose additional constraints to avoid taking wrong exits in the tracking. The implementation of such constraints is easily accounted for in our model, as we will see.  

Similarly to the plus-control variant of previous geometric curve optimization models 
\cite{duits2018optimal} we replace symmetric, anisotropic, (sub-)Riemannian geodesic tracking by \emph{asymmetric}, anisotropic Finslerian models that avoid cusps and allow for automatic placement of in-place rotations (recall Fig.~\ref{fig:cuspVsInPlaceRotation} and see \cite{duits2018optimal,Mashtakov2021extremal}). 
As a result, our model automatically accounts for bifurcations, thereby reducing the number of anchor points.

In summary, the main contributions of this article are: 
\begin{enumerate}
    \item We introduce a new geodesic tracking model that uses a crossing-preserving approach for tracking complex vasculatures in $\mathbb{M}_2$. 
    Our model uses a new 
    anisotropic fast marching algorithm to compute cusp-free data-driven geodesics.
    The induced geometric vessel tracking better adapts 
    for vessel curvature and orientation sampling biases, compared to the previous model in \cite{duits2018optimal}.
    \item We mathematically analyze these solutions (the family of all geodesics) via our data-driven version of the plus Cartan connection (Section~\ref{ch:model}) that underlies the Hamiltonian flow as we will show in Theorem~\ref{th:th}. 
    \item Finally, we demonstrate our method on highly challenging examples of retinal images with complex vasculature where adequate tracking results are obtained with only two runs of the proposed anisotropic fast marching algorithm. 
\end{enumerate}
\subsubsection*{Structure of the Article}
In Section~\ref{sec:background}, we provide background of the geometrical tools underlying our method. We explain the 
space of positions and orientations $\mathbb{M}_2$, and why it is beneficial to apply tracking in this 3D space rather than in 2D position space. 

In Sections~\ref{ch:model} and \ref{sec:ThShortVsStraight}, we describe our model.
We begin in Section~\ref{ch:model} by introducing a new data-driven Cartan connection $\nabla^U$ associated with a data-driven left-invariant metric tensor field $\mathcal{G}^U$. These geometric tools allow for curvature adaptation and correction of misalignments in existing geodesic tracking algorithms in $\mathbb{M}_2$.

In Section~\ref{sec:ThShortVsStraight} we use the data-driven Cartan connections and data-driven left-invariant metric tensor fields. We present our main theoretical result in Theorem\ref{th:th} where we 
\begin{itemize}
    \item characterize `straight curves' and `shortest curves' in data-driven left-invariant Riemannian manifolds on a finite-dimensional Lie group $G$, 
    \item analyze the 
Hamiltonian flow of all geodesics together, 
\item provide the geodesic backtracking formula of the new geodesic tracking model, 
\item address the symmetries of the geodesics and connections of the new model.
\end{itemize}
Then in Section~\ref{sec:ExtensionAFM} we employ the geometrical models and tools and present a numerical algorithm to compute the distance map for the special case where the Lie group equals the roto-translation group $G=SE(2)$. Additionally, we explain how to compute the backtracking of geodesics from end to source point. 
We present a new version of the anisotropic fast marching algorithm \cite{mirebeau2019hamiltonian} that applies to our new data-driven model. 

In Section~\ref{sec:experiments}, we report an extensive experimental evaluation of geodesic tracking in retinal images from the annotated STAR dataset \cite{Zhang2016,AbbasiSureshjani}, and show that our new model allows for adequate geometric tracking of highly complex vasculatures.

Finally, in Section~\ref{sec:conclusion}, we end with a brief discussion of future work and conclude.




\section{Background}\label{sec:background}
\subsection{Lifted space of positions and orientations $\mathbb{M}_2$}
We begin by introducing the lifted space of positions and orientations $\mathbb{M}_2$. As motivated earlier, working in this space allows for convenient ways to separate difficult crossing structures in $\mathbb{R}^2$. In this article, we specifically focus on the challenging problem of vessel tracking in retinal images having complex vasculature. 
\begin{definition}
    The space of two-dimensional positions and orientations $\mathbb{M}_2$ is defined as a smooth manifold
    \begin{align*}
        \mathbb{M}_2:=\mathbb{R}^2\rtimes S^1,
    \end{align*}
    where $S^1\equiv\mathbb{R}/(2\pi\mathbb{Z})$ using the identification 
    \begin{equation} \label{ID}
    \ul{n}=(\cos \theta, \sin \theta) \leftrightarrow \theta. 
    \end{equation}
    Elements in the homogeneous space are denoted by ordered pairs $(\x,\theta) \in \R^{2} \times S^{1}$ but to stress the semidirect product structure of the roto-translation group $SE(2):=\R^{2} \rtimes SO(2)$ acting on $\mathbb{M}_{2}$ we write $\mathbb{M}_{2}=\mathbb{R}^{2} \rtimes S^{1}$. 
\end{definition}
The space $\mathbb{M}_{2}$ is a homogeneous space under the transitive action of roto-translations given by the following mapping:
\begin{definition}[Lie group action on domain]\newline
    For each roto-translation $\g=(\y,R_{\alpha}) \in SE(2)$ the mapping of $L_\g:\mathbb{M}_2\to\mathbb{M}_2$ is given by
    \begin{equation} \label{straightL}
        L_\g (\x,\theta):=(\y,\alpha).(\x,\theta)=(\y+R_\alpha \x,\alpha+\theta),
    \end{equation} 
    where $R_\alpha\in SO(2)$ is the matrix associated with a counterclockwise rotation with rotation angle $\alpha\in S^1$.
\end{definition}
Clearly, the concatenation of two rigid body motions is again a rigid body motion and indeed one has
\[
L_{\g_2}L_{\g_1}=L_{\g_2 \g_1} \textrm{ with }  \g_2 \g_1:=(\x_2+ R_{\theta_{2}}\x_1, R_{\theta_2+\theta_1}),
\]
for all $\g_1,\g_2\in SE(2)$, i.e. 
$L$
is a group representation.

After setting a reference point $\p_0=(\mathbf{0},0) \in \mathbb{M}_{2}$
we have a 1-to-1 relation between the roto-translation $(\x,R_{\theta})$ mapping $\p_0$ to $\p=(\x,\theta)$ and the point in the homogeneous space $\p=(\x,\theta)\in\mathbb{M}_2$. In particular, $\p_0$ is then identified with the unity element $\e:=(\mathbf{0},I)$\footnote{Since $\text{Stab}_{SE(2)}(\p_0)=\{\e\}$}. In short, we identify \begin{equation*}
    \mathbb{M}_2\equiv SE(2)\text{ via }(\x,\theta)\leftrightarrow \g=(\x,R_\theta).
\end{equation*}
Under this identification, the product of two elements, say \mbox{$(\textbf{x}_1,\theta_1),(\textbf{x}_2,\theta_2)\in\mathbb{M}_2$}, in the space of positions (\mbox{$\textbf{x}_1,\textbf{x}_2\in\mathbb{R}^2$}) and orientations ($\theta_1,\theta_2\in S^1$) is given by
\begin{align*}
    (\textbf{x}_2,\theta_2)\cdot(\textbf{x}_1,\theta_1)=(\textbf{x}_2+R_{\!\theta_2}\textbf{x}_1,\theta_2+\theta_1),
\end{align*}
where \emph{henceforth} $R_{\!\theta_2}\in SO(2)$ denotes the counter-clockwise rotation matrix with rotation angle $\theta_2 \in \R/(2\pi \mathbb{Z})$.

\begin{definition}
Let $X$ be a Banach space, then $B(X)$ denotes the space of bounded linear operators on $X$.
\end{definition}
As mentioned we lift the image from $\mathbb{R}^2$ to $SE(2)$ to separate crossing structures in the corresponding orientation score (Fig.~\ref{fig:OrientationScores}).
Next we explain how this is precisely done.
\begin{definition}[Orientation Scores]\label{def:orientationScores} \\
The orientation score transform $W_\psi:\mathbb{L}_2(\mathbb{R}^2)\to \mathbb{L}_2(SE(2))$
using anisotropic wavelet $\psi$
maps an image $f$ to an orientation score $U=W_{\psi}f$. The orientation score is given by
    \begin{equation} \label{score} 
    \begin{array}{ll}
        (W_\psi f)(\x,\theta )&:=\int \limits_{\mathbb{R}^2}\overline{\psi(R_\theta^{-1}(\y-\x))}\, f(\y)\,\mathrm{d}\y \\
        &= (\mathcal{U}_{\g}\psi,f)_{\mathbb{L}_{2}(\R^2)},
        \end{array}
    \end{equation}
\end{definition}
where the rotated and translated mother wavelets are obtained by the Lie group action $\mathcal{U}:SE(2)\to B(\mathbb{L}_{2}(\R^2))$ given by\[(\mathcal{U}_{\ul{g}}\psi)(\ul{y})=
\psi(R_{\theta}^{-1}(\ul{y}-\ul{x})),
\]
for all $
 \ul{g}=(\ul{x},R_{\theta}) \in SE(2)$ and all $\ul{y}\in \R^2$.

Definition \ref{def:orientationScores} inputs an image $f\in \mathbb{R}^2$ and yields a function $U\in SE(2)$. This is achieved by taking a convolution with a rotated wavelet filter, where the canonical/mother wavelet function $\psi$ is rotated counter-clockwise with the angle $\theta$, as we can see in \eqref{score}. By varying $\theta$ over all orientations $\mathbb{R}/(2\pi\mathbb{Z})$, the image is lifted from $\mathbb{R}^2$ to $\mathbb{M}_2$.
We use the real part of the cake wavelets \cite{duits2005phd,bekkers2014multi} depicted in Fig.~\ref{fig:cakewavelets} as then the space of orientation scores is naturally embedded in $\mathbb{L}_2(SE(2))$ \cite{duits2005phd}, and gives practically informative orientation scores \cite{bekkers2017phd}. More information on the orientation score transform, its range, its invertibility, and the choice of wavelets $\psi$ can be found in previous works \cite{duits2007invertible,duits2005phd,DuitsAMS1,bekkers2017phd}.

We can rotate and translate images via $f \mapsto \mathcal{U}_\g f$. This corresponds to a left action on the orientation score: 
\begin{equation} \label{IT}
W_{\psi} \circ \mathcal{U}_{\g}= \mathcal{L}_{\g} \circ W_{\psi}
\end{equation} 
for all $\g \in SE(2)$.
Left actions are defined as follows:
\begin{definition}[Lie group action on orientation scores] \\
    The left regular representation $\mathcal{L}_\g:SE(2)\to B(\mathbb{L}_2(SE(2)))$ is given by
    \begin{equation} \label{mathcalL}
        \mathcal{L}_\g U(\h):=U(L_{\g^{-1}} \h)=U(\g^{-1}\h)
    \end{equation}
\end{definition}
for all $\g,\h\in SE(2), U\in\mathbb{L}_2(SE(2))$.

As shown in \cite[Thm.~21]{duits2005phd}, by construction of (\ref{score}),  orientation score processing must be left-invariant (i.e., equivariant with respect to left actions $\mathcal{L}_\g$) and not right-invariant. The key reason for this is the fundamental relation (\ref{IT}). This left-invariance will also be crucial in our geometric tracking which needs to be \textbf{left-invariant} (and not right-invariant). Indeed rotating and translating the image should yield to an equally rotated and translated (lifted) tracking output curve. 


\subsection{Metric Tensor Fields and Finsler Functions}
To calculate shortest paths in orientation scores, we need to establish local costs on tangents (velocities). 
We do so by assigning a metric tensor $\mathcal{G}_\p(\cdot,\cdot)$ to every point $\p=(\x,\theta)$ in the lifted space of positions and orientations. It is beneficial to design this metric tensor depending on the specific application. 
Typically, this choice of the metric tensor field establishes the geometric model. Additionally, diagonalizing this tensor at every point $\p$ provides a local frame of reference.

First, we introduce the static frame denoted by $\{\partial_x,\partial_y,\partial_\theta\}$, induced by the coordinates $x,y,\theta$ for all points in $\mathbb{M}_2$. Its dual frame, for the cotangent bundle $T^*(\mathbb{M}_2)$, is denoted by $\{\mathrm{d}x,\mathrm{d}y,\mathrm{d}\theta\}$, and can be used to express the metric tensor field $\mathcal{G}$.

It is advantageous to use left-invariant vector fields for our application, since it guarantees that tracking results are equivariant to the group of roto-translations. More specifically, tracking is independent of the roto-translation of the image, meaning that tracking on a roto-translated image is identical to tracking on the original and roto-translating the result.

\begin{definition}[Frame of Left-Invariant Vector Fields] \label{def:LIVM} \\
    The frame of left-invariant vector fields (left-invariant frame) is obtained by a pushforward of the static frame at the origin $\p_0$. We define the pushforward $(L_\g)_*:T_\h(G)\to T_{\g\h}(G)$ by
    \begin{align*}
        (L_\g)_*\left.\partial_x\right|_{\p_0} U=\left.\partial_x\right|_{\p_0}(U\circ L_\g),
    \end{align*}
    for all smooth functions $U:\mathbb{M}_2\to\mathbb{C}$.
    Then, the left-invariant frame $\{\mathcal{A}_1,\mathcal{A}_2,\mathcal{A}_3\}$ is defined by
    \begin{align*}
        &\left.\mathcal{A}_i\right|_{(x,y,\theta)}=(L_{(x,y,\theta)})_*\left.\mathcal{A}_i\right|_{\p_0} \textrm{with }\\
        &\left.\mathcal{A}_1\right|_{\p_0} =\left.\partial_x\right|_{\p_0},
        \left.\mathcal{A}_2\right|_{\p_0}=\left.\partial_y\right|_{\p_0},
        \left.\mathcal{A}_3\right|_{\p_0} =\left.\partial_\theta \right|_{\p_0}.
    \end{align*}
    After computations,
    we obtain the left-invariant vector fields \[\begin{array}{l}
    \mathcal{A}_1=\cos \theta \, \partial_x+\sin \theta \, \partial_y, \\
    \mathcal{A}_{2}=-\sin \theta \, \partial_x+\cos \theta \, \partial_y \textrm{ and } \mathcal{A}_{3}=\partial_{\theta}.
    \end{array}
    \]
    The corresponding dual frame is given by $\{\omega^i\}_{i=1}^3$ where $\omega^i\!\left(\mathcal{A}_j\right)=\delta_j^i$. A brief computation gives
    \[
    \begin{array}{l}
    \omega^1\!= \cos \theta \,{\rm d}x +\sin \theta\, {\rm d}y, \\
    \omega^2\!= -\sin \theta \,{\rm d}x +\cos \theta \,{\rm d}y \textrm{ and } \omega^3={\rm d}\theta.
    \end{array}
    \]
\end{definition}
In addition to that, we define the left-invariant metric tensor field:
\begin{definition}
Metric tensor field $\mathcal{G}$ on $\mathbb{M}_{2}$ is left-invariant iff 
\[
\mathcal{G}_{\g\cdot\p}((L_\g)_* \dot{\p}, (L_\g)_* \dot{\p})= 
\mathcal{G}_{\p}(\dot{\p}, \dot{\p})
\]
for all $\p \in \mathbb{M}_{2}$, all $\dot{\p} \in T_{\p}(\mathbb{M}_{2})$ and all  $\g \in SE(2)$.
\end{definition}

\begin{remark}[Left-invariant Metric Tensor Field] \\
    Let $\mathcal{G}$ denote a left-invariant metric tensor field on $G$. Then there exists a unique \textit{constant} matrix $[g_{ij}]_{ij}\in\mathbb{R}^{3\times 3}$ such that
    \begin{align}
        \mathcal{G}=\sum_{i,j=1}^3 g_{ij} \left.\omega^i\right.\otimes\left.\omega^j\right.,
        \label{eq:LIMTF}
    \end{align} 
    where $\otimes$ denotes the usual tensor product.
\end{remark}
In the standard left-invariant model we restrict ourselves to the case $g_{ij}=
g_{ii}\delta_{ij}$, and then $\mathcal{G}$ is diagonal with respect to the co-frame $\{\omega^i\}_i$, i.e. $\mathcal{G}=\sum \limits_{i=1}^3g_{ii}\;\omega^i\otimes\omega^i$. 

Often we do not want to work with symmetric Riemannian metric tensor fields (for instance to avoid cusps, and to ensure that fronts only move forward, see Fig.~\ref{fig:TortuousVesselTracking}), and then we resort to the general Finsler geometry as done in \cite{duits2018optimal} and \cite{mirebeau2018fast}. 
Essentially, this means that we replace the symmetric norm $\sqrt{\left. \mathcal{G} \right|_{\gamma(t)}(\dot{\gamma}(t),\dot{\gamma}(t))}$ 
in the Riemannian distance/metric:
\begin{equation}\label{RD}
d_{\mathcal{G}}(\p,\q)=\underset{{\scriptsize\begin{array}{c}\gamma\in \Gamma_1,\\
\gamma(0)=\p,\gamma(1)=\q
\end{array}}}{\inf} \int_0^1\sqrt{\mathcal{G}_{\gamma(t)}\left(\dot{\gamma}(t),\dot{\gamma}(t)\right)}\, \mathrm{d}t
                \end{equation}
by an asymmetric Finsler norm $\mathcal{F}(\gamma(t),\dot{\gamma}(t))$, given by
\begin{equation} \label{Finsler}
\begin{array}{l}
\mathcal{F}^2(\p,\dot{\p})= \mathcal{G}_\p(\dot{\p},\dot{\p}) 
+g_{11}(\varepsilon^{-2}-1) 
\min\hspace{-0.5ex}\left\{0,\omega^1_\p(\dot \p)\right\}^2
\end{array}
\end{equation}
where the relaxation parameter
$0<\varepsilon \ll 1$ punishes spatial backward motions.
For  convergence results of geodesics and distances if $\varepsilon \downarrow 0$, see \cite{duits2018optimal,Franceschiello2019geometrical}.

Later in subsection~\ref{sec:geodesictracking} we will explain all parameter settings including the choice of $g_{11},g_{22},g_{33}>0$ in $\mathcal{F}$.

Next, we define some relevant geometrical tools associated to \eqref{RD} and \eqref{Finsler}.
\begin{remark}[The Space of Curves over which we optimize] \\
To adhere to standard conventions in Riemannian geometry we optimize over the space of piecewise continuously differentiable curves in $\mathbb{M}_2$ (indexed by $T>0$): 
\begin{equation} \label{SpaceCurves}
\Gamma_T:=PC^{1}([0,T],\mathbb{M}_{2}).
\end{equation}
In (\ref{RD}) we set $T=1$, as there the choice of $T$ is irrelevant by parameterization independence of the functional.
\end{remark}

\begin{remark}[Control Sets] \\
    The control set in the tangent bundle $T(\mathbb{M}_2)$ is defined as
    \begin{align*}
        \mathcal{B}_{\mathcal{G}}(\p):=\left\{\dot{\p}\in T_\p(\mathbb{M}_2)\left|\sqrt{\mathcal{G}_\p\left(\dot{\p},\dot{\p}\right)}\leq 1\right.\right\},
    \end{align*}
    with $\p\in\mathbb{M}_2$ and $\mathcal{G}$ the underlying Riemannian metric tensor field.
    
    The corresponding asymmetric control set in the tangent bundle $T(\mathbb{M}_2)$ is defined as
    \begin{equation} \label{CSet}
        \mathcal{B}_{\mathcal{F}}(\p):=\left\{\dot{\p}\in T_\p(\mathbb{M}_2)\right.\; \left|\mathcal{F}(\p,\dot{\p})
        \leq 1 \right\},
    \end{equation}
    with $\p=(\x,\n)\in\mathbb{M}_2$ 
    and $\mathcal{F}$ the underlying Finslerian metric tensor field. In the limiting case where backward motions become prohibited as $\varepsilon \downarrow 0$ (i.e. the sub-Finslerian setting) we only get half of the Riemannian control sets 
    \begin{equation}
    \mathcal{B}_{\mathcal{F}}(\p)= \{ \dot{\p}=(\dot{\x},\dot{\n}) \in \mathcal{B}_{\mathcal{G}}(\p)|\; \dot{\mathbf{x}} \cdot \mathbf{n} \geq 0\}.
    \label{eq:RelationControlsetGandF}
    \end{equation}
\end{remark}

\subsection{Cartan Connections}\label{sec:CartanConnections}
The theory of Cartan connections was developed by \'Elie Cartan. His viewpoint on differential geometry relies on moving frames of reference (rep\`ere mobile). The idea is to connect tangent spaces by group actions on homogeneous spaces. This geometric tool allows us to understand the geodesic flow associated to \eqref{RD} and its data-driven extensions.

The homogeneous space that we will use for crossing-preserving 2D image processing is the homogeneous space of positions and orientations $\mathbb{M}_2$. 
Here, the pushforward $(L_\g)_*$ of the left-multiplication \textit{connects} $T_\e(\mathbb{M}_2)$ to $T_\g(\mathbb{M}_2)$ as it maps $T_\e(\mathbb{M}_2)$ (isometrically\footnote{W.r.t. norm induced by the metric tensor field $\mathcal{G}$ in Definition~\ref{def:LIVM}.}) onto $T_\g(\mathbb{M}_2)$ and $(L_\g^{-1})_*$, known as the Cartan-Ehresmann form, maps $T_\g(G)$ back to $T_\e(G)$.

First, we introduce the general definition of Cartan connections, after which we also introduce the Cartan plus connection \cite{duits2021Springer}. In this article, we will introduce a data-driven version of the Cartan plus connection, leading to a generalization of the existing theory on shortest and straight curves in $\mathbb{M}_2$.
\begin{definition}[Cartan Connection \cite{duits2021Springer}]\newline
    A Cartan connection on a Lie group $G$ is a tangent bundle connection with the following additional properties
    \begin{enumerate}
        \item left invariance: if $X,Y$ are left-invariant vector field then $\nabla_X Y$ is a left-invariant vector field,
        \item for any $A\in T_\e(G)$ the exponential curve is auto-parallel, i.e.
        $\nabla_{\dot{\gamma}(t)}\dot{\gamma}(t)=0$ where $\gamma(t)=\gamma(0)\, \exp(tA)$.
    \end{enumerate}
\end{definition}
We use the following special case of a Cartan connection to define shortest and straight curves. Note that this Cartan connection is easily expressed in the left-invariant frame \cite{KobayashiNomizu,Cogliati2017CartanSA,Piuze2015MaurerCartan,Momayyez2009Stochastic,PhDDaChen,duits2021Springer}.

\begin{definition}[Cartan Plus Connection \cite{duits2021Springer}]\newline
    Consider a Lie group $G$ of finite dimension $n$, with Lie brackets $[\cdot,\cdot]$ and structure constants $c_{ij}^k\in\mathbb{R}$ s.t. \[
    [\mathcal{A}_i,\mathcal{A}_j]=\sum \limits_{k=1}^n c_{ij}^k\mathcal{A}_k.
    \]
    Then the Cartan plus connection is given by
    \begin{align} \label{form}
        \nabla^{[+]}:=\sum_{k=1}^n\left(\sum_{i=1}^n\omega^i\otimes\left(\mathcal{A}_i\circ\omega^k\right)+\sum_{i,j=1}^n\omega^i\otimes\omega^j c_{ij}^k\right)\mathcal{A}_k.
    \end{align}
\end{definition}
\begin{remark}
    Note that the $\circ$ symbol denotes the composition of functions such that for example
    \begin{align*}
        &\mathcal{A}_i\circ\omega^k\left(\sum_{l=1}^n \alpha^l\mathcal{A}_l\right)=\mathcal{A}_i\left(\omega^k\left(\sum_{l=1}^n \alpha^l\mathcal{A}_l\right)\right)
        =\mathcal{A}_i\left(\alpha^k\right).
    \end{align*}
\end{remark}
For explicit coordinate expressions, see \cite{duits2021Springer}.
Next, we explain how to read and compute \eqref{form}. The covariant derivative $\nabla_X Y$ of a vector field $Y$ with respect to a vector field $X$ is again a vector field. Indeed the above formula gives
{\small\begin{equation}
    \nabla^{[+]}_X Y=\!\sum_{k=1}^n\left(\sum_{i=1}^n\omega^i(X)\left(\mathcal{A}_i\circ\omega^k(Y)\right)+\!\!\sum_{i,j=1}^n\!\!\omega^i(X)\omega^j(Y) c_{ij}^k\right)\!\mathcal{A}_k,\label{eq:formCoordinates}
\end{equation}}so that it becomes clear where $X$ and $Y$ typically enter in the open slots of the expression (\ref{form}). Note that vector field $\mathcal{A}_i$ in \eqref{eq:formCoordinates} is a differential operator applied to the smooth function $G\ni\g\mapsto \omega^k_\g(Y_\g)\in\mathbb{R}$.
\begin{remark}
The connection $\nabla^{[+]}$ is called `Cartan plus connection' as
we add the two sums between the two large round brackets. In differential geometry one also has Cartan connections with a realvalued scalar factor in front of the second term, but this does not serve our applications \cite{Duitsbookchapter}.
\end{remark}
Now that we explained Cartan connections, let us return to the core purpose of designing a geometric model such that projected geodesics follow the blood vessels.

\subsection{Geodesic Tracking in the Space of Positions and Orientations $\mathbb{M}_2$}\label{sec:geodesictracking}
Geodesic tracking methods aim to find the shortest paths following the underlying biological blood vessels in the retinal image. Such shortest paths are obtained by finding minimizing geodesics, which are defined as curves with the shortest length functionals. Typically, such length functionals are driven by a cost function that is small at locations of the blood vessels and high at all other places. Many different approaches to determine the minimizing geodesics have been proposed over the years, ranging from classical geodesic tracking in the image domain \cite{Kimmel1997GeodesicActiveControus} to tracking in higher-dimensional homogeneous spaces \cite{bekkers2017phd,bekkers2015pde,benmansour2011tubular,duits2018optimal}.

As already explained, we lift the input image $f$ from $\mathbb{R}^2$ to $\mathbb{M}_2$ using orientation scores in homogeneous spaces \cite{duits2007invertible,DuitsAMS1} (see Fig.~\ref{fig:OrientationScores}). 
The tracking process involves computing a geodesic distance map in $\mathbb{M}_2$ and then using steepest descent to find the \textit{shortest curve} in the lifted space. Finally, we project the curve back onto the input image in $\mathbb{R}^2$ to get the final tracking result, see the examples in Figure~\ref{fig:1vesselTracking}.
Over time, different models have been introduced that describe how the geodesics should behave. 
Imagine a car moving along such geodesics. Then the Reeds-Shepp car model \cite{reeds1990optimal} which describes the problem of shortest paths for cars between an initial and final point, and the Reeds-Shepp forward model \cite{duits2018optimal} turns off the reverse gear of the car. In both cases the spatially projected geodesics (optimal paths) tend to follow blood vessels in medical images well.

\subsubsection{Symmetric Reeds-Shepp Car Model}
The left-invariant metric tensor field of the symmetric Reeds-Shepp car model, $\mathcal{G}$, is given by the symmetric tensor field 
\begin{equation} \label{eq:LIMTFDD}
\begin{array}{c} 
    \mathcal{G}=C^2\left(\xi^2\omega^1\otimes\omega^1+\!\frac{\xi^2}{\zeta^2}\omega^2\otimes\omega^2+\!\omega^3\otimes\omega^3\right)   \\\Leftrightarrow\\
    \mathcal{G}_\p(\dot{\p},\dot{\p})=C(\p)^2\left(\xi^2\left|\dot{\x}\cdot\n\right|^2+\frac{\xi^2}{\zeta^2}\left\|\dot{\x}\wedge\n\right\|^2+\left\|\dot{\n}\right\|^2\right).
\end{array}
\end{equation}
for all
$\ul{p}=(\ul{x},\ul{n}),\dot{\ul{p}}=(\dot{\ul{x}},\dot{\ul{n}})$ with \mbox{$
\|\dot{\x}\wedge\n\|^2:=\|\dot{\x}\|^2-|\dot{\x}\cdot\n|^2.
$}
The anisotropy parameter $\zeta$ penalizes vectors with large sideways components.
Note that the classical sub-Riemannian model corresponds to the limit $\zeta\downarrow 0$. For formal convergence results of the Riemannian model to the sub-Riemannian model see \cite[Thm.2]{duits2018optimal}. In practice, choosing $\zeta=0.1$ usually provides a good enough approximation of the sub-Riemannian model for our purposes. The last parameter $\xi$, a weighting parameter, influences the flexibility of the tracking. It either stimulates or discourages angular movement over spatial movement \cite{duits2018optimal}.

The cost function $C:\mathbb{M}_2\to[\delta,1]$, with $\delta>0$, discourages movement at specific locations, e.g. outside vessel structures. 
In this article, the smooth costfunction $(x,y,\theta)\mapsto C(x,y,\theta)$ is typically a version of the multi-scale crossing preserving vesselness map \cite{hannink2014crossing} explained in Appendix~\ref{app:CostFunction}.
For an impression of what such a map $C$ looks like see the 3D visualization in Figure~\ref{fig:3DcalculatedCostFunction}. 
As we consider rather complex vasculatures it is often more intuitive to display their minimum projections over $\theta$, see for example Figure~\ref{fig:CostFunctionBifSeed}.

\subsubsection{Asymmetric Reeds-Shepp Car Model}\label{sec:ReedsSheppForward}
Besides the symmetric version of the left-invariant metric tensor field of the Reeds-Shepp car model, an asymmetric version has been introduced 
in \cite{duits2018optimal}. The forward gear left-invariant metric tensor field of this model is given by the asymmetric Finsler norm/function 
\begin{equation} \label{ARSCM}
\!
\begin{array}{l}
    |\mathcal{F}(\p,\dot{\p})|^2=\mathcal{G}_{\p}(\dot{\p},\dot{\p})+C(\p)^2  \left(\varepsilon^{-2}\!-\!1\right)\xi^2\left|(\dot{\x}\cdot\n)_{-}\right|^2
    \end{array}
\end{equation}
for all
$\ul{p}=(\ul{x},\ul{n}),\dot{\ul{p}}=(\dot{\ul{x}},\dot{\ul{n}})$, with $a_-:=\min\{0,a\}$. Eq.~(\ref{ARSCM}) coincides with (\ref{Finsler}) with $g_{11}=\xi^2$, $g_{22}= \xi^2/\zeta^2$, $g_{33}=1$.

The parameters $\zeta$ and $\xi$ and the cost function $C$ have the same meaning as in the symmetric model. However, we consider an extra variable $\varepsilon\in(0,1]$ in the asymmetric Reeds-Shepp car model. This parameter determines how strongly the model needs to adhere to the forward gear. 
Note that when $\varepsilon=1$, we find the symmetric Reeds-Shepp car model, and when $\varepsilon\to0$, backward movement becomes prohibited. In that case, we move from cusps to change orientation to in-place rotations visualized (cf. Fig. \ref{fig:cuspVsInPlaceRotation} in Appendix \ref{app:CostFunction}). These asymmetric Finslerian models are also highly beneficial in image segmentation as shown by Chen and Cohen \cite{chen2018fast}.

\subsection{Anisotropic Fast Marching}\label{sec:AFM}
\begin{figure*}
    \centering
    \includegraphics[width=0.44\textwidth]{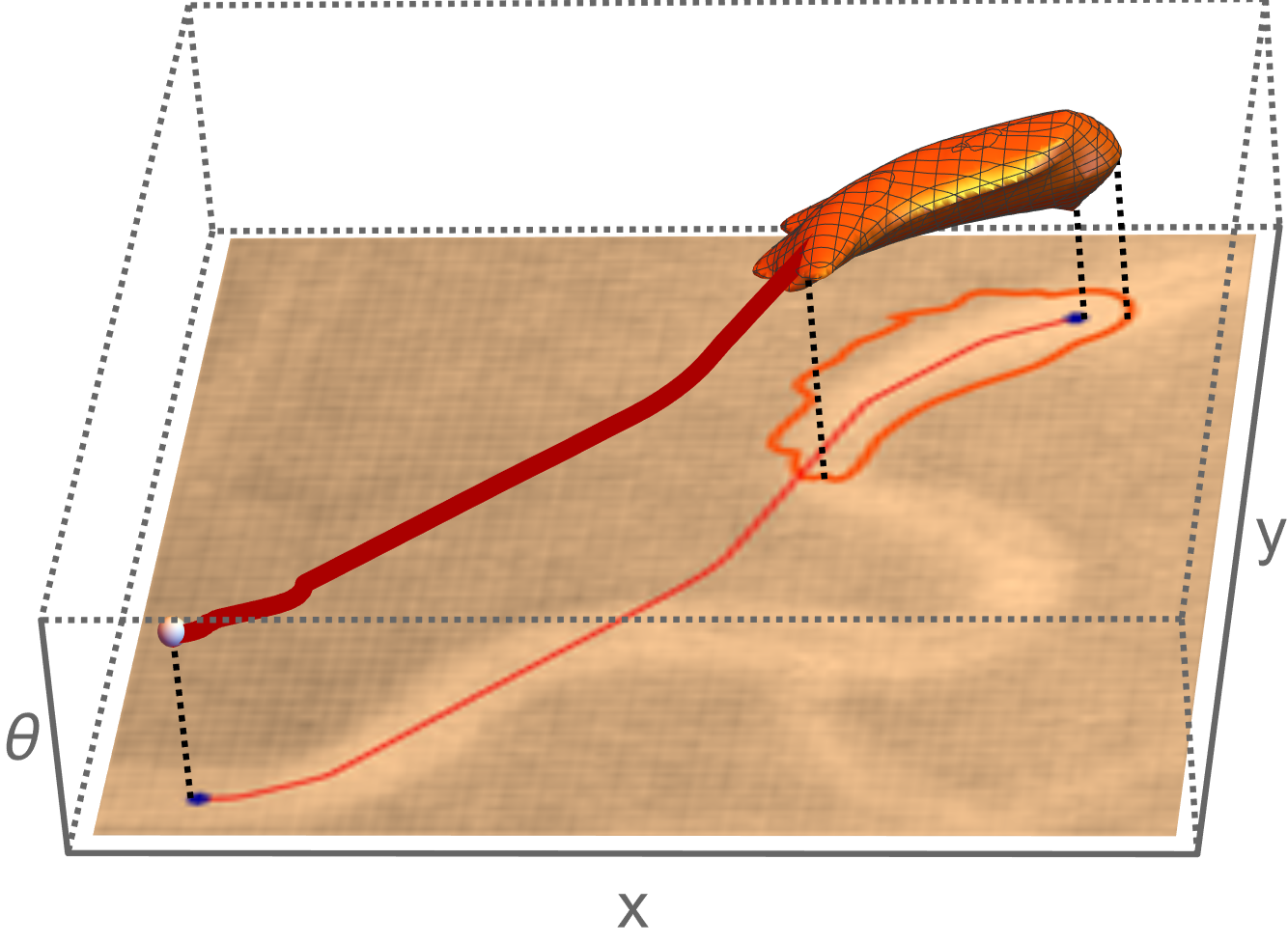}
    \hfill
    \includegraphics[width=0.44\textwidth]{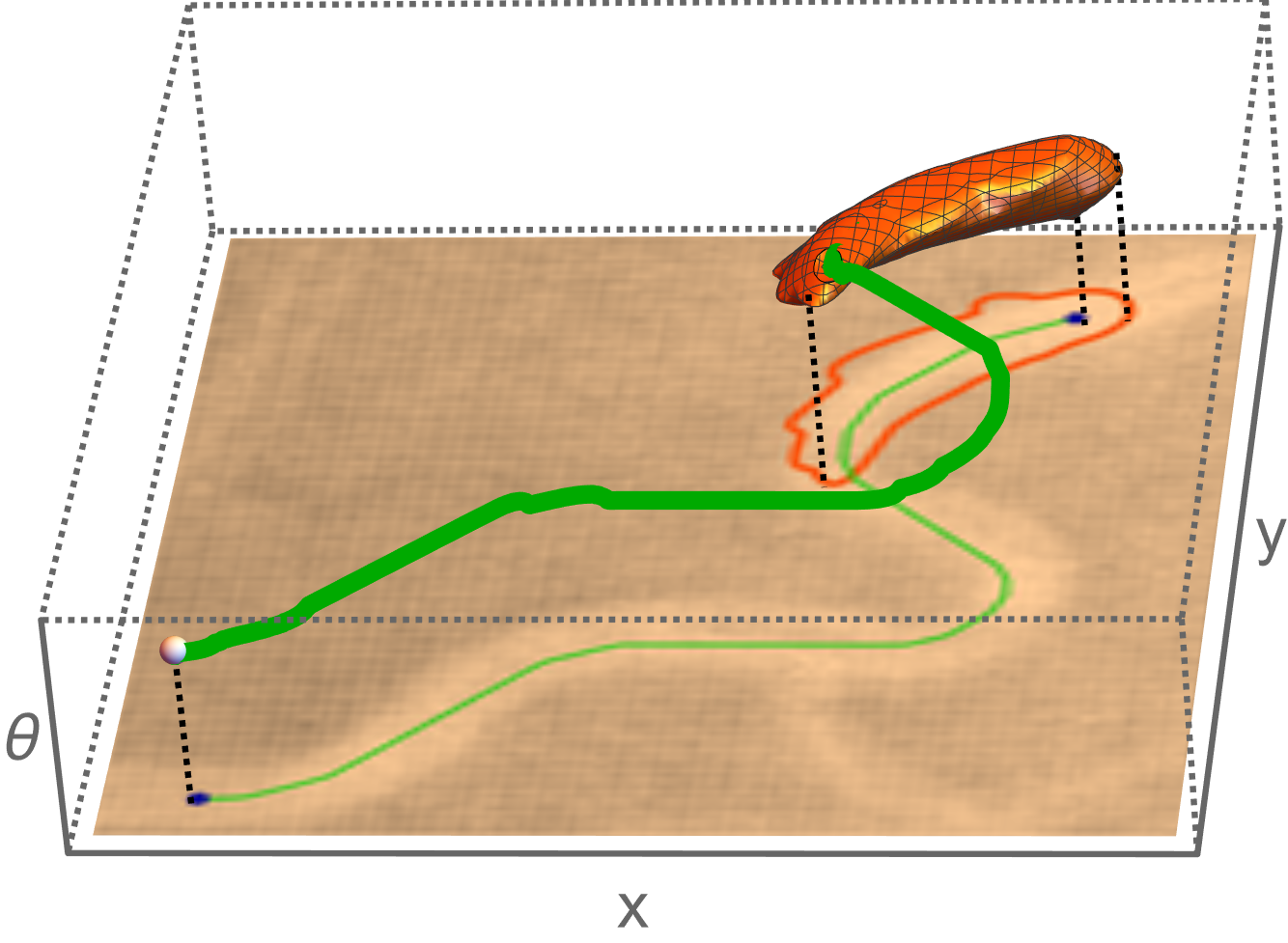}
    \includegraphics[width=0.44\textwidth]{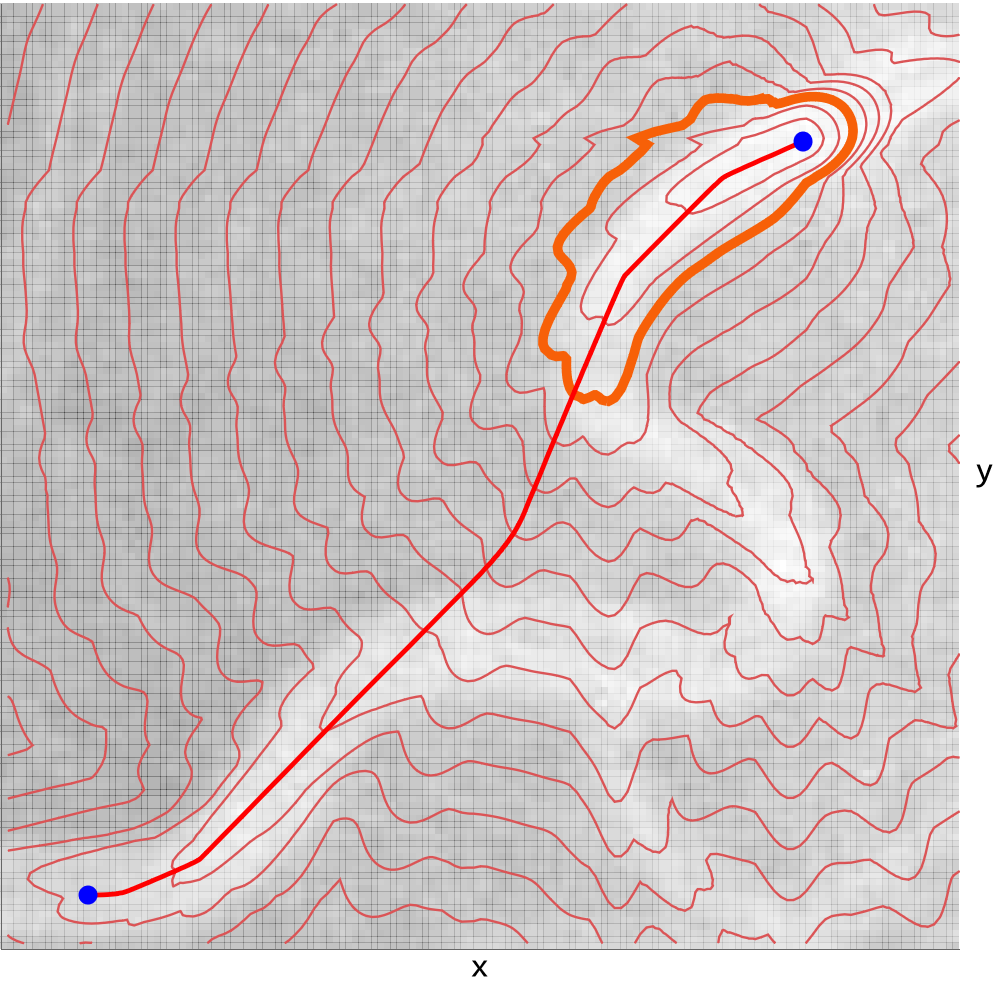}
    \hfill
    \includegraphics[width=0.44\textwidth]{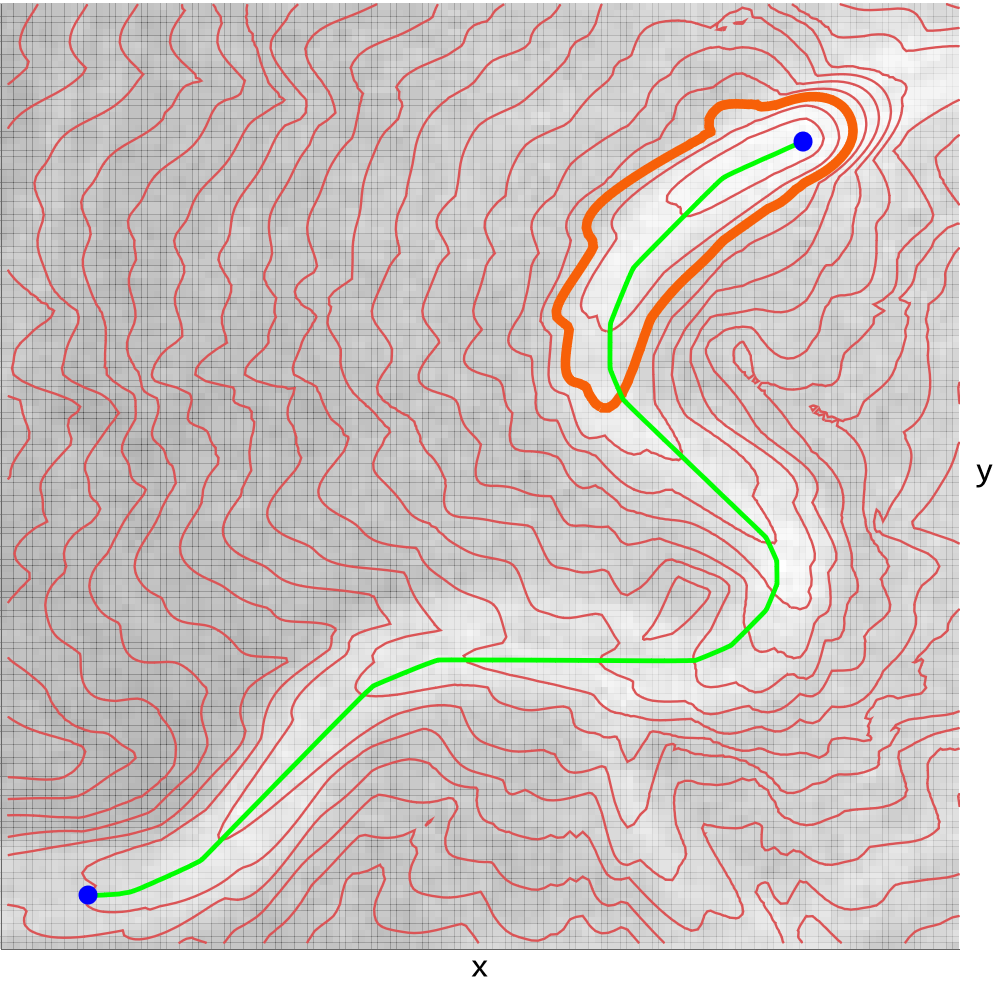}
    \caption{Top row: All points in the orange surface have the same distance to the seed. The isocontours are projected back onto the image, as depicted in the ground plane.
    Bottom row: Several isocontours are projected onto the image and a projection of the curve is visualized.
    Left: The Riemannian geodesic with parallel velocity to the Cartan plus connection $\nabla^{[+]}$ (red) takes a wrong shortcut.
    Right: The Riemannian geodesic (green) with parallel velocity to the Cartan plus connection $\nabla^{U}$ in the Riemannian manifold 
    $(\mathbb{M}_{2},\mathcal{G}^U)$ --or more precisely the Finslerian manifold $(\mathbb{M}_{2},\mathcal{F}^U)$ given by (\ref{Finslerfin})-- 
     does not, and moreover adapts for curvature in $\mathbb{M}_{2}$, cf.~Fig.~\ref{fig:ShortVsStraightCurves}). Explicit formulas for $\nabla^U$ and $\nabla^{[+]}$ will follow later in Table~\ref{tab:ToN}.}
    \label{fig:TortuousVesselTracking}
\end{figure*}

We provide here a brief overview of the partial differential equation (PDE) framework associated with geodesic distance maps, and of their numerical computation, see Section~\ref{sec:ExtensionAFM} for further details. 
We already mentioned that it is common to calculate minimizing geodesics in two steps; first calculating the geodesic distance map, then calculating the shortest curve using steepest descent. To get a first impression of how this looks like in practice, see Fig. \ref{fig:TortuousVesselTracking}. 
The geodesic distance map is characterized, in the PDE framework, as the viscosity solution of a static first-order Hamilton-Jacobi-Bellman equation, known as the Eikonal Equation.
For numerically solving the Eikonal equation, it is discretized using for instance finite differences \cite{bekkers2015pde}, leading to a coupled non-linear system of equations, which is typically solved using a front propagation method such as the fast marching algorithm (FMM).
Classical references on the FMM include \cite{peyre2010geodesic,sethian1999fast}, anisotropic variants are presented in \cite{mirebeau2019hamiltonian,mirebeau2018fast,mirebeau2019riemannian}, and 
the details related to our new model will follow in Section~\ref{subsec:ExtensionAFM_PDE}.

The fast marching algorithm is a numerical method for solving the coupled system of equations discretizing the eikonal PDE. The algorithm proceeds in only one pass over the domain hence providing significant efficiency gains, but also requiring that the numerical scheme satisfies two conditions (monotonicity and causality), see \cite[Def. 2.1]{mirebeau2019hamiltonian,mirebeau2018fast}. The proposed variant of this method uses Selling's algorithm \cite{Selling:1874Algorithm} to calculate in a preliminary step a decomposition of the quadratic forms defining the dual metric. This dual metric suitably only involves positive weights and vectors with integer coordinates, see Proposition~\ref{prop:decomp_D_v}. These ingredients are used to devise an adaptive finite differences scheme, discretizing the anisotropic Eikonal PDE and obeying the required conditions, see Section~\ref{subsec:ExtensionAFM_scheme}. The eikonal PDE is solved via an anisotropic Fast Marching algorithm, and its solution provides the desired distance map. Finally, the minimizing geodesic is calculated by solving an ordinary differential equation defined in terms of the distance map \cite{mirebeau2019hamiltonian,mirebeau2018fast}.

In previous studies of the Reeds-Shepp model and variants \cite{duits2018optimal,mirebeau2018fast}, the geodesic metric tensor matrix featured a block diagonal structure, which was exploited in the discretization. 
However, while working with data-driven metric tensor fields, this block format does not apply! Therefore we adapt the anisotropic fast marching algorithm to cope with the general setting. In this article, we will briefly discuss the changes that were necessary to solve data-driven metric tensor fields. Such data-driven geometric models (that we explain in the next section) give better tracking results than the previous model, as one can see in Fig. \ref{fig:TortuousVesselTracking}. In addition to that, using the anisotropic fast marching algorithm to calculate the geodesics, only a limited number of runs (one for a single vessel, Fig.~\ref{fig:TortuousVesselTracking}, and only two for a full vasculature, Fig.~\ref{fig:1StepTrackingCompleteVesselTree}) are needed to correctly track the vascular structures.

\subsection{Flowchart and Overview of the Methodology}
Before we dive into the details of our method, we provide a sequential flowchart of our methodology, and more information on where more details of each part will be addressed.
\begin{enumerate}
    \item Create an orientation score of input image $f$ (Eq.~\eqref{score});
    \item Calculate the Hessian (App.~\ref{app:GisDDLIV});
    \item Extract the Data-Driven frame from the Hessian (Eq.~\eqref{Finslerfin});
    \item Determine the local cost function for tracking: Vesselness Map (App.~\ref{app:CostFunction});
    \item Identify the Finsler Function (with +-control in App.~\ref{app:AdaptationAsymmetricDataDrivenFinslerFunctions}) with the special Riemannian case (Eq.~\eqref{Finsler} with $\varepsilon=1$) in Thm.~\ref{th:th};
    \item Identify the Dual Finsler Function (with +-control in Lemma~\ref{lem:asym_quad_dual} and App.~\ref{app:AdaptationAsymmetricDataDrivenFinslerFunctions});
    \item Analyze the Hamiltonian Flow of all geodesics (Thm.~\ref{th:th});
    \item Determine the Eikonal PDE distance map (symmetric case: Eq.~\eqref{eq:EikonalPDE}, with +-control: Eq.~\eqref{eq:scheme_consistent_0});
    \item Numerically solve the Eikonal PDE using Anisotropic Fast Marching:
    \begin{enumerate}
        \item Calculate the stencils using Selling matrix decomposition assuring causality for single parts (Prop.~\ref{prop:decomp_D_v});
        \item Follow the procedure of Far-Trial-Accepted points (Sec.~\ref{subsec:fmm});
    \end{enumerate}
    \item Apply steepest descent on the distance map (Sec.~\ref{sec:SteepestDescentFinslerianGeodesics});
    \item Project the geodesic spatially by $$\Pi(x(t),y(t),\theta(t))=(x(t),y(t)).$$
\end{enumerate}

\section{Data-Driven Metric \& Data-Driven Cartan Connection}\label{ch:model}
In multi-orientation image processing, it is beneficial (for vessel segmentation \cite{Zhang2016}) to rely on locally adaptive frames \cite{SavadjievPNAS2012,duits2016locally}. However, the locally adaptive frames in $\mathbb{M}_2$ typically require a stable selection of the principal eigenvector (eigenvector corresponding to the largest eigenvalue) of a symmetrized Hessian of the function $U:SE(2)\to\mathbb{R}$. Recall that a Hessian is defined by a (dual) connection. Even if one uses the left Cartan connection, selecting the principal eigenvector can be locally unstable \cite{duits2016locally} and the largest eigenvalue may not be unique. For instance, if line structures are not locally present at all. Therefore, in this article, we take a slightly different approach by creating an unconditionally stable data-driven left-invariant metric tensor field.


\begin{definition}[Data-Driven~Left-Invariant~Metric] \label{def:10} \\
Let $G$ be a Lie group. Then the metric tensor field $\mathcal{G}^U$ is data-driven left invariant when it satisfies for all $(\g,\dot{\g})\in T(G)$ and all $\q\in G$:
    \begin{equation}
        \mathcal{G}_\g^U(\dot{\g},\dot{\g})=\mathcal{G}_{\q\g}^{\mathcal{L}_\q U}((L_\q)_* \dot{\g},(L_\q)_* \dot{\g}).\label{eq:DDLIF}
    \end{equation}
  \end{definition}

Recall that in our case of interest where $G=SE(2)$ and where $U=\mathcal{W}_{\psi}f$ is an orientation score of the image $f$, the equivariance relation (\ref{IT}) holds, so roto-translation of an image $f \mapsto \mathcal{U}_\g f$ is equivalent to roto-translation $U \mapsto \mathcal{L}_\g U$ of the score.

  Consequently (as will follow in Theorem~\ref{th:th}) if a metric tensor field is data-driven left invariant then a roto-translation $\mathcal{U}_\g f$ of the input image $f$ yields a new geodesic $\gamma_{new}$ that is rotated and translated accordingly:  $\gamma_{new}(\cdot)=\g \gamma(\cdot)$. 
  
    Thus, Definition~\ref{def:10} is a valid constraint in our application as we want the vessel tracking along geodesics to be equivariant with respect to roto-translations.

By creating such a data-driven metric tensor field $\mathcal{G}^U$ on our Lie group of interest $G=SE(2)\equiv \mathbb{M}_{2}$, data-driven corrections are made for spatial and angular misalignment in existing models relying on the standard left-invariant frame \cite{franken2009crossing,duits2010left,duits2016locally}. We will see that a better fitted metric tensor field $\mathcal{G}^U$ has a significant impact on the tracking results for very tortuous vessels, as shown in Figure~\ref{fig:TortuousVesselTracking}. 
For our case of interest $\mathbb{M}_2$, a reasonable choice that satisfies the constraint, and that we use in our experiments, is given by: 
\begin{equation}
\label{Finslerfin}
\boxed{
\begin{array}{l} 
    |\mathcal{F}^U(\p,\dot{\p})|^2
    =|\mathcal{F}(\p,\dot{\p})|^2\!+\lambda \, C^2(\ul{p})\frac{\left\|\left. HU\right|_\p (\dot{\p},\cdot)\right\|_*^2}{\max\limits_{\|\dot{\q}\|=1}\left\|\left. HU\right|_\p (\dot{\q},\cdot)\right\|_*^2} 
    \\[9pt]
     \mathcal{G}^U_{\p}(\dot{\p},\dot{\p})
    =\mathcal{G}_{\p}(\dot{\p},\dot{\p})+\lambda \; C^2(\ul{p})\frac{\left\|\left. HU\right|_\p (\dot{\p},\cdot)\right\|_*^2}{\max\limits_{\|\dot{\q}\|=1}\left\|\left. HU\right|_\p (\dot{\q},\cdot)\right\|_*^2}, \\[9pt]
    \textrm{where  }\mathcal{G} \textrm{ and }\mathcal{F} \textrm{are given in \eqref{eq:LIMTFDD} and \eqref{ARSCM}}, 
\end{array}}
\end{equation}
with $\ul{p}=(\ul{x},\ul{n}) \in \mathbb{M}_{2}$ and $\dot{\ul{p}}=(\dot{\ul{x}},\dot{\ul{n}}) \in T_\p(\mathbb{M}_{2})$. \\[6pt]
Here 
the Hessian field $HU$ is defined in Lemma \ref{def:Hessian} in Appendix \ref{app:GisDDLIV}, and $\|\cdot\|_*$ the dual norm corresponding to the primal norm given by $\sqrt{\mathcal{G}_\p(\dot{\p},\dot{\p})}$ with $\zeta=1$ in \eqref{eq:LIMTFDD}. 


Parameter $\lambda>0$ regulates inclusion of data-driven 2nd order line-adaptation to the orientation score data $U$, cf.~Fig.~\ref{fig:OrientationScores}. 

Finally, the data-driven left-invariant metric tensor field relies on the usual Reeds-Shepp car models $\mathcal{G}$ respectively $\mathcal{F}$ with external smooth cost $C(\ul{p})$ satisfying:
\begin{equation}\label{eq:boundariesCost}
    0<\delta\leq C\leq 1,
\end{equation}
computed from the orientation score $U$, 
as explained in Appendix~\ref{app:CostFunction}. There we combine ideas on crossing-preserving vesselness maps from \cite{duits2018optimal,hannink2014crossing,Zhang2016}.
\begin{remark} 
Within $\mathcal{G}$ and $\mathcal{F}$
in (\ref{Finslerfin}) we set $\zeta^2=0.01=g_{11}/g_{22}$
as relative costs for sideward motion, recall (\ref{eq:LIMTFDD}). Ideally we want this to be high, but as we will prove in the numerics section (Subsection~\ref{subsec:ExtensionAFM_scheme}), a spatial anisotropy of 
$\zeta^2=0.01$ still guarantees numerical accuracy. We follow \cite{bekkers2015pde,duits2018optimal} and set bending stiffness $\xi^2=g_{11}=0.01$ and $g_{33}=1$.
\end{remark}
\begin{proposition}
Metric tensor field $\mathcal{G}^U$ given by Eq.~\eqref{Finslerfin} is indeed  data-driven left-invariant (i.e. satisfying (\ref{eq:DDLIF})).
\end{proposition}
\begin{proof}
See Lemma~\ref{lemma:DDLIF} in Appendix~\ref{app:GisDDLIV}.
\end{proof}

\begin{figure}
    \centering
    \begin{subfigure}[b]{0.23\textwidth}
        \includegraphics[width=\textwidth]{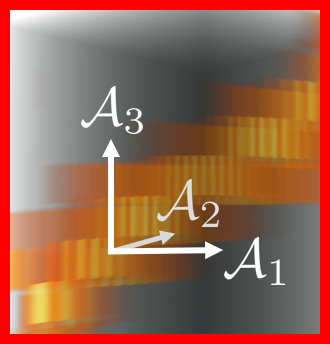}
        \caption{Left-invariant frame.}
        \label{fig:FramesLIF}
    \end{subfigure}
    \hfill
    \begin{subfigure}[b]{0.23\textwidth}
        \includegraphics[width=\textwidth]{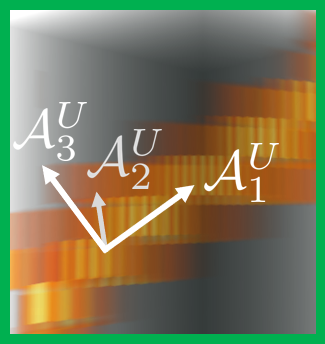}
        \caption{Data-Driven Left-invariant frame.}
        \label{fig:FramesDDLIF}
    \end{subfigure}
    \caption{Visualizations of the left-invariant frame and the data-driven left-invariant frame in $\mathbb{M}_2$. Locally along one of the spirals in the orientation score depicted in Fig.~\ref{fig:OrientationScores}. In Fig.~\ref{fig:FramesLIF}, the main direction of $\mathcal{A}_1$ is not properly aligned with the underlying 3D structure, whereas in Fig.~\ref{fig:FramesDDLIF} $\mathcal{A}_1^U$ is.}
    \label{fig:Frames}
\end{figure}

Next, we list a few remarks that underpin and explain our specific choice of metric tensor field.
\begin{remark}
    In geometric image analysis \cite{Romeny2003FrontEnd}, eigenvectors of the Hessian typically provide a local coordinate frame along lines. 
    In orientation scores, this is not different \cite{franken2009crossing}. In $\mathbb{M}_2$, Hessians $HU=\nabla^{[+],*}\mathrm{d}U$ are not symmetric and we rely on a singular value decomposition via the dual norm in \eqref{Finslerfin} which only relies on the symmetric product, see Remark~\ref{rem:Hessian}.
\end{remark}
\begin{remark}
Formally speaking, the (old) metric tensor fields $\mathcal{G}$ and  asymmetric version $\mathcal{F}$ are also data-driven if it comes to scalar-adaptation via cost function $C$, but as they do not adapt for any kind of directional data-adaptation (as illustrated in Fig.~\ref{fig:Frames}) we do \emph{not} refer to them as `the data-driven model'.    \end{remark}
\begin{remark}\label{rem:Hessian}
Via the identification
(\ref{ID}) we write $\ul{p}=(x,y,\theta)$ for short.
Then in fixed coordinates 
on $\R^2 \times \R/(2\pi \mathbb{Z})$
one may write $\dot{\p}=(\dot{x},\dot{y},\dot{\theta})^\top$. 
The dual norm expression $\left\|\left. HU\right|_\p (\dot{\p},\cdot)\right\|_*^2$ in (\ref{Finslerfin}) then boils down to a straightforward Euclidean norm:
\begin{align}
\!\left\|\!\left. HU\right|_\p \! (\dot{\p},\cdot)\right\|_*^2\!=\!\left\|M_\xi\!{\small\left(\begin{array}{ccc}
U_{xx}(\p) & U_{xy}(\p) & U_{x\theta}(\p)+U_{y}(\p)  \\
U_{yx}(\p) & U_{yy}(\p) & U_{y\theta}(\p)-U_{x}(\p)  \\
U_{\theta x}(\p) & U_{\theta y}(\p) & U_{\theta\theta}(\p)
\end{array}
\right)^{\!\!\top\!}}
\dot{\p}\right\|^2\label{eq:dualnorm}
\end{align}
where $M_\xi=\textrm{diag}\,(\xi^{-1},\xi^{-1},1) \in \R^{3\times 3}$.
\end{remark}
For details of Hessians of functions on manifolds with a connection, see Appendix~\ref{app:GisDDLIV}. For now let us focus on the notion of data-driven left-invariant frames, where we improve upon the `Locally Adaptive Derivatives (LADs)' in \cite{Zhang2016,duits2016locally}.

\begin{definition}[Data-Driven Left-Invariant Frame]
\label{def:DDLIV}\\
    Any data-driven metric tensor field $\mathcal{G}^U$ can be diagonalized:
    \begin{align}
    \label{eq:DDMTFDiagonalization}
        \mathcal{G}^U=\sum_{i=1}^3 \alpha_i^U(\cdot)\;\omega_U^i \otimes \omega_U^i   
    \end{align}
    and this defines the positively oriented data-driven left-invariant co-frame $\left\{\omega_U^i\right\}_{i=1}^3$, dual to the primal frame $\left\{\mathcal{A}_j^U\right\}_{j=1}^3$
    related by 
    $\langle \omega_U^i,\mathcal{A}_j^U\rangle=\delta_j^i$.
\end{definition}
\begin{remark}[Advantages of our data-driven metric and frame]\newline
    The local frame of reference $\{\mathcal{A}_{i}^U\}$ depends on the image data, cf.~Fig.~\ref{fig:Frames}. In fact, Eq.~\eqref{eq:DDMTFDiagonalization} is used to define the dual of the data-driven left-invariant frame via diagonalisation. This deviates from LADs in previous work \cite{duits2021Springer,Zhang2016}. 
    We now have the advantage of coercivity  
\begin{equation} \label{coerc}
    \mathcal{G}^U\geq  \mathcal{G}\geq \delta>0,
\end{equation}
    recall \eqref{eq:boundariesCost}, independent of the orientation score data $U$, which makes the tracking algorithms unconditionally stable. Furthermore, we now have another advantageous property over LADs, namely that $\mathcal{A}_i^U=\mathcal{A}_i$ for $U$ constant.
\end{remark}


    
In order to calculate distances using the new data-driven metric tensor field, we need to introduce the data-driven Riemannian distance.
\begin{definition}[Data-Driven Riemannian Distance]\label{def:DDRiemannianDistance}\newline
    The data-driven Riemannian distance $d_{\mathcal{G}^U}$ from a point $\p\in\mathbb{M}_2$ to a point $\q\in\mathbb{M}_2$ is given by
    \begin{align*}
        d_{\mathcal{G}^U}(\p,\q)=
        \underset{{\scriptsize\begin{array}{c}\gamma\in \Gamma_1,\\
\gamma(0)=\p, \\ \gamma(1)=\q
\end{array}}}{\inf}
        \int_0^1\sqrt{\mathcal{G}_{\gamma(t)}^U\left(\dot{\gamma}(t),\dot{\gamma}(t)\right)}\;\mathrm{d}t
        \numberthis\label{eq:distance}
    \end{align*}
    where $\Gamma_1$ was defined in (\ref{SpaceCurves}), and $\dot{\gamma}(t):=\frac{\mathrm{d}}{\mathrm{d}t}\gamma(t)$.
\end{definition}
\begin{remark}
    \mbox{If image $U$ is constant, 
    then 
$\mathcal{G}^U=\mathcal{G}$, $d_{\mathcal{G}^U}=d_{\mathcal{G}}$.}
\end{remark}
\begin{remark}
    Note that this distance can always be transformed to a quasi-distance when we are working in with the forward gear version of the model:
    \[
      d_{\mathcal{F}^U}(\p,\q)=
        \underset{{\scriptsize\begin{array}{c}\gamma\in \Gamma_1,\\
\gamma(0)=\p, \ \gamma(1)=\q
\end{array}}}{\inf}
        \int_0^1\mathcal{F}^U(\gamma(t),\dot{\gamma}(t))\,\mathrm{d}t.
        \]
    
\end{remark}
Using the new data-driven metric frame, recall Def.~\ref{def:DDLIV}, we introduce the data-driven Cartan plus connection, which will be used to express `short' and `straight' curves 
in Section~\ref{sec:ThShortVsStraight}.
\begin{definition}[Data-Driven Cartan Plus Connection]\\
    The data-driven Cartan plus connection is given by
    \begin{align*}
        \nabla^U:=\!\sum_{k=1}^n\left(\sum_{i=1}^n\!\omega_U^i\otimes\left(\mathcal{A}_{i}^U\circ\omega_U^k\right)\!+\!\!\sum_{i,j=1}^n\!\!\omega_U^i\otimes\omega_U^j c_{ij}^k\right)\mathcal{A}_k^U.
    \end{align*}
\end{definition}
Explicit coordinate expressions will follow in Lemma \ref{lemma:ConnectionInComponents}.

In Table~\ref{tab:prevVsCurrentWork}, an overview of the notation used for the new concepts introduced in this work and concepts introduced in earlier work is given.
\begin{table}
    \centering
        \begin{tabular}{c|ccc}
        \multirow{3}{*}{} & Metric   & \multirow{3}{*}{Diagonalization} & Cartan \\
        &Tensor&&Connection\\
        &Field&&\\\hline
        Earlier work & $\mathcal{G}$ & $\left\{\omega^i\right\}_{i=1}^n$ & $\nabla^{[+]}$ \\
        Current work & $\mathcal{G}^U$ & $\left\{\omega_U^i\right\}_{i=1}^n$ & $\nabla^{U}$ \\
    \end{tabular}
    \caption{Comparison of (notation of) current and previous work. Diagonalization is w.r.t. dual frame associated to the frames depicted in Fig.~\ref{fig:Frames}.}
    \label{tab:prevVsCurrentWork}
\end{table}

\begin{figure*}
    \begin{center}
    \begin{minipage}{0.49\textwidth}
    \begin{center} \large{Cartan Plus Connection $\nabla^{[+]}$}
    \end{center}
    \vspace*{\baselineskip}
    \end{minipage}
    \begin{minipage}{0.49\textwidth}
    \begin{center} \large{Data-driven Cartan Plus Connection $\nabla^{U}$}
    \end{center}
    \vspace*{\baselineskip}
    \end{minipage}
    \vspace*{\baselineskip}
    \includegraphics[width=0.3\textwidth,angle=180,origin=c]{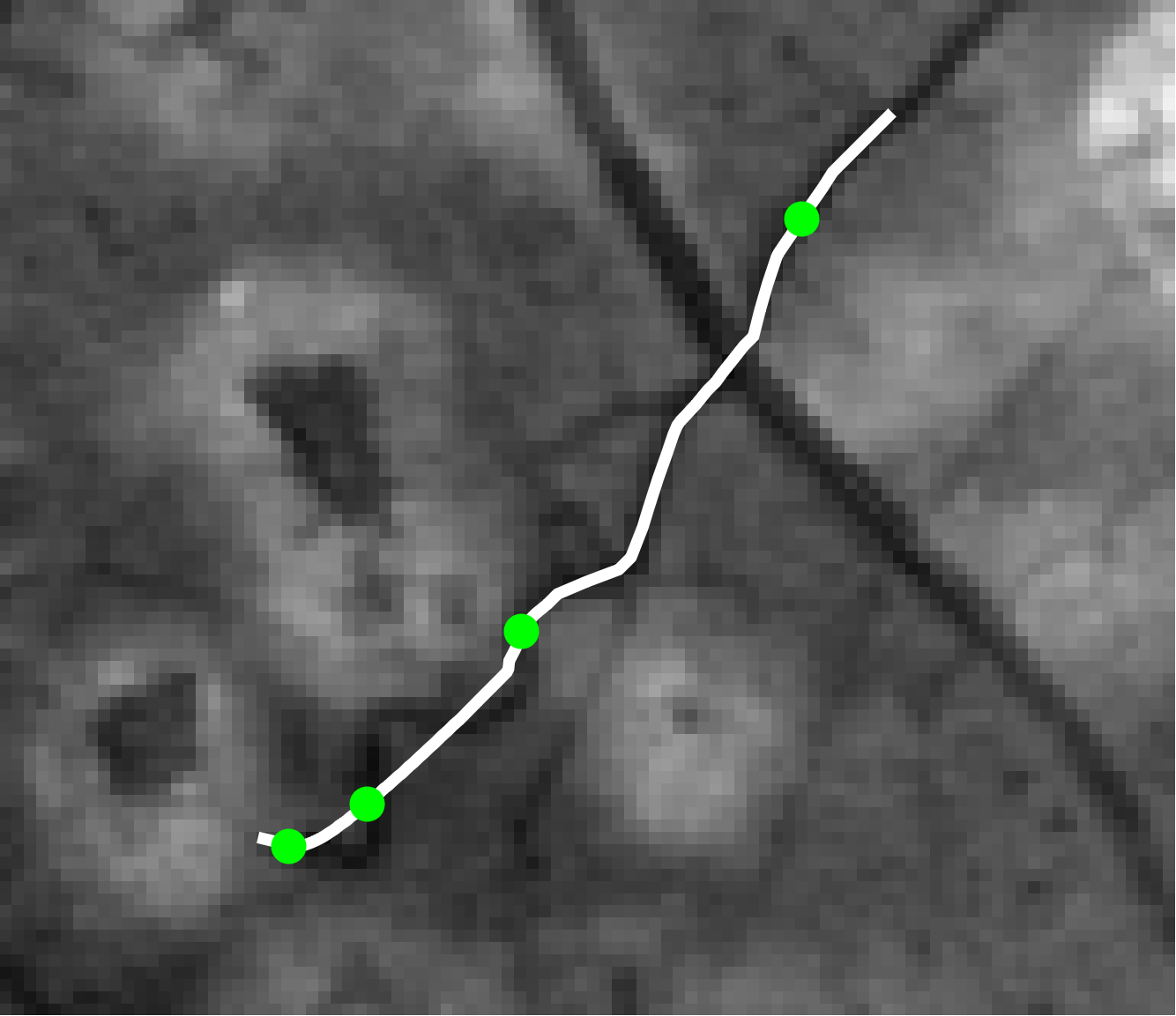}\hspace{0.2\textwidth}
    \includegraphics[width=0.3\textwidth,angle=180,origin=c]{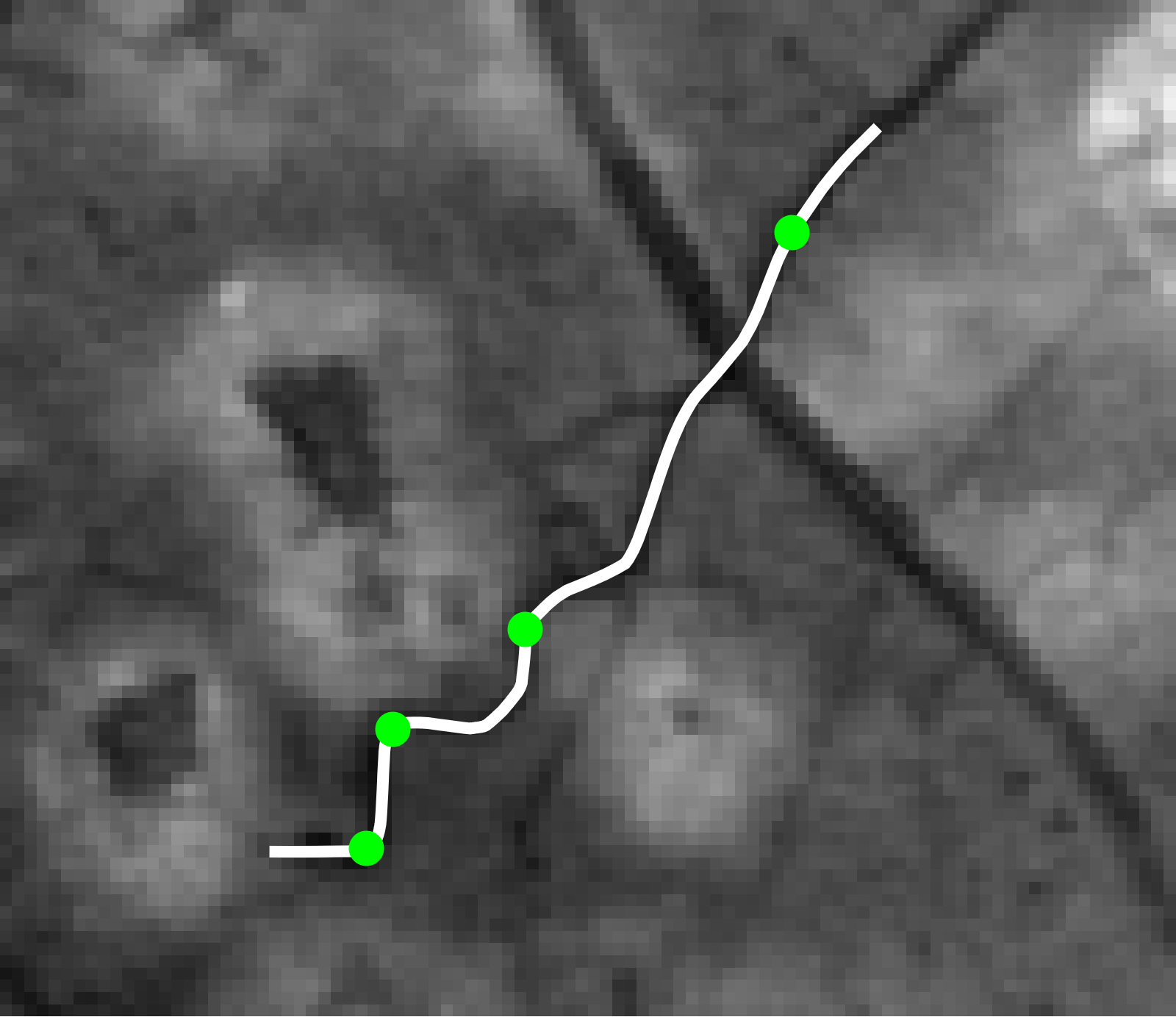}
    \includegraphics[width=0.38\textwidth]{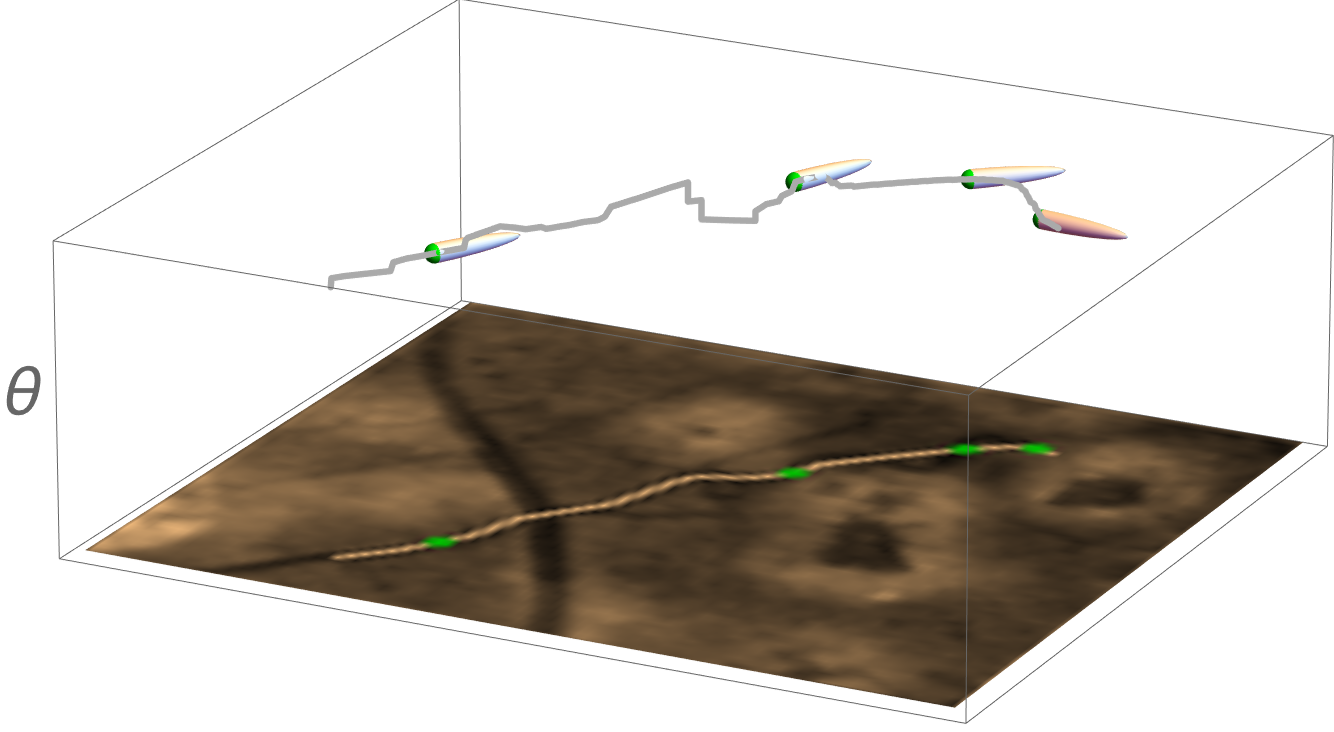}\hspace{0.1\textwidth}
    \includegraphics[width=0.38\textwidth]{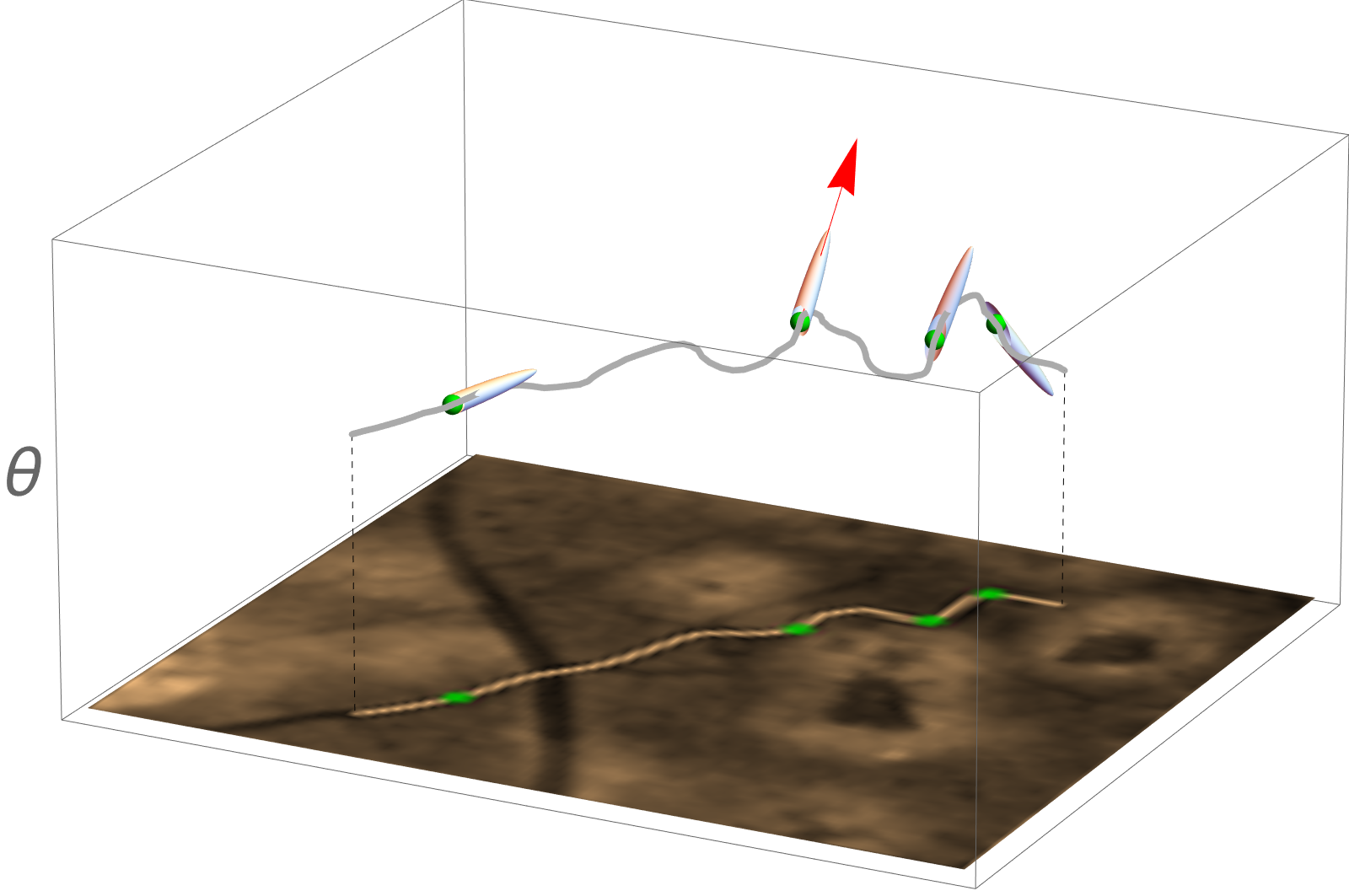}
    \includegraphics[width=0.38\textwidth]{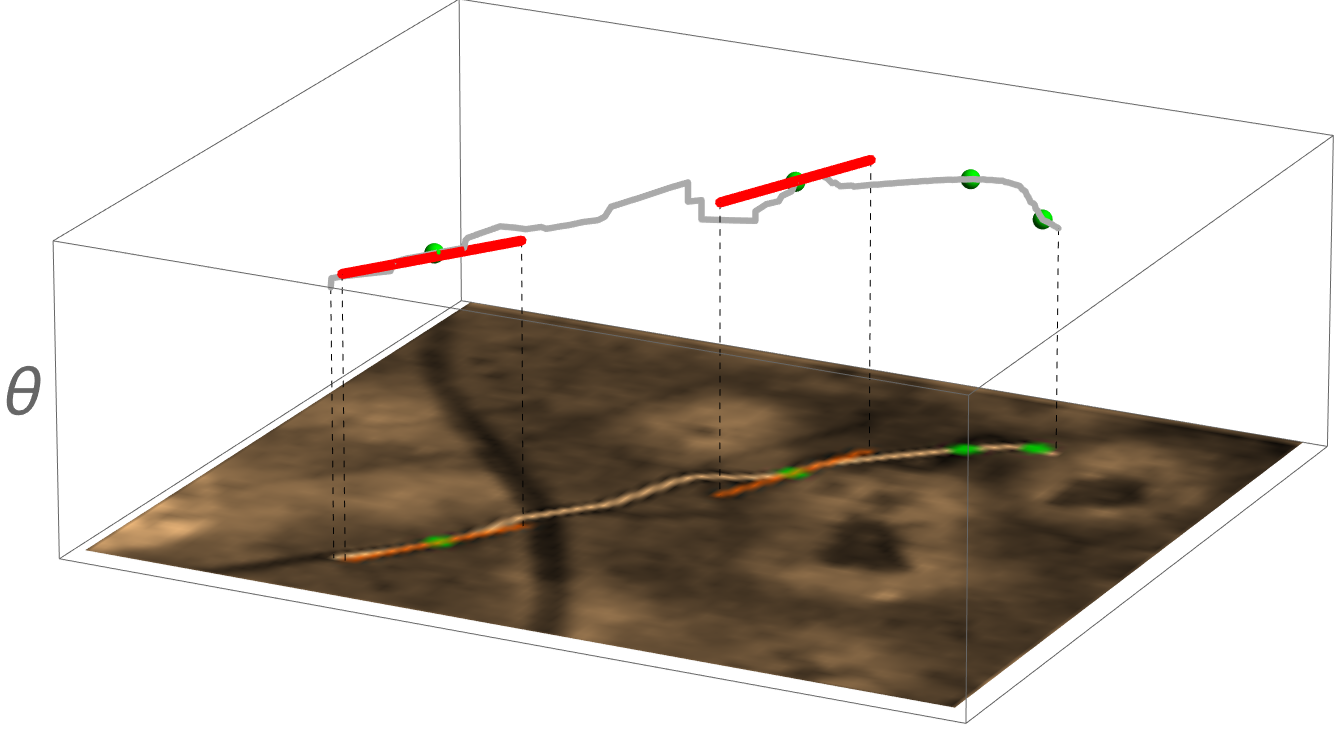}\hspace{0.1\textwidth}
    \includegraphics[width=0.38\textwidth]{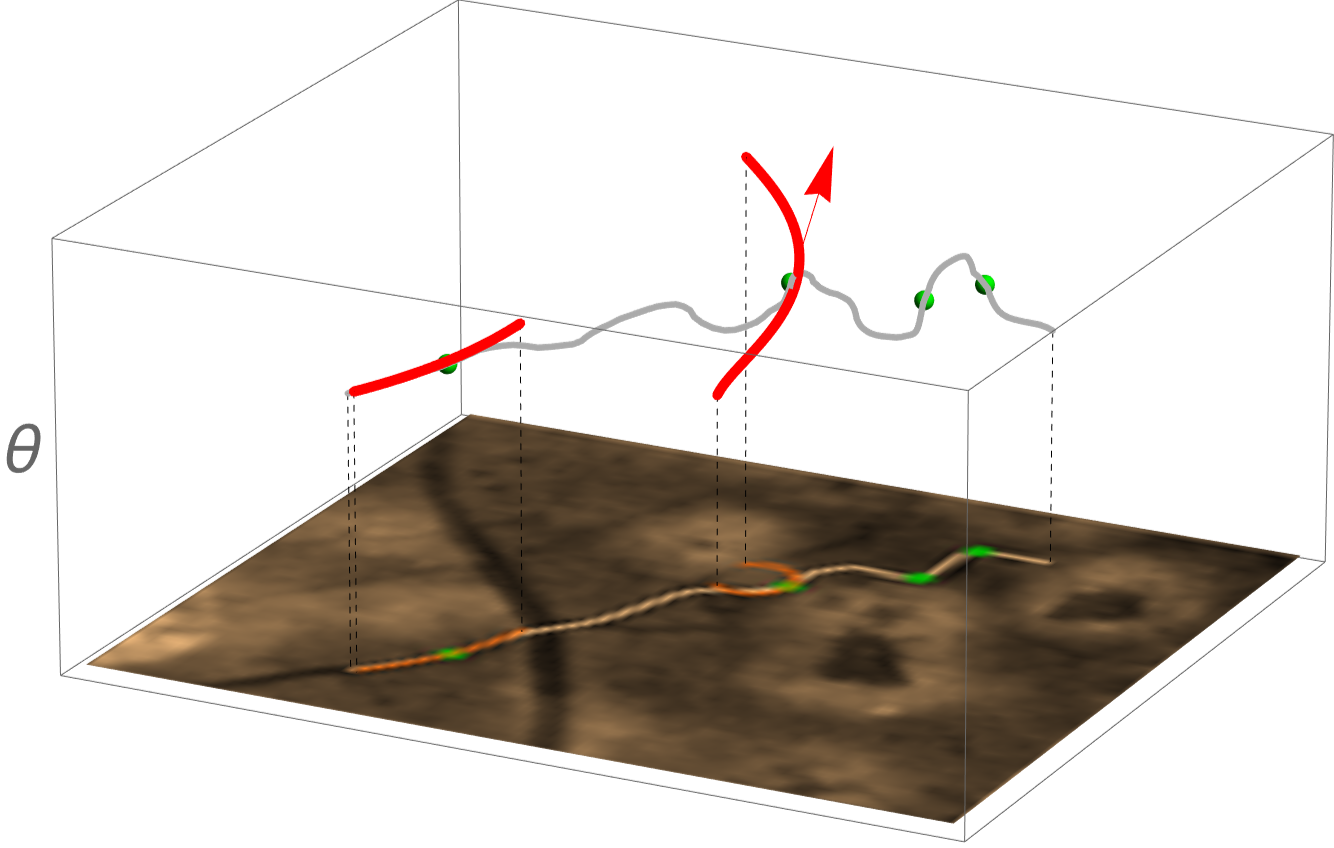}
    \end{center}
    \caption{
    The advantage of using  data-driven Cartan connections $\nabla^{U}$ instead of the non-data-driven Cartan plus connection $\nabla^{[+]}$.
    In grey, a shortest curve (geodesic) between two points in $\mathbb{M}_{2}$ is visualized, along with its spatial projection. The left geodesic has parallel momentum w.r.t. $\nabla^{U}$ (cf.Thm~\ref{th:th}) and the right w.r.t. $\nabla^{[+]}$ \cite[Thm.1]{duits2021Springer}.  The new geodesic better adapts for curvature (and spatial alignment). This is also visible in the corresponding control sets \eqref{eq:RelationControlsetGandF} depicted by the white closed surfaces above at several green points on the geodesics. The red arrow indicates the principal direction of the local metric tensor (left: $\mathcal{G}$, right: data-driven $\mathcal{G}^U$). The control sets belonging to $\nabla^{[+]}$ are only aligned to the underlying structure in the spatial domain, whereas the control sets belonging to $\nabla^U$ align with the appropriate curvature in the tangent space as well. In the bottom row, we depict exponential curves through the green points with a tangent in the principal direction (left of $\mathcal{G}$, right of $\mathcal{G}^U$). They are straight-curves of $\nabla^{[+]}$ (left) and $\nabla^{U}$ (right),  and `steer' the geodesic tracking as we will show in Thm~\ref{th:th}.}
    \label{fig:ShortVsStraightCurves}
\end{figure*}

In Figure~\ref{fig:ShortVsStraightCurves}, the exponential curves and the control sets for both discussed Cartan connections, $\nabla^{[+]}$ and $\nabla^U$ are visualized. In addition to that, the tracking results relying on different models are plotted. One sees that the data-driven Cartan connection better adapts for curvature leading to more accurate tracking results. 

\section{The New Geometric Tracking Model: Asymmetric Finsler Functions steered by Locally Adaptive Frames}\label{sec:ThShortVsStraight}

We discuss a new data-driven version of the Cartan connection. This result applies to all Lie-groups $G$ of finite dimension $\textrm{dim}(G)=n$. Note that $SE(2)\equiv\mathbb{M}_2$, but the result does not apply to all homogeneous spaces (like $\mathbb{M}_d$ for $d>2$). The notation used in this section is summarized in Table~\ref{tab:ToN}. 
\newcounter{tblEqCounter}
\setcounter{tblEqCounter}{\theequation}
\begin{table*}[ht]
    \centering
    \begin{tabular}{|llc|}\hline
     & & \\
        $T(G):=\left\{(g,\dot{g})|\dot{g}\in T_g (G)\right\}$ bases & $\left\{\partial_{x^i}\right\}_{i=1}^n;\qquad\quad\left\{\mathcal{A}_i\right\}_{i=1}^n;\qquad\quad\left\{\mathcal{A}_i^U\right\}_{i=1}^n$.&\\
        $T^*(G)$ bases & $\left\{\mathrm{d}x^i\right\}_{i=1}^n;\qquad\quad\left\{\omega^i\right\}_{i=1}^n;\qquad\quad\left\{\omega_U^i\right\}_{i=1}^n$.&\\[6pt]
        Notation duals&$\langle\hat{a},b\rangle=\hat{a}(b)$.&\\[6pt]
        Relation between bases $T(G)$ and $T^*(G)$&$\langle\mathrm{d}x^i,\partial_{x^j}\rangle=\mathrm{d}x^i(\partial_{x^j})=\delta_j^i;\quad
        \langle\omega^i,\mathcal{A}_j\rangle=\delta_j^i;\quad
        \langle\omega_U^i,\mathcal{A}_j^U\rangle=\delta_j^i$.&\numberTblEq{eq1}\label{duality}\\&&\\
        LIV metric tensor field&$\mathcal{G}=\sum\limits_{i=1}^n g_{ii}\omega^i\otimes\omega^i$&\\
        &For specific choices of coefficients $g_{ii}$, see \eqref{eq:LIMTFDD}.&\\&&\\
        Data-driven LIV metric tensor field & $\mathcal{G}^U = \sum \limits_{i=1}^n \alpha_i^U(\cdot) \omega^{i}_{U} \otimes \omega^{i}_U$ & \\
        & \textrm{that must satisfy }\eqref{eq:DDLIF}\textrm{. We choose to use }(\ref{Finslerfin}) and write $\alpha_i:=\alpha_i^U$.&\\ 
         &  \textrm{For its asymmetric Finslerian extension }$\mathcal{F}^U$ \textrm{ see also }\eqref{Finslerfin}. & \\ & & \\
        Velocity components &$\dot{\gamma}=\sum\limits_{i=1}^n\dot{\overline{\gamma}}^i\partial_{x^i}=\sum\limits_{i=1}^n\dot{\gamma}^i\mathcal{A}_i=\sum\limits_{i=1}^n\dot{\tilde{\gamma}}^i\mathcal{A}_i^U$.&\\
        &&\\
        Momentum components &$\lambda=\sum\limits_{i=1}^n\overline{\lambda}_i\mathrm{d}x^i=\sum\limits_{i=1}^n\lambda_i\omega^i=\sum\limits_{i=1}^n\tilde{\lambda}_i\omega_U^i$.&\\
        &&\\
        Lagrangian & $\mathcal{L}(\gamma,\dot{\gamma})=\frac{1}{2}|\mathcal{G}^U_{\gamma}(\dot{\gamma},\dot{\gamma})|^2= \sum_{i=1}^n \alpha_{i}^U(\gamma) |\dot{\tilde{\gamma}}^i|^2$.&\\&&\\
        Hamiltonian & $\mathfrak{h}(\lambda)=\frac{1}{2}\sum\limits_{i=1}^n\tilde{\lambda}_i\tilde{\lambda}^i=\frac{1}{2}\sum\limits_{i=1}^n\tilde{\lambda}_i\alpha^i(\cdot)\tilde{\lambda}_i$ $\textrm{\ \ \  with }\alpha^i=\alpha_i^{-1}$.&\numberTblEq{eq2}\label{eq:Hamiltonian} \\&&\\
        Fundamental Symplectic form &$\sigma\equiv\sum\limits_{i=1}^n\mathrm{d}x^i\wedge\mathrm{d}\overline{\lambda}_i=\sum\limits_{i=1}^n\omega^i\wedge\dbar \lambda_i=\sum\limits_{i=1}^n\omega_U^i\wedge\dbar\tilde{\lambda}_i,$&\\& where $\langle\dbar\lambda_i,\partial_{\lambda_j}\rangle=\delta_j^i$, but $\dbar\lambda_i\neq\mathrm{d}\lambda_i$.&\\&&\\
        Hamiltonian Flow $\sigma(\overrightarrow{\mathfrak{h}},\cdot)=\mathrm{d}\mathfrak{h}$ on $T^*M$ &
        {$
            \begin{cases}
            \begin{matrix}
                \dot{\nu}=\overrightarrow{\mathfrak{h}}(\nu)\\
                \nu(0)=\nu_0,
            \end{matrix}
            \end{cases}
        $} &\\
        &where $\nu=(\gamma,\lambda)$, with $\gamma$ denoting the geodesic and $\lambda$ the momentum along it. &\\&&\\
        Lie Bracket &  $[X,Y]=XY-YX$.&\numberTblEq{eq3}\label{eq:LieBracket}\\&&\\
        Poisson Bracket & $\left\{f,g\right\}=\sum\limits_{i=1}^n\partial_{x^i}g\frac{\partial f}{\partial\overline{\lambda}_i}-\frac{\partial g}{\partial\overline{\lambda}_i}\partial_{x^i}f=\sum\limits_{i=1}^n(\mathcal{A}_i^U g)\frac{\partial f}{\partial\tilde{\lambda}_i}-(\mathcal{A}_i^U f)\frac{\partial g}{\partial\tilde{\lambda}_i}.$&\\&&\\
        Structure Functions &
        {$
        [\mathcal{A}_i^U,\mathcal{A}_j^U]=\sum\limits_{k=1}^n\tilde{c}_{ij}^k(\cdot)\, \mathcal{A}_k^U
        $.}&\numberTblEq{eq4}\label{eq:StructureConstants}\\&&\\
        Structure Constants&{$
        [\mathcal{A}_i,\mathcal{A}_j]=\sum\limits_{k=1}^n c_{ij}^k\, \mathcal{A}_k
        $}, where $c_{ij}^k=-c_{ij}^k$.&\numberTblEq{eq5}\label{eq:StructureConstants2}\\&&\\
        Data Driven Cartan Connection  & $        \nabla^{U}=\sum\limits_{k=1}^n\left(\sum\limits_{i=1}^n\omega_U^i\otimes\mathcal{A}_i^U\circ\omega_U^k(\cdot)+\sum\limits_{i,j=1}^n\left(\omega_U^i\otimes\omega_U^j\right)\tilde{c}_{ij}^k(\cdot)\right)\mathcal{A}_k^U$.& see also \eqref{connectcomp}\\&&\\
        Dual Data Driven Cartan Connection & $
        \left(\nabla^{U}\right)^*=\sum\limits_{i=1}^n\left(\sum\limits_{j=1}^n\omega_U^j\otimes\left(\mathcal{A}_j^U\circ\mathcal{A}_i^U\right)+\sum\limits_{j,k=1}^n\left(\omega_U^j\otimes\mathcal{A}_k^U\right)\tilde{c}_{ij}^k(\cdot)\right)\omega_U^i$. & see also \eqref{dualconnectcomp} \\ & & \\\hline
    \end{tabular}
    \caption{Table of Geometric Tools and Notations}
    \label{tab:ToN}
\end{table*}
\setcounter{equation}{\thetblEqCounter}


We consider a locally adaptive frame $\left\{\mathcal{A}_{i}^U\right\}_{i=1}^n$ with dual frame $\left\{\omega_U^i\right\}_{i=1}^n$. This can be any well-defined frame that depends on the underlying data. 
The (data-driven) metric tensor field that is considered, is given by \eqref{eq:DDMTFDiagonalization}.
The data-driven terms can adapt for curvature and deviation from horizontality where the direction of the left-invariant frame deviates from the underlying line structure.

\subsection{Combine Optimally Straight and Short: A new Data-Driven Version $\nabla^{U}$ of the Cartan Connection}\label{sec:Cartan}
In previous works, the Cartan plus connection, which relies on the left-invariant frame, has been used to describe straight and shortest curves in Lie groups \cite{Duitsbookchapter}. 
However, this frame is not always adequate in multi-orientation image processing as it
does not always align perfectly with the underlying line structures in the orientation scores (see Figure~\ref{fig:Frames}). 
To improve the tracking results, we, therefore, switch to using a data-driven Cartan connection associated with the data-driven metric tensor field $\mathcal{G}^U$ given by (\ref{Finslerfin}).
Let us first define what we mean by a `data-driven Cartan connection'.
\begin{definition}
  The data-driven Cartan connection and its corresponding dual are given by
    {\small\begin{equation}
    \begin{array}{l}
    \nabla^{U}= 
       \sum\limits_{k=1}^n\left(\sum\limits_{i=1}^n\omega_U^i\otimes\left(\mathcal{A}_{i}^U\circ\omega_U^k\right)+\sum\limits_{j=1}^n\left(\omega_U^i\otimes\omega_U^j\right)\tilde{c}_{ij}^k\right)\mathcal{A}_k^U,\numberthis\label{eq:connection}\\
    \end{array}
    \end{equation}
    \begin{equation}
    \begin{array}{l}
    \left(\nabla^{U}\right)^*= 
    \sum\limits_{i=1}^n\left(\sum\limits_{j=1}^n\omega_U^j\otimes\left(\mathcal{A}_j^U\circ\mathcal{A}_{i}^U \right)+\sum\limits_{k,j=1}^n\left(\omega_U^j\otimes\mathcal{A}_k^U\right)\tilde{c}_{ij}^k\right)\omega_U^i,\numberthis\label{eq:dualConnection}
    \end{array}
    \end{equation}}
    where $\left(\nabla^U\right)_X^*\lambda:=\left(\nabla^U\right)^*(X,\lambda)$ and $\nabla_X^U Y:=\nabla^U(X,Y)$.
\end{definition}
\begin{remark}
The relation between $\nabla$ and its dual $\nabla^*$ is 
\begin{equation} \label{dualCON}
\langle 
\nabla^*_X \lambda, Y \rangle := X \langle \lambda,Y \rangle  -\langle 
\lambda, \nabla_X Y \rangle 
\end{equation}
for all vector fields $X,Y$ and all covector fields $\lambda$ on $G$, which may be interpreted as a product rule for the pairing between the vectors and co-vectors. In particular for $X=\mathcal{A}^U_i$, $Y=\mathcal{A}_j^U$ and $\lambda=\omega_U^k$ we get
\[
\begin{array}{ll}
\tilde{c}^{k}_{ji}&\overset{\eqref{eq:dualConnection}}{=}\langle 
(\nabla^{U})^*_{\mathcal{A}_i^U}
\omega^{k}_U,
\mathcal{A}_j^U \rangle\\
 &\overset{\eqref{dualCON}}{=}
\mathcal{A}_i^U
\langle 
\omega^k_U,\mathcal{A}^U_j
\rangle 
 -\langle \omega^{k}_U, \nabla^U_{\mathcal{A}_i^U} \mathcal{A}_j^U\rangle\\
 &=
\mathcal{A}_i^U(\delta_{j}^k) 
 -\langle \omega^{k}_U, \nabla^U_{\mathcal{A}_i^U} \mathcal{A}_j^U\rangle \\
 &=-\langle \omega^{k}_U, \nabla^U_{\mathcal{A}_i^U} \mathcal{A}_j^U\rangle \overset{\eqref{eq:connection}}{=}
-\tilde{c}_{ij}^k.
\end{array}
\]
\end{remark}
In the next lemma, we will express the data-driven Cartan connection and its corresponding dual  explicitly in coordinates, which will provide us an expression on which we will build in the proof of our main theorem, Theorem~\ref{th:th}.

\begin{lemma}\label{lemma:ConnectionInComponents}
When expressing Eq.~\eqref{eq:connection} and \eqref{eq:dualConnection} more explicitly in data-driven left-invariant frame components (gauge frame components for short), one finds
\begin{equation} \label{connectcomp}
\left(\nabla^{U}\right)_XY=\sum_{k=1}^n\left(\dot{\tilde{y}}^k+\sum_{i,j=1}^n \tilde{c}_{ij}^k(\cdot)\tilde{x}^i\tilde{y}^j\right)\mathcal{A}_k^U,
\end{equation}
and for the dual connection
\begin{equation} \label{dualconnectcomp}
    \left(\nabla^{U}\right)_X^* \lambda=\sum_{i=1}^n\left(\dot{\tilde{\lambda}}_i+\sum_{k,j=1}^n\left(\tilde{x}^j\tilde{\lambda}_k\right)\tilde{c}_{ij}^k(\cdot)\right)\omega_U^i,
\end{equation}
where $X=\sum\limits_{i=1}^n\tilde{x}^i\mathcal{A}_{i}^U|_\gamma$
, $Y=\sum\limits_{i=1}^n \tilde{y}^i\mathcal{A}_{i}^U|_\gamma$ and $\lambda=\sum\limits_{i=1}^n\tilde{\lambda}_i\omega_U^i$, 
and where derivations of the components of $Y$ and $\lambda$ equal
\[
\begin{array}{rll}
\dot{\tilde{y}}^k(t) &:=\frac{\mathrm{d}}{\mathrm{d}t}\tilde{y}^k(\gamma(t))&=\left(X\left(\tilde{y}^k\right)\right)(\;\gamma(t)\;), \\
\dot{\tilde{\lambda}}_i(t)&:=\frac{\mathrm{d}}{\mathrm{d}t}\tilde{\lambda}_i(\gamma(t))&=\left(X\left(\tilde{\lambda}_i\right)\right)\, (\;\gamma(t)\;),
\end{array}
\] 
along a flow-line\footnote{A curve $\gamma$ satisfying $\dot{\gamma}(t)=X_{\gamma(t)}$.} $\gamma:[0,1]\to \mathbb{M}_2$ of smooth vector field $X$.
\end{lemma}
\begin{proof}
    See Appendix~\ref{app:lemma}.
\end{proof}
\subsection{Main Theorems}
 Our goal is to analyze and structure the Hamiltonian flow belonging to the new data-driven geometric model determined by a data-driven metric tensor field $\mathcal{G}^U$. For convenience, we restrict ourselves in our main theorem to the case where the homogeneous space equals a full finite-dimensional Lie group $G$ as the base manifold.
 We are mainly 
 interested in the case 
 $G=SE(2) \equiv \mathbb{M}_{2}$ with $n=3$ and with data-driven metric tensor field $\mathcal{G}^U$ given by (\ref{Finslerfin}).

In geometric control curve optimization problems, the Hamiltonian flow 
is a powerful mechanism \cite{marsden2013introduction,montgomery2002tour,AgrachevSachkov,Agrachev2019introductionSubRiemannianGeometry,moiseev2010maxwell}. It typically allows us to analyze the behavior of \emph{all} geodesics (and their momentum) together, see e.g. \cite{moiseev2010maxwell}. In the left-invariant (non-data-driven) setting, the Hamiltonian flows were characterized \cite[Thm.1\&2]{duits2021Springer} via the plus Cartan connection, where shortest curves (geodesics) have parallel momentum. 
It has led to exact formulas \cite{duits2014association,duits2016sub} for specific settings (uniform cost case in the left-invariant model $\mathcal{G}$ given in \eqref{eq:LIMTFDD}, $C=1$ which was introduced on $\mathbb{M}_2\equiv SE(2)$ by Citti and Sarti \cite{Citti} and deeply analyzed by Sachkov \cite{Sachkov}). Remarkably, the Cartan plus connection $\nabla^{[+]}$ has torsion resulting in different `straight curves' (parallel velocity) and `shortest curves' (parallel momentum), and it underlies many multi-orientation image analysis applications \cite{duits2021Springer}.  

Before stating the main theoretical result, that generalizes \cite[Thm.1\&2]{duits2021Springer} to \emph{data-driven} metric tensor fields $\mathcal{G}^U$,
we introduce the concepts of parallel momentum and parallel velocity. They are now determined by
the data-driven Cartan connection $\nabla^U$ and its dual $\left(\nabla^U\right)^*$.
\begin{definition}[Parallel momentum]\label{def:parallelMomentum}
Let $\gamma:[0,1]\to G$ be a smooth curve. Then, the curve $\gamma(\cdot)$ has $\nabla^{U}$ parallel momentum $\lambda(\cdot)$
    \begin{align} \label{eq:parallelMomentum}
        \Leftrightarrow
        \begin{cases}
            \left(\nabla^{U}\right)_{\dot{\gamma}}^*\lambda=0\\
            \mathcal{G}^U \dot{\gamma}=\lambda.
        \end{cases}
    \end{align}
\end{definition}
\begin{definition}[Parallel velocity]\label{def:parallelVelocity}
    Let $\gamma:[0,1]\to G$ be a smooth curve. Then, the curve $\gamma(\cdot)$ has parallel velocity $\dot{\gamma}(\cdot)$ w.r.t. connection $\nabla^{U}$
        \begin{align}
            \Leftrightarrow\left(\nabla^{U}\right)_{\dot{\gamma}}\dot{\gamma}=0.\label{eq:parallelVelocity}
        \end{align}
\end{definition}
\begin{remark}
Using the antisymmetry of the structure functions (\ref{eq:StructureConstants}) and (\ref{connectcomp}) in Lemma~\ref{lemma:ConnectionInComponents} we can rewrite Eq.~\eqref{eq:parallelVelocity} to
    \begin{align*}
    \left(\nabla^{U}\right)_{\dot{\gamma}}\dot{\gamma}=0 \Leftrightarrow
       \exists_{ (c^1,\ldots,c^n) \in \R^n} \text{ constant s.t. }\dot{\gamma}=\sum_{i=1}^n c^i \mathcal{A}_{i}^U|_\gamma.
    \end{align*}
\end{remark}

Next, we formulate the main theoretical results where we generalize the main results \cite[Thm.1\&2]{duits2021Springer} from the left-invariant setting
$(G,\mathcal{G})$ with connection $\nabla^{[+]}$, 
to the new data-driven geometric models $(G,\mathcal{G}^U)$ with connection $\nabla^U$. 
In more detail, the next theorem shows 
\begin{enumerate}
    \item how to compute globally optimal distance minimizers in a geometry that relies on data-driven left-invariant frames: These geodesics have parallel momentum w.r.t. connection $\nabla^{U}$ (Def. \ref{def:parallelMomentum}).
    \item that the locally optimal straight curves are the straight curves w.r.t. connection $\nabla^{U}$: These curves have parallel velocity (i.e. are auto-parallel) w.r.t. $\nabla^{U}$ (Def. \ref{def:parallelVelocity}).
\end{enumerate}
Note that the equation for geodesics of the new data-driven model $(\mathbb{M}_{2},\mathcal{G}^U)$ gives a wild expression in the left-invariant frame. In the fixed frame it is even worse. 
However, our new tool of the connection $\nabla^U$ expressed in the locally adaptive frame $\omega_{i}^U$ allows us to describe these
geodesic equations (and the Hamiltonian flow) concisely and intuitively by the next theorem, using the tools listed in Table \ref{tab:ToN}.
\begin{theorem}[straight and shortest curves: parallel velocity and momentum w.r.t. connection $\nabla^{U}$]\label{th:th} \\[6pt]
    In a Riemannian manifold $(G, \mathcal{G}^U)$ with data-driven left-invariant metric tensor field $\mathcal{G}^U$ satisfying \eqref{eq:DDLIF}, and with induced Riemannian metric $d_{\mathcal{G}^U}$  \eqref{eq:distance}, we have:
    \begin{itemize}
        \item $\gamma$ is a straight curve with respect to $\nabla^U$  $\overset{\text{def.}}{\Leftrightarrow}$
        $\nabla_{\dot{\gamma}}^{U}\dot{\gamma}=0$
        \begin{align*}
            \Leftrightarrow\exists\ (c^1,\ldots,c^n) \in \R^n \text{ constant s.t. }\dot{\gamma}=\sum_{i=1}^n c^i \mathcal{A}_{i}^U|_\gamma.
        \end{align*}
        \item $\gamma$ is a shortest curve with respect to $\nabla^U$ 
        $\Rightarrow$ 
        the curve-momentum pair $\nu=(\gamma, \lambda): [0,1] \to T^*(G)$ satisfies the Hamiltonian flow 
        \begin{equation} \label{relevant}     
            \dot{\nu}=\overrightarrow{\mathfrak{h}}(\nu)
            \Leftrightarrow\begin{cases}
                (\nabla^{U})_{\dot{\gamma}}^*\lambda=0\\
                \mathcal{G}^U \dot{\gamma}=\lambda,
            \end{cases}
        \end{equation}
where 
one has the following pull-back symmetry\footnote{For the definition of pullback of a dual connection, see Remark~\ref{remark:Pullback} in Appendix~\ref{app:th}.} of the data-driven Cartan connection
        \begin{equation} \label{symmetry}
        (L_{\q})^* (\nabla^{\mathcal{L}_{\q} U})^* =(\nabla^U)^* \textrm{ for all }\q  \in G,
        \end{equation}
        with left actions $L$ and $\mathcal{L}$ given by (\ref{straightL}) and (\ref{mathcalL}) respectively.
    \end{itemize}
    The shortest curve $\gamma:[0,1] \to G$ with $\gamma(0)=\g$ and $\gamma(1)=\g_0$ may be computed by steepest descent backtracking on the distance map $W(\g)=d_{\mathcal{G}^U}(\g,\g_0)$
    \begin{align} \label{BT}
        \gamma(t):=\gamma^U_{\g,\g_0}(t)=
      \textrm{Exp}_{\g}(t\; v(W))   
      \qquad t \in [0,1],
    \end{align}
    where \textrm{Exp} integrates the following vector field on $G$:\\ $v(W):=-W(\g) ({\mathcal{G}^U})^{-1}
    {\rm d}W=-W(\g) \sum \limits_{k=1}^n |\alpha_k|^{-1} \mathcal{A}_k^U(W) \mathcal{A}_{k}^U$  and where  
    $W$
is the 
viscosity solution of the eikonal PDE system
    \begin{align}\label{eq:EikonalPDE}
    \begin{cases}
        \|\textrm{grad}_{\mathcal{G}^U}W\|=1, \\
        W(\g_0)=0,
    \end{cases}
    \end{align}
    where we assume $\g$ is neither a 1st Maxwell-point nor a conjugate point.
    As $v(W)$ is data-driven left-invariant, the geodesics carry the symmetry
    \begin{equation} \label{symgeods}
    \gamma_{\q \g, \q \g_0}^{\mathcal{L}_\q U}(t) = \q \;\gamma_{\g,\g_0}^U(t) \textrm{  for all }\q,\g,\g_0 \in G, t \in [0,1].
    \end{equation}
\end{theorem}
\begin{proof}
See Appendix~\ref{app:th}.
\end{proof}
\begin{remark} \mbox{\emph{(Assumptions on  point $\g$ in backtracking (\ref{BT}))}} \\
For the geodesic backtracking formulated above, we need a differentiable distance map along the path and a well-posed Hamiltonian flow (i.e. the mapping from 
$\nu(0)$ to $\nu(t)$ arising from (\ref{relevant}) must have a non-vanishing Jacobian) along the path.
If $\g$ would be a first Maxwell-point (i.e. two distinct geodesics meet for the first time at $\g$) the distance map is not differentiable at $\g$. 
If $\g$ would be a conjugate point (often limits of first Maxwell points \cite{bekkers2015pde}) then the Hamiltonian flow breaks down.  In the latter case, local optimality is lost.
In the first case, global optimality is lost. Fortunately the viscosity property of the viscosity solution \cite{evans1998partial} of (\ref{eq:EikonalPDE}) kills non-optimal fronts \cite{bekkers2015pde} and one may resort to multi-valued backtracking via sub-gradient backtracking.
\end{remark}

The next 3 remarks show how incredibly simple the Hamiltonian flow, the eikonal PDE, and the backtracking of geodesics become when expressed in the data-driven left-invariant frame.  
\begin{remark}
    Eq.~\eqref{eq:parallelMomentum} is in gauge frame components simply:
    \begin{align*}
        \begin{cases}
            \dot{\tilde{\gamma}}^i=\tilde{\lambda}^i&i=1,\ldots,n\\
            \dot{\tilde{\lambda}}^i=\sum \limits_{j=1}^n\sum \limits_{k=1}^n \tilde{c}_{ji}^k(\gamma(\cdot))\;\tilde{\lambda}_k\tilde{\lambda}^j&i=1,\ldots,n.
        \end{cases}
    \end{align*}
\end{remark}
\begin{remark}
    Eq.~\eqref{eq:EikonalPDE} is in gauge frame components simply:
    \begin{align*}
        \begin{cases}
            \sum \limits_{j=1}^n\alpha_j^U(\cdot)^{-1}(\mathcal{A}_j^U W)^2=1\\
            W(\gamma(0))=0.
        \end{cases}
    \end{align*}
    
\end{remark}
\begin{remark}\label{rem:proofeq29toEq31}
Eq.~\eqref{BT} is in gauge frame components simply:
\begin{equation} \label{BTsimple}
\dot{\tilde{\gamma}}^k=\frac{1}{W(\g)}|\alpha_{k}^U|^{-1} (\mathcal{A}_k^U W)(\tilde{\gamma}), \qquad k=1,\ldots, n.
\end{equation}
    This explains the definition of $v(W)$ below \eqref{BT}. A more explicit integration formula for (\ref{BT}) can be obtained in a similar way as in \cite{duits2016sub,duits2014association} (where exact solutions are derived for $C=U=1$) via momentum preservation laws. 
\end{remark}

\subsection{Asymmetric Finsler Extension to Automatically Deal with Bifurcations \label{ch:extendFU}} 
One can always decide to include a positive control variant, to avoid cusps in the spatial projection of geodesics in $G=SE(2)$. 
This is done by considering the geodesics in the asymmetric Finslerian manifold $(\mathbb{M}_{2},\mathcal{F}^U)$, recall (\ref{Finslerfin}), rather than the geodesics in the Riemannian manifold $(\mathbb{M}_{2},\mathcal{G}^U)$. 

It is not too hard in practice: a slight adaptation of the eikonal PDE (taking the positive part of one momentum component) will guarantee that all optimal geodesic wave-fronts propagate in the preferred forward direction around the source point, as can be observed in Figure~\ref{fig:TortuousVesselTracking} where the fronts initially move `down-left' (and not `up-right'). 

The asymmetric Finslerian model $(\mathbb{M}_{2},\mathcal{F}^U)$ is still well-posed (controllable and piecewise regular geodesics) even if $\zeta \downarrow 0$.
In fact, such asymmetric Finslerian geodesics are by construction piecewise concatenations of Riemannian geodesics and in-place rotations.
These observations follow by a straightforward 
generalization of  \cite[Thm.1,2,4]{duits2018optimal} to the data-driven setting, where the control set formulation of the geodesic distances, still applies: 
\begin{equation} \label{eqdef:themetric}
\begin{array}{l}
\!
d_{\mathcal{F}^U}(\mathbf{p},\mathbf{q})=\\[5pt]
 {\small \inf \!\small\left\{ T \geq 0 \,|\, \exists \gamma \in \Gamma_T,  \gamma(0)\!=\!\mathbf{p},  \gamma(T)\!=\!\mathbf{q}, 
 \dot{\gamma} \!\in\! \mathcal{B}_{\mathcal{F}^U}(\gamma) \small\right\}}
 \end{array}
\end{equation}
where $\Gamma_T$ was defined in (\ref{SpaceCurves}). Moreover, the field of control sets given by
$
\mathbf{p} \mapsto \mathcal{B}_{\mathcal{F}^U}(\mathbf{p})
$ recall (\ref{CSet}) remains continuous when using $\mathcal{F}^U$ or $\mathcal{G}^U$ (instead of $\mathcal{F}$ or $\mathcal{G}$),
which directly follows by \cite[Lemma~6]{duits2018optimal} in conjunction with the important coercivity property (\ref{coerc}).

The nice thing in practice is that in-place rotations are automatically placed at optimal locations by the eikonal PDE system (solved by the anisotropic fast-marching algorithm that we discuss in the next section). It is not surprising that, when using a reasonable cost function $C$ (see Figs.~\ref{fig:CostFunctionTipBif} \& \ref{fig:CostFunctionBifSeed}), these in-place rotations are automatically placed at bifurcations in complex vascular trees as can be observed in the upcoming Figure~\ref{fig:TrackingPerTreeLIFvsDDLIF}.

\section{Numerical Solutions to the Eikonal PDE System: Extension of the Anisotropic Fast-Marching Algorithm that allows for all Left-Invariant Data-Driven Metrics}\label{sec:ExtensionAFM}


\def\bM{{\mathbb M}}
\def\bR{{\mathbb R}}
\def\bS{{\mathbb S}}
\def\bZ{{\mathbb Z}}
\def\cF{{\mathcal F}}
\def\cG{{\mathcal G}}
\def\cO{{\mathcal O}}
\def\gF{{\mathfrak F}}
\def\Id{\mathrm{Id}}
\def\<{\langle} \def\>{\rangle}
\def\ve{\varepsilon}
\def\cot{\hat}

In this section we describe the computation of globally minimizing geodesics for the new models $\mathcal{F}^U$ considered in this paper, whose fundamental ingredient is the numerical solution to an anisotropic eikonal PDE. 
The \emph{Reeds-Shepp forward} optimal control model $\mathcal{F}$, of which a variant $\mathcal{F}^U$ is considered in this paper, has already been addressed numerically using a Semi-Lagrangian \cite{duits2018optimal} and Eulerian \cite{mirebeau2018fast} discretization of the corresponding eikonal PDE. Both works however take advantage of the fact that the original geodesic model $\mathcal{F}$ regards the tangent spaces to the physical $\bR^2$ and the angular $S^1$ domains as orthogonal to each other. However, in our case of interest (with model $\mathcal{F}^U$ given by (\ref{Finslerfin})), we cannot expect such a block-matrix structure in
the fixed coordinate system $(x,y,\theta)$.

To overcome this problem, we describe below the extension of \cite{mirebeau2018fast} to the adaptive frame setting considered here, where this orthogonality relation is lost; in contrast, \cite{duits2018optimal} could not be generalized in an efficient manner. 

\begin{remark}[Convenience notations for the numerical section]
\label{rem:notations_numerics}
Throughout this section, we label the dependence w.r.t.\ the relaxation parameter $\ve \in (0,1]$,
so as to analyze it more easily, and to investigate the limit $\ve \to 0$. 
In contrast, we often drop the dependence w.r.t.\ the data $U$, which is regarded as fixed.

We also take advantage of the fact that the manifold $\bM_2 := \bR^2 \rtimes S^1 \equiv \bR^2 \times \bR/(2\pi \bZ)$ has a trivial tangent bundle: $T_{\p}(\mathbb{M}_{2}) \equiv \bR^2 \times \bR \equiv \bR^3$ canonically for any $\p \in \bM_2$, and likewise $T^*_\p(\mathbb{M}_{2}) \equiv \bR^3$.
As a result, by identifying co-vectors and vectors by their components in $\R^3$, the functional brackets $\langle \cdot, \cdot \rangle$ below boil down to the ordinary dot product on $\R^3$. Similarly, the tensor product $\omega\otimes \omega$  boils down to 
$\omega \omega^\top$. 
\end{remark}




\subsection{Asymmetric quadratic eikonal PDE}
\label{subsec:ExtensionAFM_PDE}

The Reeds-Shepp \emph{forward} model, is defined through a sub-Finslerian quasi-metric\footnote{I.e. the cost associated with a velocity at a point is non-Riemannian, non-symmetric, and possibly infinite.},
relaxed by a small parameter $\ve>0$, recall $\mathcal{F}$ was given by Eq.~(\ref{Finsler}) and its data-driven version $\mathcal{F}^U$ was given by  (\ref{Finslerfin}). 
Throughout this section, and in our vessel tracking experiments, we are concerned with the case where sideward motions and backward motions become very expensive: we set spatial anisotropy parameter $\zeta=\varepsilon$ with $0<\varepsilon\ll 1$ in the Finsler norm $\mathcal{F}^U$ given by (\ref{Finslerfin}).

The generic form of the data-driven Finsler metric function considered in this paper \eqref{Finslerfin} reads:
\begin{equation}
\label{eq:rsf_geodesic_metric}
\begin{array}{l}
	\cF_{\varepsilon}(\p,\dot \p)^2 = C(\p)^2\cdot\\[5pt]
	\left(\<\dot \p, M^0(\p) \dot \p\> + \epsilon^{-2}\<\omega^2(\p),\dot \p\>^2 + \epsilon^{-2} \<\omega^1(\p),\dot \p\>_-^2\right),
	\end{array}
\end{equation}
for any point $\p \in \overline \Omega$, within a given bounded connected domain $\Omega \subset \bM_2$, and any tangent vector $\dot \p \in T_{\p}(\mathbb{M}_{2}) \equiv \bR^3$, and where the two small parameters $\epsilon,\varepsilon$ relate via
\begin{equation}\label{eq:epsilon}
    \epsilon^{-2}:=(\ve^{-2}-1)\xi^2.
\end{equation}
In the following analysis, we only use the property that the tensor field $M^0$ is pointwise positive definite, that the differential forms $\omega^1$ and $\omega^2$ are pointwise linearly independent, 
and that $M^0 : \overline \Omega \to S_3^{++}$, and $\omega^1, \omega^2 : \overline \Omega \to \bR^3$ (following the conventions of Remark \ref{rem:notations_numerics}) have Lipschitz regularity. Here $S_{3}^{++}$ denotes the space of $3\times 3$ symmetric positive definite real-valued matrices.


\begin{remark}
In the normal left-invariant setting $\mathcal{F}^{U=1}=\mathcal{F}$
the asymmetric metric
expressed in the fixed frame $(\dot{x},\dot{y},\dot{\theta})$, of the tangent space at any coordinates $(x,y,\theta)$,
has a block diagonal structure. In contrast the data-driven metrics $\mathcal{G}^U$, $\mathcal{F}^U$, in general, does \emph{not} have a block-matrix structure, as the independent elements $\{\omega^i_U\}_{i=1}^3$  \emph{may point anywhere} due to their data-driven nature, as can be seen in Figure~\ref{fig:Frames}, keeping in mind the duality (\ref{duality}). Therefore, 
we must introduce a new modification of the anisotropic fast-marching algorithm. 
\end{remark}
The purpose of the second
term in (\ref{eq:rsf_geodesic_metric})
is to increase the cost of sideways motions, whereas the final  term prevents motions in reverse gear; both are excluded in the genuine Reeds-Shepp forward car model obtained in the limit (akin to \cite[Thm.2]{duits2018optimal}) as $\ve \to 0$, which is non-holonomic.

The distance map $W : \overline \Omega \to \bR$ from a given point $\p_0 \in \Omega$ and w.r.t.\ the Finsler function $\cF_\ve$, is the unique viscosity solution to the following anisotropic eikonal system 
\begin{equation}
\label{eq:rsf_eikonal}
\left\{
\begin{array}{l}
	\cF_\varepsilon^*(\p,{\rm d} W(\p)) = 1, \qquad \p \in \mathbb{M}_2\\
	W(\p_0)=0
	\end{array}
	\right.
\end{equation}
where the dual Finsler function equals by definition
\begin{equation*}
    \cF_\varepsilon^*(\p, \cot \p) := \sup\{\<\cot \p, \dot \p\>\;|\; \, \dot{\p}\in T_{\p}(\mathbb{M}_{2})\textrm{ and } \cF_\varepsilon(\p, \dot \p) \leq 1\},
\end{equation*}
with $\hat{\p}\in T_\p^*(\mathbb{M}_2)$. 

The structure of the metric \eqref{eq:rsf_eikonal}, referred to as \emph{asymmetric quadratic},
allows to compute a closed form expression of the dual metric $\cF_\epsilon^*$, and thereby the eikonal PDE \eqref{eq:rsf_eikonal}, as we will see below.

\begin{lemma}
\label{lem:asym_quad_dual}
	Let $M \in S_3^{++}$ and $\omega \in \bR^3$, and define 
\begin{equation*}
\label{eq:asym_quad}
	\cF(\p,\dot \p)^2 := \<\dot \p,M \dot \p\> + \<\omega,\dot \p\>_-^2.	
\end{equation*}
	Then $\cF$ is a quasi-norm (i.e.\ an asymmetric norm), whose dual quasi-norm reads for all $\cot \p \in  \bR^3$
\begin{equation}
\label{eq:asym_quad_dual}
	\cF^*(\p,\cot \p)^2 = \<\cot \p,D \cot \p\> + \<\cot \p,\eta\>^2_+,
\end{equation}
with $D := (M +\omega \, \omega^\top)^{-1}$ and $\eta := M^{-1} \omega/\sqrt{1+\omega^{\top} M^{-1}\omega}$.
\end{lemma}
\begin{proof}
See {\cite[Lemma 4]{duits2018optimal}}.
\end{proof}

For concreteness, we apply Lemma \ref{lem:asym_quad_dual} to our Finsler function $\cF_\epsilon$ of interest \eqref{eq:rsf_eikonal}, defined pointwise by the parameters
\begin{equation}
\label{eqdef:Meps}
    M_\epsilon := M^0+\epsilon^{-2} \omega^2 (\omega^2)^\top, \textrm{ and }
	\omega_\epsilon  := \epsilon^{-1}\omega^1.
\end{equation}
This then results in the following dual Finsler functions:
\begin{lemma}[Dual Finsler Functions]
\label{corol:rsf_eps_coefs} \\[4pt]
With our choice \eqref{eq:rsf_geodesic_metric} of Finsler function 
$\mathcal{F}_{\varepsilon}$ used in \eqref{eq:rsf_eikonal}, the dual Finsler function $\mathcal{F}_{\varepsilon}^*$ is given for all 
$\cot \p \in T_{\p}^*(\mathbb{M}_2) \equiv \bR^3$ by 
\begin{align}
\label{eq:D_eta_eps}
\cF_\varepsilon^*(\p, \cot \p)^2 &= \<\cot \p,D_\epsilon \cot \p\> + \<\cot \p,\eta_\epsilon\>^2_+, \textrm{ with } \\
	D_\epsilon & = 
	\frac{\cA \cA^\top}{\cA^\top M^0 \cA}
+ \cO(\epsilon^2), 
\label{eq:Depsilon}\\
	\eta_\epsilon &= 
	\frac{M^{-1}(\omega_U^1-\alpha \omega_U^2)}{\sqrt{ (\omega_U^1)^\top M^{-1}(\omega_U^1-\alpha \omega_U^2)}}
	+ \cO(\epsilon^2),\label{eq:etaepsilon}
\end{align}
where we used the shorthand notation $M^{-1} := (M^0)^{-1}$, the cross product $\cA := \omega_U^1 \times \omega_U^2$, and the orthogonalization coefficient $\alpha:=(\omega_U^2)^\top M^{-1} \omega_U^1 / (\omega_U^2)^\top M^{-1} \omega_U^2$.
\end{lemma}


\begin{proof}
Follows by Lemma~\ref{lem:asym_quad_dual} and Taylor expansion, for details, see Appendix~\ref{app:proofcor1}.
\end{proof}
Note that by \eqref{eq:epsilon}, $\mathcal{O}(\epsilon^2)=\mathcal{O}(\varepsilon^2)$ for small values of $\varepsilon$. 
%

Lemma \ref{corol:rsf_eps_coefs} shows that one can define an ideal sub-Finsler function $\cF_0^*$ that arises in the limiting case $\epsilon \downarrow 0$, and that 
\begin{equation}
\label{eq:FinslerFunctionApprox}
    \cF_\varepsilon^*(\p, \cot \p) = \cF_0^*(\p, \cot \p) + \cO(\epsilon^2).    
\end{equation}
Our goal, achieved in Sections \ref{subsec:ExtensionAFM_PDE} and \ref{subsec:ExtensionAFM_scheme} below, is to design a numerical scheme that is consistent with the sub-Finslerian eikonal PDE $\cF_0^*(\p, {\rm d} W(\p))=1$, and which satisfies the conditions that make the fast marching algorithm applicable.

\subsection{Discretization and consistency}
\label{subsec:ExtensionAFM_scheme}

We discretize the eikonal PDE \eqref{eq:rsf_eikonal}, which has an asymmetric quadratic structure \eqref{eq:asym_quad_dual}, using adaptive finite differences on the Cartesian grid $\Omega_h := \Omega \cap h \bZ^3$ of the domain $\Omega$, where $h>0$ denotes the grid scale. Note that $2\pi/h$ must be an integer, for consistency with the periodic boundary conditions in the angular coordinate. 
The numerical scheme construction relies on Selling's decomposition of positive definite matrices \cite{Selling:1874Algorithm} 
and on a corollary related to the approximation of the squared positive part of a linear form. The versions of these results presented in \cite[Corollary 4.12, Corollary 4.13]{mirebeau2018fast} are gathered in the following proposition. 

We denote $\mu(D) := \sqrt{\|D\| \|D^{-1}\|}$, for any $D \in S_3^{++}$, and $a_+ := \max \{0,a\}$, for any $a \in \bR$. 
\begin{proposition}[Selling matrix decomposition, 
see \cite{mirebeau2018fast}]  
\label{prop:decomp_D_v}
For any $D \in S_3^{++}$, there exists $\dot{\ul e}_1,\cdots,\dot{\ul e}_I \in \bZ^3$ and $\lambda_1,\cdots,\lambda_I \geq 0$, such that for all $\cot\p \in \bR^3$
\begin{equation*}
	\<\cot\p,D \cot\p\> = \sum_{1 \leq i \leq I} \lambda_i \<\cot\p,\dot{\ul e}_i\>^2. 
\end{equation*}
	For any $\eta \in \bR^3$, $\epsilon>0$, there exists $\dot{\ul f}_1,\cdots,\dot{\ul f}_I \in \bZ^3$ and $\mu_1,\cdots,\mu_I \geq 0$, such that for all $\cot\p \in \bR^3$
\begin{equation*}
	\<\cot\p,\eta\>_+^2 \leq \!
	\sum_{1 \leq i \leq I} \mu_i \<\cot\p,\dot{\ul f}_i\>_+^2 \leq \<\cot\p,\eta\>_+^2 + \epsilon^2 (\|\cot\p\|^2 \|\eta\|^2 - \<\cot\p,\eta\>^2).
\end{equation*}
In addition $\|\dot{\ul e}_i\|,\cdots,\|\dot{\ul e}_I\| \leq C \mu(D)$ and $\|\dot{\ul f}_i\|,\cdots,\|\dot{\ul f}_I\| \leq C\ve^{-1}$. The above holds with the constants $I := 6$, $C := 4\sqrt{3}$.
\end{proposition}
\begin{remark}
A key aspect of Proposition~\ref{prop:decomp_D_v} is that the vectors $(\dot{\ul e}_i)$ and $(\dot{\ul f}_j)$ have \emph{integer} coordinates, hence we avoid (off-grid) interpolations in our discretization. 
\end{remark}
Proposition \ref{prop:decomp_D_v} suggests the following discretization of \\ $\cF^*({\rm d}W(\p))$,  as in \eqref{eq:asym_quad_dual}, for arbitrary $D \in S_3^{++}$, $\eta \in \bR^3$, $\ve>0$:
\begin{equation} \label{gothicZFdef} 
\begin{array}{l}
	\sum \limits_{1 \leq i \leq I} \lambda_i \max \left\{0,\frac{W(\p)-W(\p-h \dot{\ul e}_i)} h,\frac{W(\p)-W(\p+h\dot{\ul e}_i)}h\right\}^2 \\
	+ \sum \limits_{1 \leq j \leq I} \mu_j \max\left\{0,\frac{W(\p)-W(\p-h\dot{\ul f}_i)}h\right\}^2 =: \gF W(\p),
	\end{array}
\end{equation}
with suitable boundary conditions.
This numerical scheme falls within the Hamiltonian fast-marching framework \cite{mirebeau2019hamiltonian}, and thus can be efficiently solved numerically, see Section \ref{subsec:fmm}. Before that, we discuss its consistency with the eikonal equation: 
inserting a first order Taylor expansion in \eqref{gothicZFdef}, we obtain (using Proposition \ref{prop:decomp_D_v}):
\begin{equation}
	\label{eq:asym_quad_consistency}
	\begin{array}{ll}
	\gF W(\p) &=
 \< \nabla W(\p), D \nabla W(\p)\> + \<\nabla W(\p),\dot \eta\>_+^2 \\
 
	&\qquad + \cO( R h + \epsilon^2),
	\end{array}
\end{equation}
where $R := \max\{\mu(D),\epsilon^{-1}\}>0$ denotes the stencil radius. 


We next analyze the approximation error towards the ideal model as $\epsilon \to 0$ and $h\to 0$ suitably. 
To this end we denote by $\gF^h_\epsilon$ the finite differences scheme on $\Omega_h$ associated with the parameters $D_\epsilon$ and $\eta_\epsilon$ of our application \eqref{eq:D_eta_eps}. 
Note that the stencil radius is $R_\epsilon=\max\{\mu(D_{\epsilon}),\epsilon^{-1}\} = \cO(\epsilon^{-1})$, since $\mu(D_\epsilon) = \mu(M_\epsilon) = \cO(\epsilon^{-1})$ in view of \eqref{eqdef:Meps}.  Now combining \eqref{eq:FinslerFunctionApprox} with  \eqref{eq:asym_quad_consistency}, we obtain the overall consistency error 
\begin{equation}
\label{eq:scheme_consistent_0}
	\begin{cases}
	    \gF^h_\epsilon W(\p) = \cF_0^*(\p, {\rm d} W(\p))+ \cO(\epsilon^{-1} h + \epsilon^2)=1,\\
	    \gF_\epsilon^h W(\p_0)=0.
	\end{cases}
\end{equation}
The optimal order of the consistency error $\cO(h^\frac 2 3)$ is achieved with the scaling $\epsilon = h^{\frac 1 3}$. 
In our practical experiments however, we rather use the small fixed value $\zeta=\ve=0.1$ which by \eqref{eq:scheme_consistent_0} 
indeed yields a sufficiently accurate scheme \cite{mirebeau2018fast}!

\subsection{Anisotropic Fast-Marching for Computing Distance Maps of Data-driven Left-invariant Finsler Models}
\label{subsec:fmm}

In fast-marching methods (FMM), one divides the grid points into 3 categories: Far, Trial, and Accepted. 
In each step of the FMM, one selects the trial point $\p$ whose function value $W(\p)$ is minimal. The point $\p$ is moved into the accepted set, and $W(\p)$ is frozen. In contrast, all the trial or far points whose stencil contains $\p$ see their function value \emph{updated}, and they are moved into the trial set. 
This procedure generalizes and abstracts the classical FMM \cite{sethian1999fast}, for details see \cite{PhDDaChen}.
When all points have moved to the accepted category, the FMM terminates, and a geodesic can be easily calculated by steepest descent as described in Section \ref{sec:SteepestDescentFinslerianGeodesics}.

The update of a single function value $W(\p)$ 
is defined as follows: isolate $W(\p)$ in the numerical scheme expression \eqref{gothicZFdef}, and express it by the values of its neighbors so as to obey $\gothic{F}W(\ul{p})=1$. The latter equation admits by construction a unique solution, which is obtained as the largest root of a quadratic equation. 


There are two crucial properties of the discretization $\gF$:
\begin{itemize}
	\item \underline{Discrete Degenerate ellipticity}:\\[3pt] $\gF W(\p)$ is a non-decreasing function of the finite differences 
	in (\ref{gothicZFdef}).
	\item \underline{Causality}: \\[3pt] $\gF W(\p)$ only depends on the positive part of all finite differences in (\ref{gothicZFdef}).
\end{itemize}
These are the two key\footnote{
There are two other minor assumptions, the existence of a sub- and super-solution to the scheme, and a perturbation property, which can be established 
analogously to the Riemannian case in \cite[Proposition 2.5]{mirebeau2019riemannian}.
} 
assumptions of \cite[Theorem 2.3]{mirebeau2019riemannian}, implying that \emph{the discretized PDE \eqref{eq:scheme_consistent_0} admits a unique solution $W_h^\epsilon$, which is computable in a single pass over the discretization domain} $\Omega_h$, using 
anisotropic fast-marching.

Following the steps of the proof \cite{mirebeau2018fast} associated with the standard Reeds-Shepp forward model, 
and with straightforward adaptations, we obtain that $W_h^\epsilon$ converges uniformly as $\epsilon \to 0$ and $h/\epsilon \to 0$ to the solution $W$ of the sub-Finslerian Eikonal PDE $\cF_0^*(\p, {\rm d} W(\p))=1$. 
\\
\\










\subsection{Steepest Descent for the Finslerian Geodesics}
\label{sec:SteepestDescentFinslerianGeodesics}

In previous work \cite[Thm.4]{duits2018optimal}, standard Riemannian steepest descent tracking on the distance map $W$, recall (\ref{BT}) in Theorem~\ref{th:th}, is generalized from the symmetric Riemannian setting to the (possibly asymmetric) Finsler model setting. That idea also transfers to the new data-driven left-invariant model as the Finslerian back-tracking 
\cite[App.B]{duits2018optimal} still applies. 
Steepest descent tracking  (\ref{BT}) from $\p$ to source point $\p_0$ again becomes (using the embedding of $\mathbb{M}_2\subset\mathbb{R}^3$)
\begin{align}
\gamma(t)=\p -\frac{1}{W(\p)}\int \limits_{0}^{t}
{\rm d}\mathcal{F}^{*}_{\varepsilon}(\, \gamma(s)\,,\, {\rm d}W(\gamma(s))\,)\, {\rm d}s, \textrm{ \ }t\in [0,1],\label{eq:BacktrackingSteepestDescent}
\end{align}
with $0<\epsilon \ll 1$ and
where the derivative is taken with respect to the second entry of the dual Finsler function, and where $W$ is the viscosity solution of the 
eikonal PDE system (\ref{eq:rsf_eikonal}).  

In the practical implementations we use a second order Runge-Kutta method for time integration approximations and at time $t=1$ we arrive at the source-point $\p_0$. 


This geodesic backtracking in $(\mathbb{M}_{2},\mathcal{F}^U)$ again boils down to piecewise concatenations of  
\begin{enumerate}
\item\label{it:1} symmetric Riemannian geodesics in $\mathbb{M}_{2}$
with metric tensor field $\mathcal{G}^U_{\p}(\dot{\p},\dot{\p})$ recall (\ref{Finslerfin}), and
\item symmetric Riemannian geodesics in $\mathbb{M}_2$ with metric tensor field $\mathcal{G}^U_{\p}(\dot{\p},\dot{\p})+(\epsilon^{-2}-1)|\omega^{1}_U(\dot{\p})|^2$ that are in-place rotations, at locations where $\mathcal{A}_1^U \approx \mathcal{A}_{1}$ if $0<\varepsilon \ll 1$. 
\end{enumerate}
\begin{remark} \label{rem:niceinpractice}
    Taking the negative part of the final term in (\ref{eq:rsf_geodesic_metric}) means that we switch
    between two Riemannian manifolds (one with the final term active and with the final term non-active). 
    At locations where $\omega^{3}_U \approx \omega^3$, for instance at locations where
    $\mathcal{A}_1^U\approx \mathcal{A}_1$ this means that
    one switches between anisotropic geodesics and spherical geodesics (in-place rotations). In the vessel tracking applications we want such in-place rotations to appear above bifurcations in the vasculature.
\end{remark}
A closely related situation is discussed in \cite[Thm.4]{duits2018optimal}, but now it is applied to the new data-driven model $\mathcal{F}^U$ (\ref{Finslerfin}) with dual 
$(\mathcal{F}^U)^*=\lim \limits_{\epsilon \downarrow 0} \mathcal{F}_{\epsilon}^*$. 

By Theorem~\ref{th:th} the Riemannian geodesics have parallel momentum
and their transitions
 towards spherical geodesics is like $C$ continuously  differentiable.
The surface $\mathcal{S}$ where Finslerian geodesics of $\mathcal{F}^U$ in $\mathbb{M}_{2}$ switch from one Riemannian manifold to the other is given by
\[
\mathcal{S}=
\{\p \in \mathbb{M}_{2}\;|\; \omega^{1}_U(\nabla W(\p))= 
(\mathcal{A}_{1}^U W)(\p)=0\}.
\]
Interestingly, in the mixed-model $\mathcal{F}^M$ that we will explain later in Section~\ref{sec:experiments}, the condition in Remark~\ref{rem:niceinpractice} above is satisfied at bifurcations. Then in-place rotations are indeed automatically placed at the bifurcations in the tracking of blood vessels, as can be seen in Figure~\ref{fig:TrackingPerTreeLIFvsDDLIF}. 

\section{Experiments}\label{sec:experiments}
We choose the data-driven left-invariant metric tensor field with forward gear $\mathcal{F}^U$ as given in Eq.~\eqref{Finslerfin}.  An elaboration on the calculation of the cost function can be found in Appendix~\ref{app:CostFunction}. We will discuss the new model's ability to adapt for curvature. Additionally, we show and discuss some full vascular tree tracking results.





\subsection{Curvature Adaptation}
\begin{figure}
    \centering
    \includegraphics[width=0.48\textwidth]{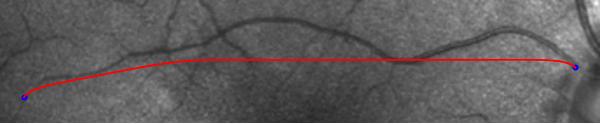}
    \includegraphics[width=0.48\textwidth]{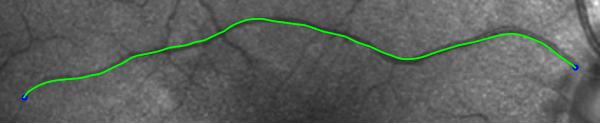}
    \caption{ {\bf Influence of \emph{data-driven} metric tensor fields:} (top) Tracking with the vanilla left-invariant metric tensor field from \eqref{eq:LIMTFDD}. (bottom) Tracking with the proposed \emph{data-driven} left-invariant metric tensor field from \eqref{Finslerfin}. To isolate its effect in the tracking process and record the effect of only directional adaptation of the underlying metric, we have set the cost function $C=1$. We observe that the \emph{data-driven} nature of our model allows for a better fidelity to the underlying vascular structure. The parameters are given by $g_{11}=0.01,g_{22}=1,g_{33}=1,\lambda=100$. }
    \label{fig:TrackingNoCost}
\end{figure}

The data-driven left-invariant metric tensor field $\mathcal{G}^U$ and its asymmetric variant $\mathcal{F}^U$ both consists of a ``standard'' left-invariant metric tensor field to which a data-driven term is added, as introduced in \eqref{Finslerfin}. The idea of the second data-driven term in this equation is that it pushes the main direction of the model into the direction of the underlying vessel, as illustrated in Fig.~\ref{fig:TrackingNoCost}, where no data-dependent cost function $C=1$ was used to generate the tracking result. We see that the data-driven left-invariant metric tensor field takes the image data into account and steers the tracking in the direction of the underlying vascular structure, even when the cost function does not contain information about vessel locations and curvature.

\begin{figure}
    \centering
    \begin{subfigure}[b]{0.48\textwidth}
        \centering
        \includegraphics[width=0.48\textwidth]{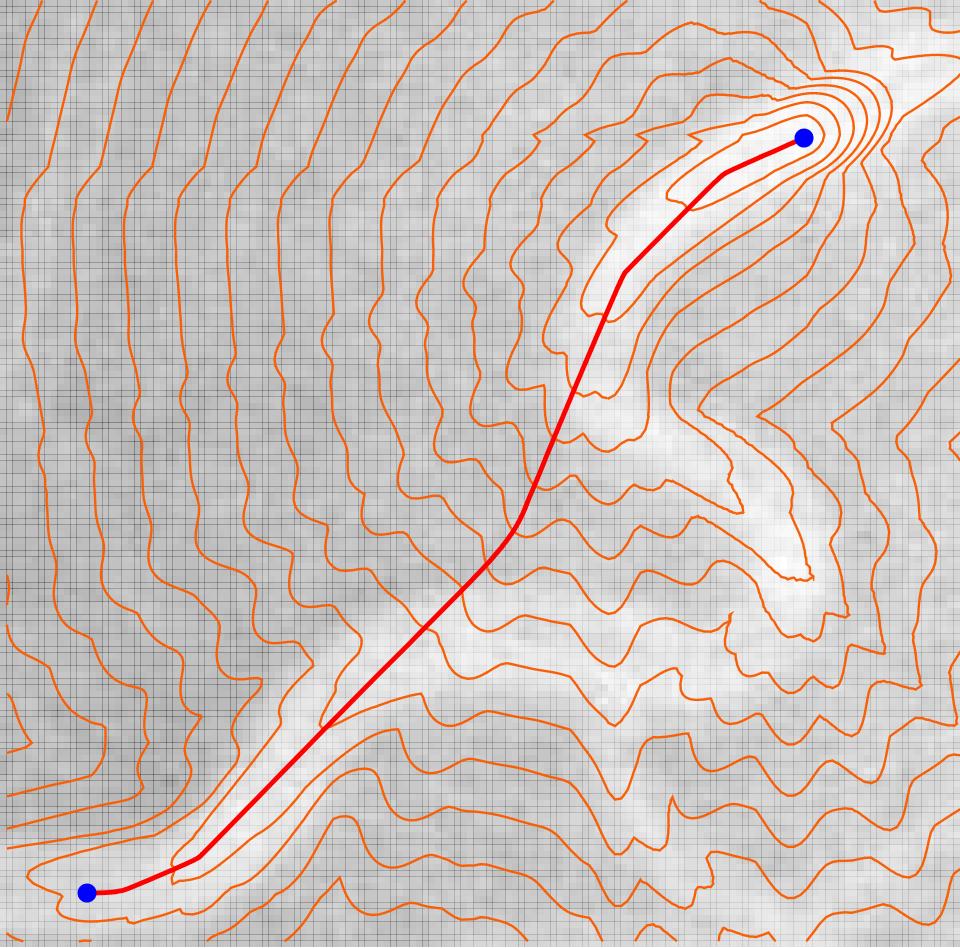}
        \hfill
        \includegraphics[width=0.48\textwidth]{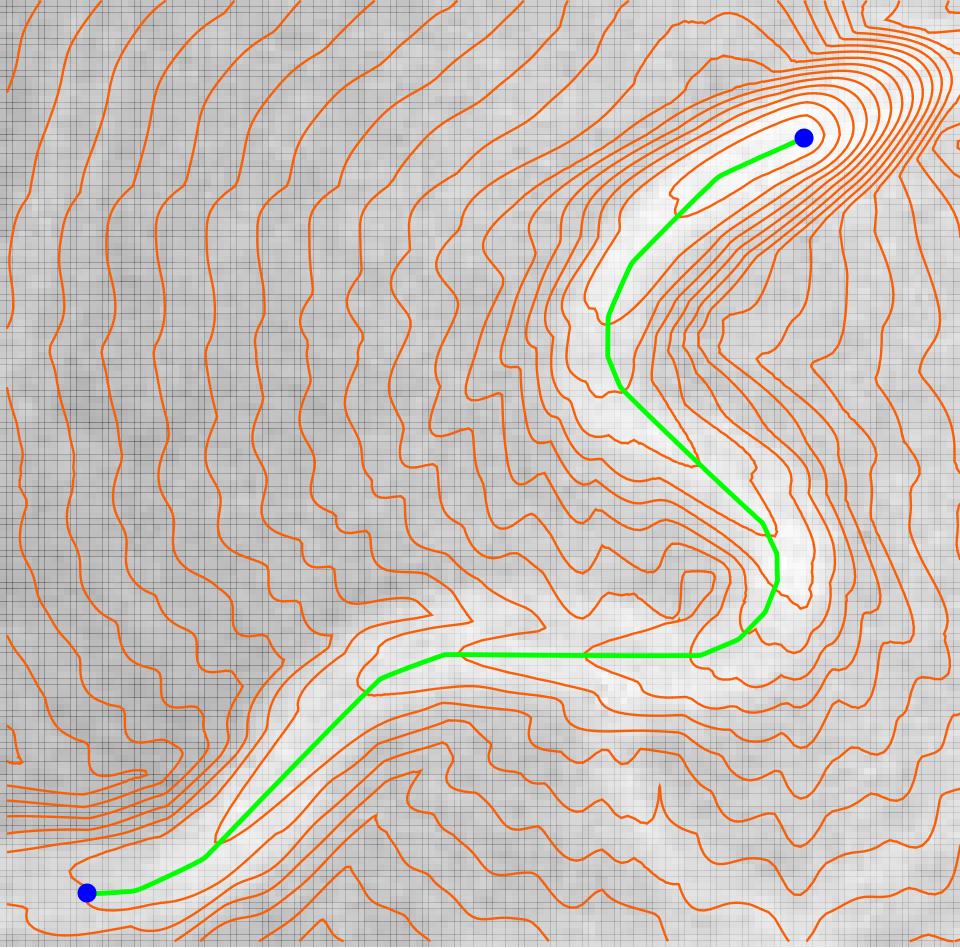}
        \caption{Tracking results using 16 orientations}
        \label{fig:S-curve16orientations}
    \end{subfigure}
    \begin{subfigure}[b]{0.48\textwidth}
        \centering
        \includegraphics[width=0.48\textwidth]{S-curveLIFFGbeta5zeta0125eps01trackingDistances.png}
        \hfill
        \includegraphics[width=0.48\textwidth]{S-curveGFFGbeta5alpha08eps01trackingDistances.png}
        \caption{Tracking results using only 8 orientations}
        \label{fig:S-curve8orientations}
    \end{subfigure}
    \caption{\textbf{Comparison tracking results of left-invariant and data-driven left-invariant metric tensor field:} Tracking results for left-invariant metric tensor field $\mathcal{F}$ (left) and data-driven left-invariant metric tensor field $\mathcal{F}^U$ (right). The parameters for the (data-driven) left-invariant metric tensor field are given by $g_{11}=0.01,g_{22}=0.16,g_{33}=1,\lambda=100$. The cost function is given by $C=1/(1+ 200 
    |U_f|^2)$. We see that the iso-contours of the data-driven metric tensor field follow the vessel structure better, and the tracking is correct (even with 8 orientations).}
    \label{fig:S-curve}
\end{figure}

\begin{figure}
    \centering
    \begin{subfigure}[b]{0.45\textwidth}
        \centering
        \includegraphics[width=\textwidth]{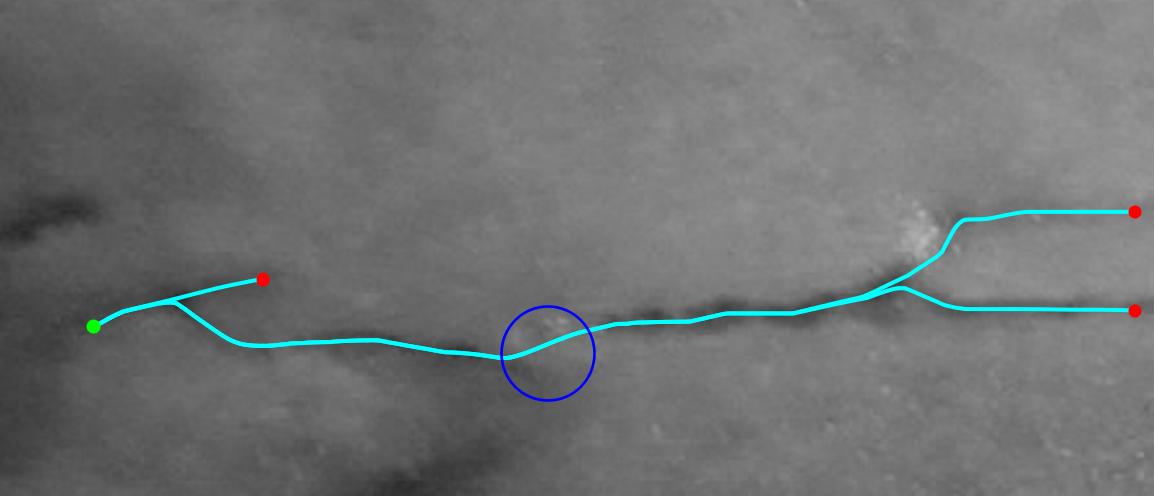}
        \caption{Tracking results with the left-invariant metric tensor field.}
    \end{subfigure}
    \begin{subfigure}[b]{0.45\textwidth}
        \centering
        \includegraphics[width=\textwidth]{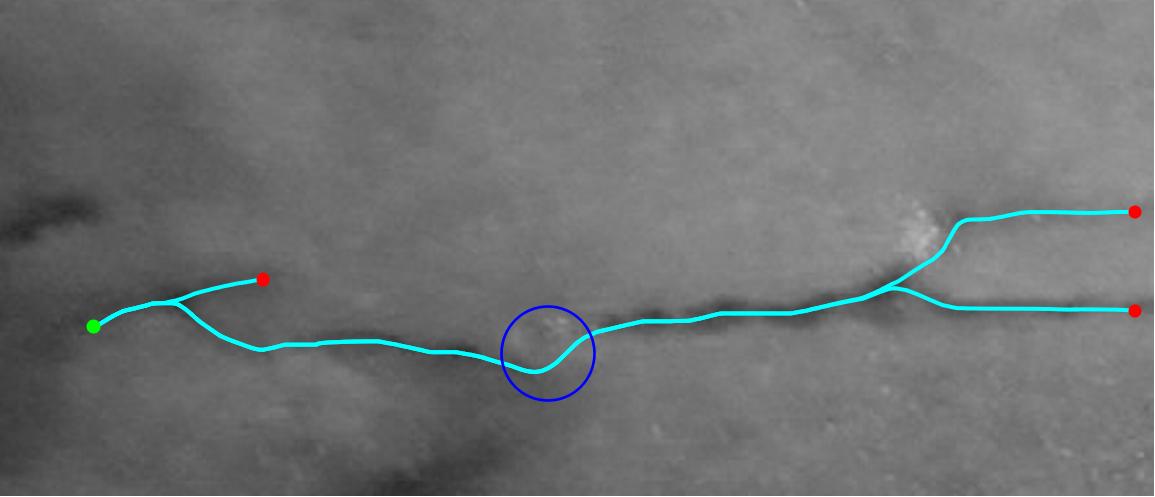}
        \caption{Tracking results with the data-driven left-invariant metric tensor field with $\lambda=100$.}
    \end{subfigure}
    \caption{\textbf{Application in crack detection:} Tracking results for the left-invariant and data-driven left-invariant metric tensor field on an image of cracks in a building. The presented results are calculated using 32 orientations, and parameter settings $g_{11}=0.01$, $g_{22}=1$, $g_{33}=1$. In regions with high curvature, the data-driven model adapts more gradually for curvature to get more data evidence than the left-invariant model which tends to prefer in-place rotations. }
    \label{fig:ImageAndrii}
\end{figure}

The data-driven term also leads to better tracking results of very tortuous vascular structures as we see in Fig.~\ref{fig:S-curve}. In Fig.~\ref{fig:S-curve16orientations}, the tracking results relying on (solely) the left-invariant metric tensor field $\mathcal{F}$ fail to follow the underlying vessel correctly, contrary to the new data-driven left-invariant model $\mathcal{F}^{U}$ (\ref{Finslerfin}) which follows the vascular structure correctly. In addition, one sees that when using the left-invariant model, the wave fronts (indicated in orange) suffer from the discretization. In the data-driven left-invariant model, these discretization artifacts are gone, and the wavefronts follow the underlying vascular structure correctly. When only considering 8 orientations, as in Fig.~\ref{fig:S-curve8orientations}, the data-driven left-invariant frame is still able to follow the vascular structure correctly. Although the discretization is clearly visible in the tracking results, the data-driven left-invariant metric tensor field is still able to follow the vessel correctly. It is important to note that both models use the same cost function. The differences in the tracking results are due to the data-driven left-invariant metric tensor fields' ability to  better adapt for:
\begin{enumerate} 
    \item gradual curvature change of blood vessels. The same applies to other applications such as the detecting of cracks, see Figure~\ref{fig:ImageAndrii},
    \item orientation biases by orientation score data- alignment as depicted in Figure~\ref{fig:Frames}. 
\end{enumerate}

\subsection{Complete Vascular Tree Tracking}

In the previous section, we discussed the curvature adaption feature of the new (asymmetric) data-driven left-invariant metric tensor field $\mathcal{F}^U$. This model also can automatically place `in-place' rotations in globally optimal geodesics which are typically placed at bifurcations.

\begin{figure}
    \centering
    \begin{subfigure}[b]{0.45\textwidth}
    \centering
        \includegraphics[width=\textwidth]{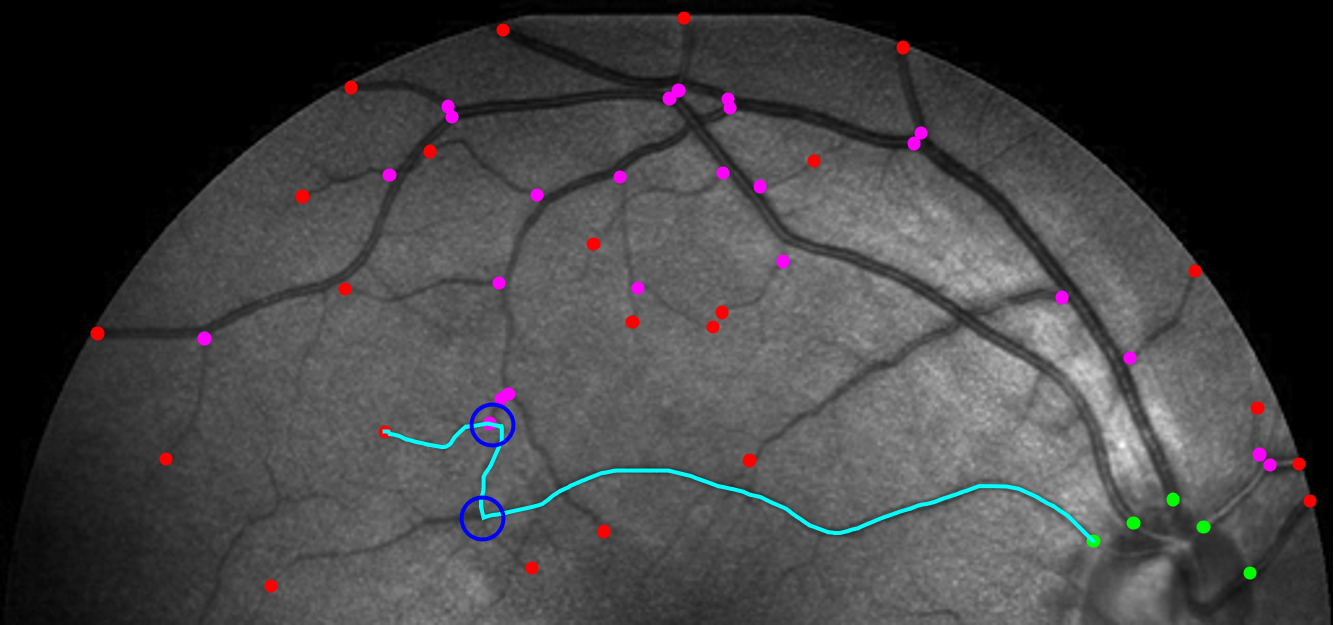}
        \caption{\textbf{Challenging case with wrong tracking result} via model $\mathcal{F}^U$.}
    \end{subfigure}
    \begin{subfigure}[b]{0.35\textwidth}
    \centering
        \includegraphics[width=.25\textwidth]{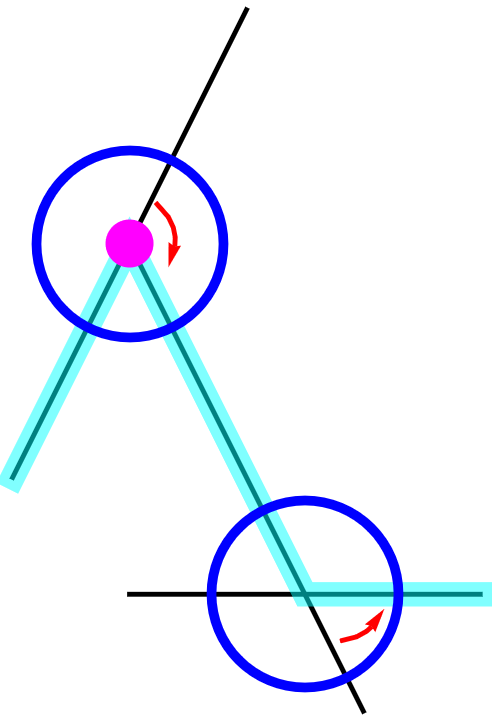}
        \caption{\textbf{Schematic representation:} Two vascular structures are visualized (black). One of them has a bifurcation point (purple). The tracking result (light blue) moves in the wrong direction at the bifurcation, after which it switches vascular structure at the crossing (rotation directions indicated in red).}
        \label{fig:SchemeticRepresentationRisk}
    \end{subfigure}
    \begin{subfigure}[b]{0.45\textwidth}
    \centering
        \includegraphics[width=\textwidth]{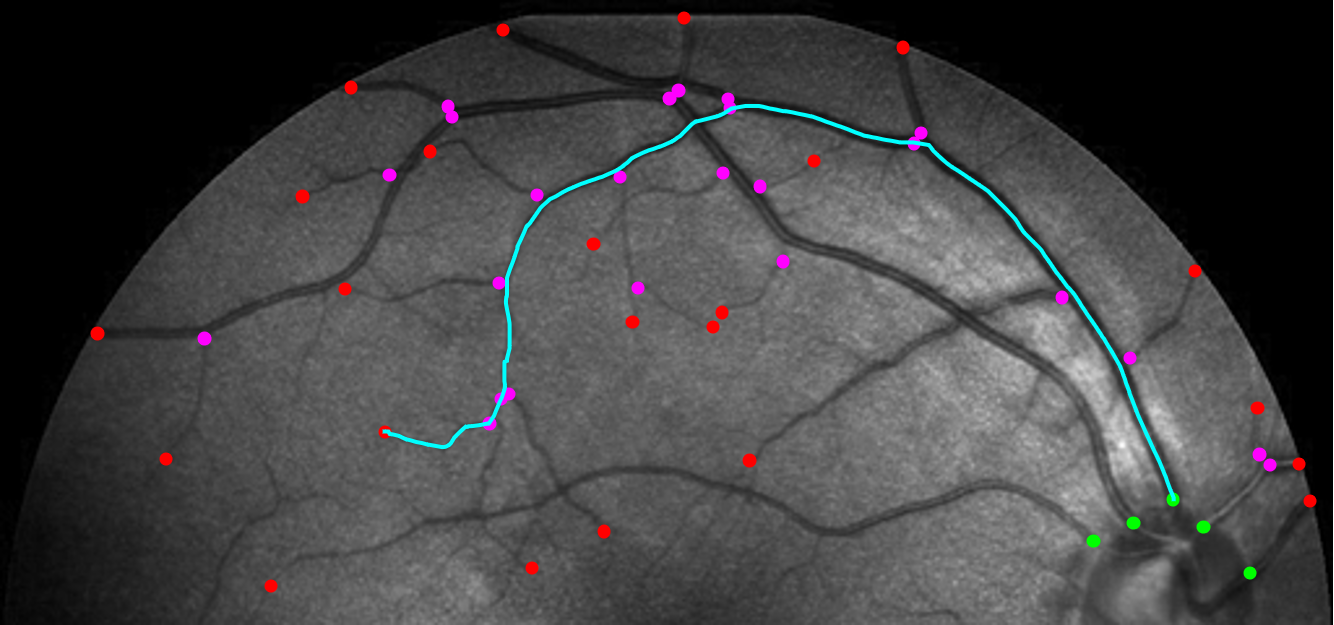}
        \caption{\textbf{Correct tracking result} via model $\mathcal{F}^M$}
        \label{fig:TortuousVesselTrackingMM}
    \end{subfigure}
    \caption{\textbf{Motivation mixed model:} Too much curvature adaptation at crossings is dangerous in extreme cases. The mixed model, introduced in Eq.~\eqref{eq:MixedLIMTF}, is preferable as it only adapts for curvature (like in Fig.~\ref{fig:S-curve}) \emph{in between} those complex structures, and indeed provides correct tracts everywhere, as can be seen in Fig.~\ref{fig:TortuousVesselTrackingMM}. The geodesics of both models are computed using $\lambda=100$. 
    }
    \label{fig:TrackingWrongAtCrossings}
\end{figure}

However, these valuable abilities of the model may also lead to extremely complex structures to some undesirable behavior. In full vascular tree tracking, we see that the data-driven term may steer the tracking away from the actual vessel at crossings in extreme cases, as can be seen in Fig.~\ref{fig:TrackingWrongAtCrossings}.

To overcome this problem (see item~c in Fig.~\ref{fig:TrackingWrongAtCrossings}), and to still take advantage of improved data-adaptation (like in Figs.\ref{fig:TortuousVesselTracking}, \ref{fig:TrackingNoCost} and \ref{fig:S-curve}) we introduce a (new) mixed model that prevents this undesirable behavior at (extreme) complex structures, where it locally relies on the old model. Then in between (extreme) complex structures we still benefit from the directional data adaptation in the orientation score.  

The mixed metric tensor field $\mathcal{G}^M$ (and its asymmetric version $\mathcal{F}^M$) is given by the data-driven left-invariant metric tensor field away from the crossing structures, and by the left-invariant metric tensor field otherwise:
\begin{align*}
    &\mathcal{G}_\p^{M}(\dot{\p},\dot{\p})=\kappa(\x)\,\mathcal{G}_\p(\dot{\p},\dot{\p})+(1-\kappa(\x))\; \mathcal{G}_\p^{U}(\dot{\p},\dot{\p}),\numberthis\label{eq:MixedLIMTF} \\ 
     &\mathcal{F}^{M}(\p,\dot{\p})^2=\kappa(\x)\, \mathcal{F}(\p,\dot{\p})^2+(1-\kappa(\x))\, \mathcal{F}^{U}(\p,\dot{\p})^2, 
\end{align*}
for all $\p=(\x,\n) \in \mathbb{M}_2$, and all $\cot{\mathbf{p}}=(\dot{\x},\dot{\n})  \in T_{\p}(\mathbb{M}_2)$, and
where $\kappa(\x)= \mathbbm{1}_A * G_{\sigma}(\x)$ with $A=\cup_{i=1}^N[\x_i-a,\x_i+a]$, where $\x_i$ representing $N \in \mathbb{N}$ crossing locations in the image. 

The results are typically not sensitive to the choice of $a$ and $\sigma$ in our application as long as $a>2$. Therefore we always set $a=5$ and $\sigma= 1$ pixel-size in our experiments. 
\\This construction of the metric tensor field ensures that the metric tensor field is not tempted to move in the wrong direction in extreme cases where vessels cross each other. The tracking result relying on the mixed metric tensor field is visualized in Fig.~\ref{fig:TortuousVesselTrackingMM}, and does not show the earlier mentioned undesirable behavior, as shown in Fig.~\ref{fig:SchemeticRepresentationRisk}. Therefore, this new model will be used in all full vascular tree tracking results. All results presented in this section are calculated using  parameters 
$g_{11}=0.01, g_{22}=g_{33}=1$. For the curve optimisation this is the same as setting $\varepsilon=\zeta=0.1$ in (\ref{ARSCM}) used in (\ref{Finslerfin}). Even for such extreme anisotropy settings, our numerical algorithm is appropriate as motivated in Section~\ref{sec:ExtensionAFM}. 
We always observed that tracking is stable with respect to small variations in these parameters, so there was no point in fine-tuning them.

\subsubsection{Asymmetric Double Step}
\begin{figure}
    \centering
    \begin{subfigure}[b]{0.45\textwidth}
    \centering
        \includegraphics[width=\textwidth]{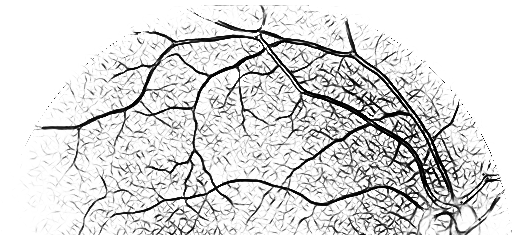}
        \caption{Image 1}
    \end{subfigure}
    \begin{subfigure}[b]{0.45\textwidth}
    \centering
        \includegraphics[width=\textwidth]{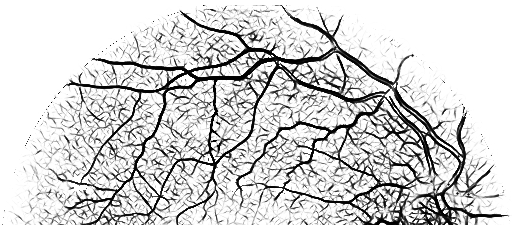}
        \caption{Image 2}
    \end{subfigure}
    \caption{{\bf Projected cost functions for tracking in two steps  - from bifurcation to tip:} Cost used to connect tips to the nearest bifurcation. Black and white mean low and high costs respectively. This cost function supports movement along the thin vessels very well. The multiscale vesselness is computed as explained in App.~\ref{app:CostFunction}, and the considered spatial scales is $\sigma_s=1$.
    }
    \label{fig:CostFunctionTipBif}
\end{figure}

\begin{figure}
    \centering
    \begin{subfigure}[b]{0.49\textwidth}
    \centering
        \includegraphics[width=\textwidth]{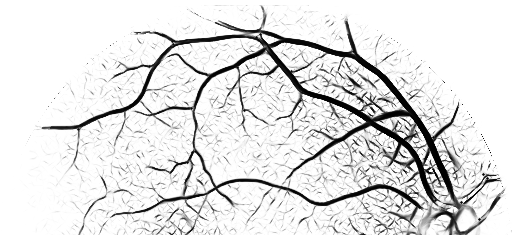}
        \caption{Image 1}
    \end{subfigure}
    \begin{subfigure}[b]{0.49\textwidth}
    \centering
        \includegraphics[width=\textwidth]{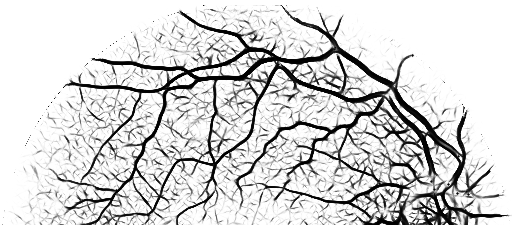}
        \caption{Image 2}
    \end{subfigure}
    \caption{{\bf Projected cost functions for tracking in two steps - from seed to bifurcation:} Cost used to connect bifurcations to optic nerve. Black and white mean low and high costs respectively. This cost function supports movement along the thick vessels very well. The multiscale vesselness is computed as explained in App.~\ref{app:CostFunction}, and the considered spatial scales are $\sigma_s\in\{1,2\}$.}
    \label{fig:CostFunctionBifSeed}
\end{figure}

\begin{figure}
    \centering
    \begin{subfigure}[b]{0.45\textwidth}
    \centering
    \includegraphics[width=\textwidth]{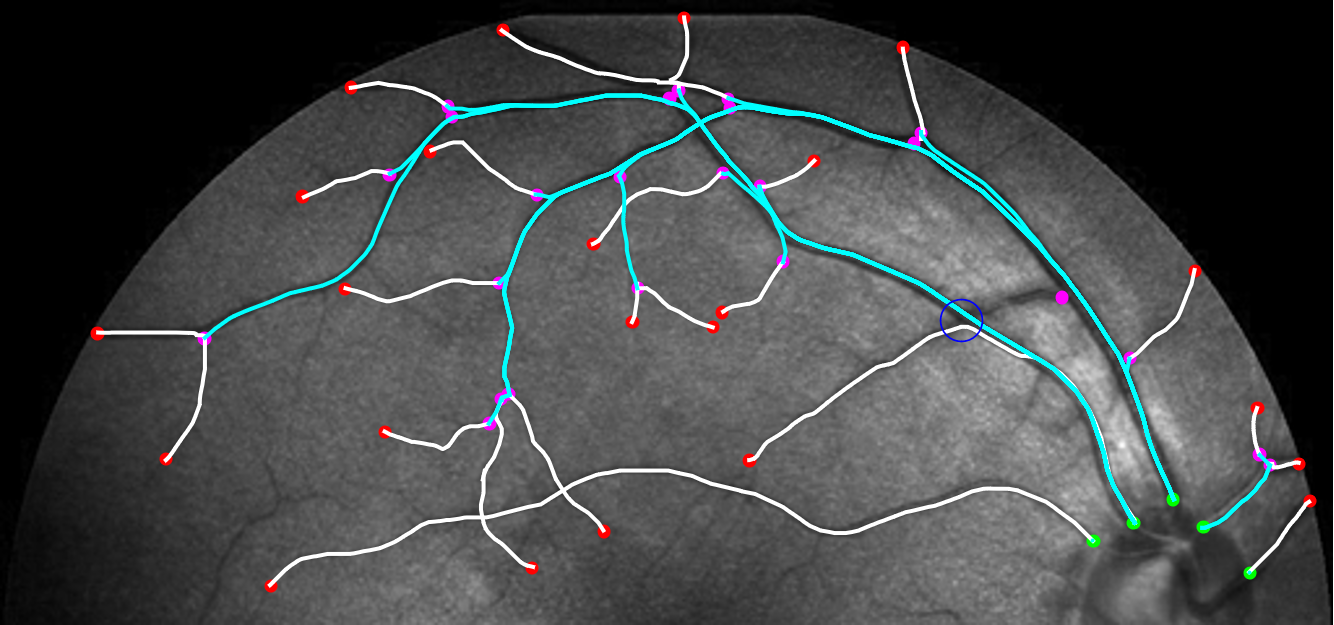}
        \caption{Image 1: Tracking with the mixed model with $\lambda=20$.}
        \label{fig:2StepMethodImage1}
    \end{subfigure}
    \begin{subfigure}[b]{0.45\textwidth}
    \centering
    \includegraphics[width=\textwidth]{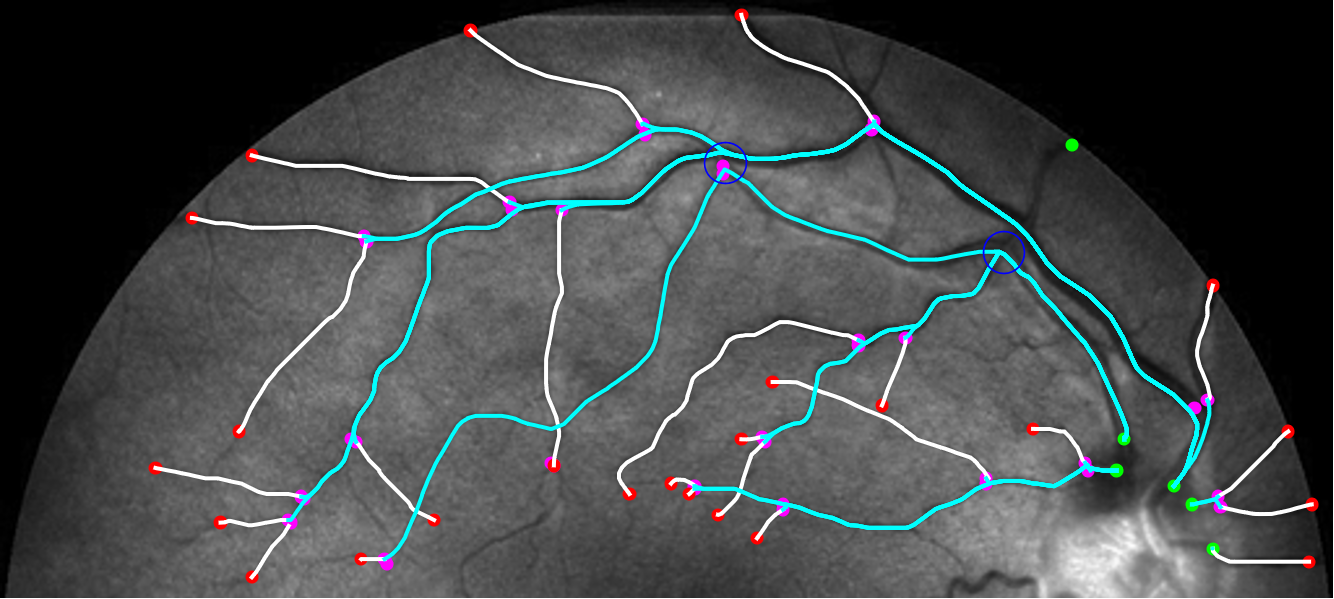}
        \caption{Image 2: Tracking with the mixed model with $\lambda=15$.}
    \end{subfigure}
    \caption{{\bf Two step tracking of Vascular Tree structures:} Tracking with mixed model $(\mathbb{M}_2,\mathcal{F}^M)$ proposed in \eqref{eq:MixedLIMTF}.  
    The first step involves connecting the tips (marked in red) to the nearest bifurcation (marked in purple) using the cost function depicted in Fig.~\ref{fig:CostFunctionTipBif}. Second, these bifurcations are now tracked to the seeds (marked in green) using the cost function depicted in Fig.~\ref{fig:CostFunctionBifSeed}.} 
    \label{fig:AsymmetricDoubleStepLIFvsDDLIF}
\end{figure}

The tracking results were computed in two steps; first tips are connected to the nearest bifurcation/seed (cost function visualized in Fig.~\ref{fig:CostFunctionTipBif}), after which those points are connected to the nearest seed (cost function visualized in Fig.~\ref{fig:CostFunctionBifSeed}). The used cost functions (from tip to nearest bifurcation and from bifurcation to seed) support movement along the thin and thick vessels respectively. The tracking results that correspond to these cost functions are depicted in Fig.~\ref{fig:AsymmetricDoubleStepLIFvsDDLIF}. The calculated geodesics are all correct, except for 2 difficulties when: 
\begin{enumerate}
\item crossings and bifurcations are very close to 
each other, 
\item vascular structures are kissing. \end{enumerate}

\begin{figure}
    \centering
    \begin{subfigure}[b]{0.45\textwidth}
    \centering
    \includegraphics[width=\textwidth]{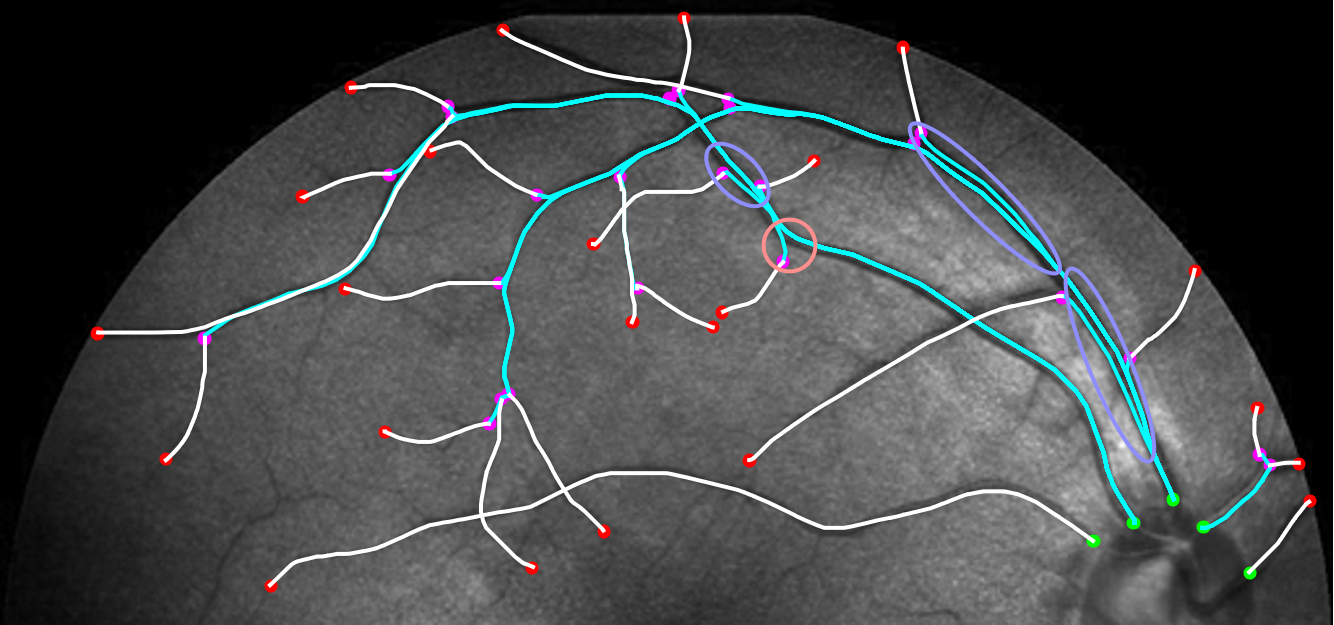}
        \caption{Tracking with the left-invariant model $\mathcal{F}$.} 
        \label{fig:LIFTracking}
    \end{subfigure}
    \begin{subfigure}[b]{0.45\textwidth}
    \centering
    \includegraphics[width=\textwidth]{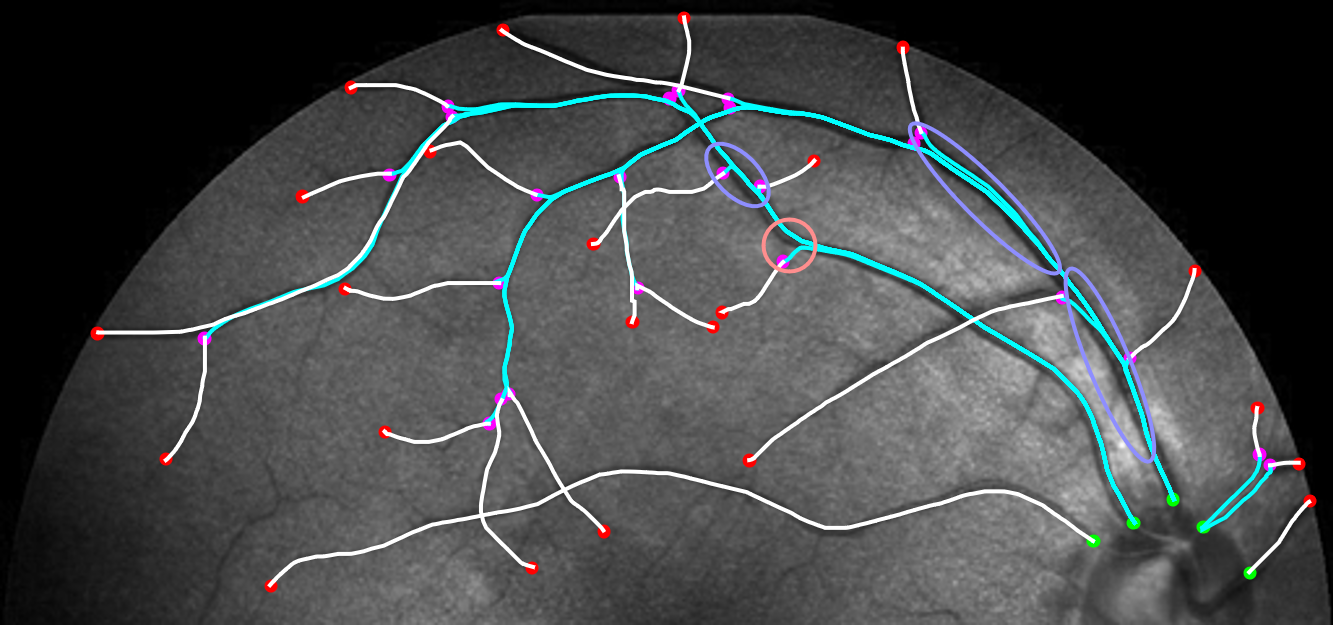}
        \caption{Tracking with the mixed model $\mathcal{F}^M$ proposed in \eqref{eq:MixedLIMTF} with $\lambda=65$.}
    \end{subfigure}
    \caption{{\bf Two-step tracking of vascular tree structures:}   
    The first step involves connecting the tips (marked in red) to the nearest bifurcation (marked in purple). Second, these bifurcations are now tracked back to the seeds (marked in green). The cost functions used in the first (white) and second (blue) step are  given in \cite{duits2018optimal}, with $\sigma=800$ and $p=4$, and Fig.~\ref{fig:CostFunctionBifSeed} respectively. The main differences between both models indicate that the mixed model $\mathcal{F}^M$ follows the vascular structure better and is better able to follow the centerline of a given vascular tree (marked by pink and blue circles respectively).} 
    \label{fig:LIFvsMM2Step}
\end{figure}

Next, we compare the results of the new mixed model $\mathcal{F}^M$ and the left-invariant model $\mathcal{F}$ in Fig.~\ref{fig:LIFvsMM2Step}. There are some visible differences between both tracking methods, marked in pink and blue. First, we see that the tracking results relying on the mixed method ensure that the centerline is better followed, and multiple geodesics are at approximately the same place (in blue). Second, we see the ability of the mixed method to adapt to the direction of the vascular structure (in pink). Instead of moving towards a bifurcation point away from the seed in the left-invariant metric tensor field, the curvature adaptation ensures that the tracking results immediately move towards the seed it is connected to.

\subsubsection{Asymmetric Single Step with Prior Classification of Seeds and Tips}
\begin{figure*}
    \centering
    \begin{subfigure}[b]{\textwidth}
    \centering
    \includegraphics[width=\textwidth]{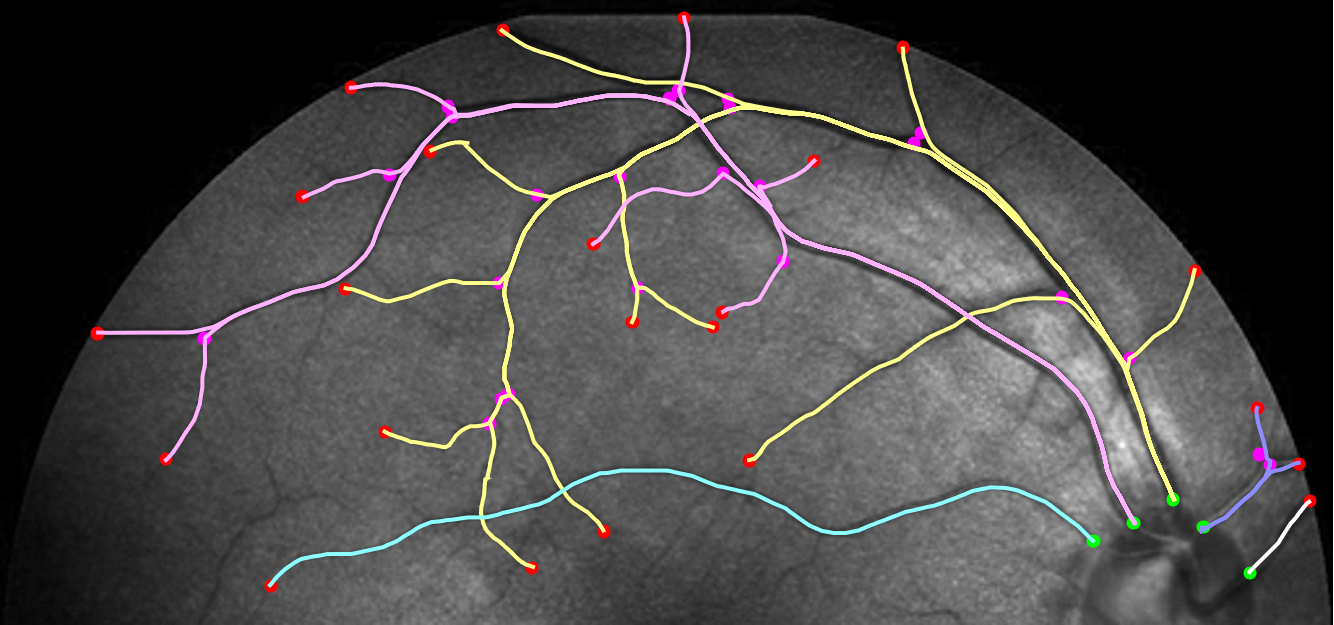}
        \caption{Image 1: Tracking with the mixed model with $\lambda=50$.}
    \end{subfigure}
    \begin{subfigure}[b]{\textwidth}
    \centering
    \includegraphics[width=\textwidth]{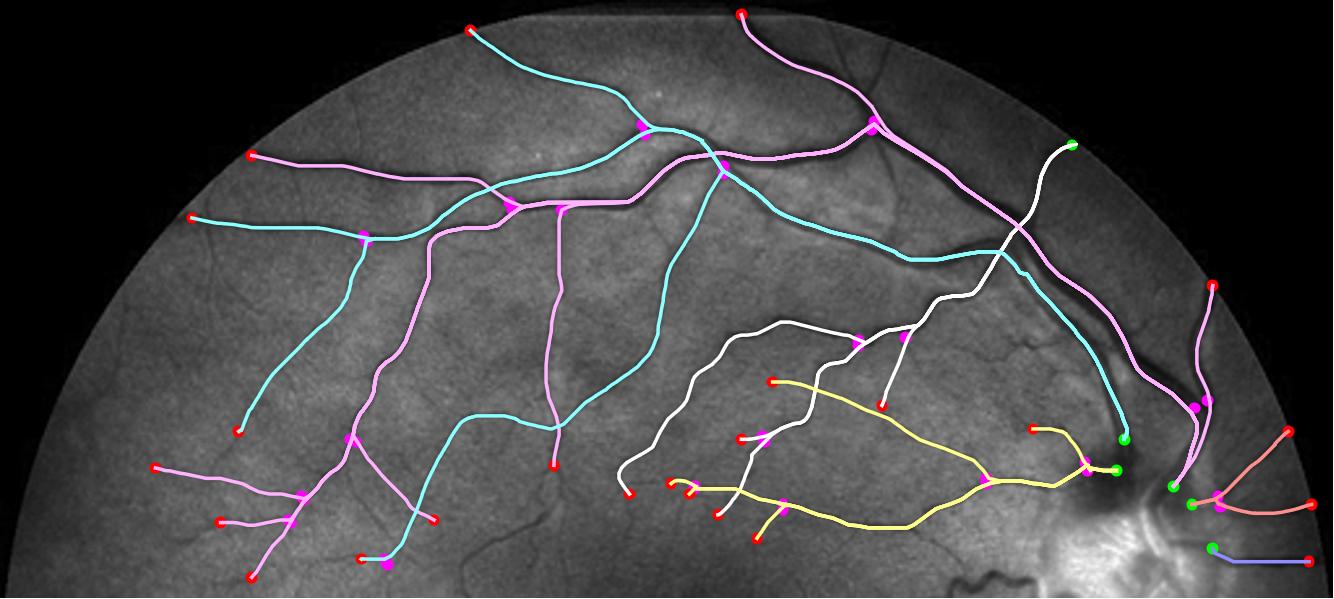}
        \caption{Image 2: Tracking with the mixed model with $\lambda=15$.}
    \end{subfigure}
    \caption{\textbf{Tracking of Vascular Tree per Seed:} Tracking with the mixed model $(\mathbb{M}_2,\mathcal{F}^M)$ proposed in \eqref{eq:MixedLIMTF}. Prior grouping of tips (in red) and seeds (in green) results in perfect tracking of the vessel tree, using only one efficient anisotropic fast-marching run via the numerical method in Section~\ref{sec:ExtensionAFM}. Both results are calculated with the cost function visualized in Fig.~\ref{fig:CostFunctionBifSeed}. At the purple points, we have bifurcations and our tracking is solely based on the mixed model produced (spatially projected) geodesics $\gamma$ with (automatic) in-place rotations at such bifurcations.}
    \label{fig:TrackingPerTreeLIFvsDDLIF}
\end{figure*}

Common practical setups in vascular tracking of retinal images include the prior knowledge of the locations of tips and seeds of vessel structures. We implemented our data-driven model using a prior classification of the  connected tips and seeds. More specifically, in every run of the fast marching algorithm, one of the seeds is considered together with its corresponding tips, and the connecting vessel structures are tracked. Figure~\ref{fig:TrackingPerTreeLIFvsDDLIF} shows our result in this setup and demonstrates that our approach can determine the geodesics that accurately follow the vascular structure in the retinal image.   
\subsection{Accuracy of the Model}
We now present some quantitative evaluations to measure the accuracy of our data-driven metrics for geodesic tracking. We measure the mistake ratio $E$ for the images in the STAR dataset. For these images, the ground truth of the vessels is known, which allows us to calculate the percentage of the vessel that is not on the correct vascular structure, where 
$$E=\frac{\text{\# pixels not on correct vessel}}{\text{\# pixels of all geodesics}}.$$

We have calculated this accuracy for images of the STAR dataset where we connect the tips to their nearest bifurcation, since one should use the new model away from crossing structures. The results are presented in Fig.~\ref{fig:PerformanceDDLIFvsLIF}. We see that for most tracks, the performance improves when switching to the new data-driven model, and in the cases where there is no improvement, the results do not get significantly worse. On average, we find an improvement of 10.7\% of the calculated tracks for the considered images.
\begin{figure}
    \centering
         \centering
        \scalebox{.8}{\begin{tikzpicture}
\begin{axis}[
    xlabel = {LIF},
    ylabel= {Mixed Model},
    xtick={0,.1,...,1},
    ytick={0,.1,...,1},
    scale only axis=true,
    axis equal image,
    scatter/classes={%
        a={mark=o,draw=black},
        b={mark=x,draw=black}}]
    \addplot[scatter,only marks,%
        scatter src=explicit symbolic]%
    table[meta=label] {
    LIF Mixed label
    .326934 .210185 a
    .193913 .193609 a
    .385554 .386654 b
    .445976 .416438 a
    .385127 .285812 a
    .400008 .372329 a
    .476843 .445372 a
    .294477 .238829 a
    .478078 .484434 b
    .413168 .37883 a
        };
    \path[name path = A] (axis cs:0,0) -- (axis cs:1,1);
    \path[name path = B] (axis cs:1,0) -- (axis cs:1,1);
    \path[name path = C] (axis cs:0,0) -- (axis cs:0,1);
    \addplot+ [green, fill opacity=0.2] fill between [of=A and B,soft clip={domain=0:1}];
    \addplot+ [red, fill opacity=0.2] fill between [of=A and C,soft clip={domain=0:1}];
    \draw[dotted] (0,0)--(1,1);
\end{axis}
\end{tikzpicture}}
        \caption{Visualization of the scatterplot of the accuracy of the mixed model vs. the left-invariant model applied on images in the STAR dataset (1,2,8,9,13,15,16,24,38,48). The accuracy is calculated on the calculated tracks between the tips and the nearest bifurcation for all vessels in one single run. The red area marks where the former left-invariant model performs better than the new mixed model, incidences indicated by an `x'. The green area marks where the new mixed model performs better than the former left-invariant model, incidences are indicated by an `o'. Most measurements show the improved performance of the mixed models compared to the left-invariant model (LIF).}
        \label{fig:PerformanceDDLIFvsLIF}
\end{figure}




\section{Conclusion and Future Work}
\label{sec:conclusion}
In this article, we introduced the concept of a data-driven left-invariant metric tensor field $\mathcal{G}^U$ and its asymmetric variant $\mathcal{F}^U$. The metric tensor field is defined by the underlying image, where movement along line structures is encouraged by its design in (\ref{Finslerfin}). In addition, a data-driven version $\nabla^U$ of the plus Cartan connection, relying on $\mathcal{G}^U$, was introduced. 

We used these geometrical tools to formulate a challenging data-driven version of \cite[Thm.1]{duits2021Springer}, which was stated and proved in Theorem~\ref{th:th}. In this theorem, `straight' and `short' curves are described with respect to the data-driven Cartan connection. In particular, it describes the entire Hamiltonian flow of the new Riemannian manifold model $(\mathbb{M}_{2},\mathcal{G}^U)$ in terms of the new data-driven Cartan connection $\nabla^U$, and explains the backtracking procedure for backtracking data-driven left-invariant geodesics in $(\mathbb{M}_{2},\mathcal{G}^U)$. As subsequently explained this can be extended to the asymmetric Finsler model $(\mathbb{M}_{2},\mathcal{F}^U)$ that
often yields the same geodesics, but also automatically places in-place rotations. The latter is beneficial at bifurcations in complex vasculature when using crossing-preserving vesselness costs for the costfunction $C$.

The diagonalization of the new data-driven left-invariant models $\mathcal{G}^U$ and $\mathcal{F}^U$ provides locally adaptive frames that are beneficial over previous approaches to locally adaptive frames in $\mathbb{M}_2$ \cite{duits2016locally,Zhang2016,smets2021total} in the sense that:\begin{enumerate}
    \item they coincide with the usual left-invariant frame if the data is locally constant,
    \item they are more stable as they are constructed by coercive metric tensor fields, recall (\ref{coerc}).
\end{enumerate}

To calculate the minimizing geodesics efficiently, an adaptation to the efficient anisotropic fast marching algorithm was required and presented in Section~\ref{sec:ExtensionAFM}. The metric tensor component matrix was no longer of block form in the fixed coordinate system, and the necessary changes to overcome this have been discussed and analyzed in Section~\ref{sec:ExtensionAFM}. 
We also provide an asymptotic error analysis of our numerical scheme. 

To show the performance of the data-driven metric tensor field and the mixed metric tensor field, we have tested them on 2D images of the retina. All experiments confirm that the new model is better able to adapt for curvature. In addition to that, for the tracking of a single vessel, a low number of orientations is sufficient to find the correct minimizing geodesic, as can be seen in Figure~\ref{fig:S-curve}. Full vascular tree tracking needs to be handled with care at difficult crossings structures, which is done in the mixed model $\mathcal{F}^M$ introduced in \eqref{eq:MixedLIMTF}. 

In general, the tracking results perform very well in the discussed two-step approach (see e.g. Figure~\ref{fig:LIFvsMM2Step}), where tips are first connected to the nearest bifurcation, after which the geodesics connecting these bifurcations to the corresponding seeds are calculated. After prior classification of seeds and tips belonging to the same vascular tree, the tracking results follow the vessels perfectly, recall Figure~\ref{fig:TrackingPerTreeLIFvsDDLIF}.

Despite some very appealing theoretical and practical advantages of our model, we still require considerable computation and runtime (tripling the overall processing time) to make the data-driven metric-tensor field and distance maps. Therefore, the exact usage of the proposed data-driven metric depends on the specific context of the tracking requirements.

For future work, it would be interesting to look into the possibilities to train the cost function $C$ using PDE-G-CNNs \cite{smets2020pde}, which is now geometrically computed as explained in Appendix~\ref{app:CostFunction}. In the past, this method had promising results for vessel segmentation. Besides using PDE-G-CNNs to construct the cost function, it would be worth looking into the possibilities to use neural networks to calculate the distance function as was done in \cite{DeepEikonalSolvers}.

\section*{Acknowledgements}

We gratefully acknowledge former TU/e master student A.J.~Wemmenhove\footnote{See A.J.~Wemmenhove's master thesis: {\scriptsize \url{https://pure.tue.nl/ws/portalfiles/portal/173216892/Wemmenhove_Jelle.pdf}}} (supervised by R.~Duits and B.M.N.~Smets) for preliminary 
investigations in support of Appendix~\ref{app:th}. 


\renewcommand{\appendixname}{Appendix}
\begin{appendices}

\section{Proof of Lemma~\ref{lemma:ConnectionInComponents}}\label{app:lemma}
Writing out the definition \eqref{eq:connection} gives
\begin{align*}
    \left(\nabla^{U}\right)_X Y&=\sum_{k=1}^n\left(\sum_{i=1}^n\omega_U^i(X)\mathcal{A}_{i}^U\omega_U^k(Y)\right.\\
    &\qquad\qquad\left.+\sum_{i,j=1}^n \omega_U^i(X)\omega_U^j(Y)\tilde{c}_{ij}^k(\cdot) \right)\mathcal{A}_k^U\\
    &=\sum_{k=1}^n\left(\sum_{i=1}^n \tilde{x}^i\mathcal{A}_{i}^U\tilde{y}^k+\sum_{i,j=1}^n\tilde{x}^i\tilde{y}^j\tilde{c}_{ij}^k(\cdot)\right)\mathcal{A}_k^U\\
    &=\sum_{k=1}^n\left(\dot{\tilde{y}}^k+\sum_{i,j=1}^n\tilde{c}_{ij}^k(\cdot)\tilde{x}^i\tilde{y}^j\right)\mathcal{A}_k^U,
\end{align*}
where the last equality holds since
\begin{align*}
    \dot{\tilde{y}}^k(t)&=\frac{\mathrm{d}}{\mathrm{d}t}\tilde{y}^k(\gamma(t))=\sum_{i=1}^n\dot{\tilde{\gamma}}^i\left(\mathcal{A}_{i}^U\tilde{y}^k\right)\!\left(\gamma(t)\right)\\&=\sum_{i=1}^n\tilde{x}^i\left(\mathcal{A}_{i}^U\tilde{y}^k\right)\!\left(\gamma(t)\right)=X\left(\tilde{y}^k\right)\!(\gamma(t)).
\end{align*}
Similarly, we have for \eqref{eq:dualConnection}
\begin{align*}
    \left(\nabla^{U}\right)_X^* \lambda&=\sum_{i=1}^n\left(\sum_{j=1}^n\omega_U^j(X)\mathcal{A}_j^U\mathcal{A}_{i}^U(\lambda)\right.\\
    &\qquad\qquad\left.+\sum_{j,k=1}^n\omega_U^j(X)\mathcal{A}_k^U(\lambda)\tilde{c}_{ij}^k(\cdot)\right)\omega_U^i\\
    &=\sum_{i=1}^n\left(\sum_{j=1}^n\tilde{x}^j\mathcal{A}_j^U\tilde{\lambda}_i+\sum_{j,k=1}^n\tilde{x}^j\tilde{\lambda}_k\tilde{c}_{ij}^k(\cdot)\right)\omega_U^i\\
    &=\sum_{i=1}^n\left(\dot{\tilde{\lambda}}_i+\sum_{j,k=1}^n\tilde{c}_{ij}^k\tilde{x}^j\tilde{\lambda}_k\right)\omega_U^i,
\end{align*}
where the last equality holds since
\begin{align*}
    \dot{\tilde{\lambda}}_i(t)&=\frac{\mathrm{d}}{\mathrm{d}t}\tilde{\lambda}_i(\gamma(t))=\sum_{j=1}^n\dot{\tilde{\gamma}}^j\left(\mathcal{A}_j^U\tilde{\lambda}_i\right)\!(\gamma(t))\\
    &=\sum_{j=1}^n\tilde{x}^j\left(\mathcal{A}_j^U\tilde{\lambda}_i\right)\!(\gamma(t))=X\left(\tilde{\lambda}_i\right)\!(\gamma(t)).\qquad\qquad\qed
\end{align*}

\section{Proof of Theorem~\ref{th:th}}\label{app:th}
{\color{blue} Firstly, we address the characterisation of straight curves.} \\[6pt]
$\left(\nabla^{U}\right)_{\dot{\gamma}}\dot{\gamma}=0 \Leftrightarrow $ \\    $\ddot{\tilde{\gamma}}^k+\sum \limits_{i,j=1}^{n} {c}_{ij}^k(\gamma(t))\dot{\tilde{\gamma}}^i\dot{\tilde{\gamma}}^j=0$ for $k=1,\ldots,n$.
    
    Since $[\mathcal{A}_{i}^U,\mathcal{A}_j^U]=-[\mathcal{A}_j^U,\mathcal{A}_i^U]$ due to \eqref{eq:LieBracket} and \eqref{eq:StructureConstants}, we have $\tilde{c}_{ij}^k=-\tilde{c}_{ji}^k$ and $\tilde{c}_{ii}^k=0$. \\
    Consequently, we are left with $\ddot{\tilde{\gamma}}^k=0\Rightarrow\ \dot{\tilde{\gamma}}^k-c^k=0$.
    Summarizing we have
    \[
    \left(\nabla^{U}\right)_{\dot{\gamma}}\dot{\gamma}=0\Leftrightarrow\dot{\gamma}=\sum \limits_{k=1}^n\dot{\tilde{\gamma}}^k\mathcal{A}_k^U=\sum \limits_{k=1}^n c^k\mathcal{A}_k^U.
    \]  
{\color{blue} Secondly, we address the characterisation of shortest curves.} \\[6pt]    
    The Pontryagin Maximum Principle is given by (Hamiltonian flow on co-tangent bundle $T^*(G)$)
    \begin{align*}
        &\begin{cases}
            \dot{\nu}=\overrightarrow{\mathfrak{h}}(\nu)\\
            \nu(0)=(\gamma(0),\lambda(0)),
        \end{cases}\numberthis\label{eq:PMP}\\
        &\text{where }\nu(t)=(\gamma(t),\lambda(t))\in T^*(G),\lambda(t)\in T_{\gamma(t)}^*(G).
    \end{align*}
    For details see \cite{AgrachevSachkov}.
    \begin{remark}[Liouvilles Theorem]
        By the chain law, for any smooth function $f:T^*(G)\to\mathbb{R}$, one has 
    \begin{align*}
        \frac{\mathrm{d}}{\mathrm{d}t}f(\gamma(t),\lambda(t))
        &=\nabla_{\gamma} f(v(t))\cdot\nabla_{\lambda}\mathfrak{h}(v(t))\\
        &\qquad-\nabla_{\lambda}f(v(t))\cdot\nabla_{\gamma}\mathfrak{h}(v(t))\\
        &=\left.\left\{\mathfrak{h},f\right\}\right|_{\gamma(t)},\numberthis\label{eq:Liouville}
    \end{align*}
    where $\mathfrak{h}$ denotes the Hamiltonian flow, and where $\nabla_{\gamma}$ and $\nabla_{\lambda}$ denote the gradient with respect to $\gamma$ and $\lambda$ respectively.
    \end{remark}
    The result of Eq.~\eqref{eq:PMP} can be expressed in coordinates
    \begin{align*}
        \begin{cases}
            \gamma(t)=(x(t),y(t),\theta(t))\\
            \lambda(t)=\sum \limits_{i=1}^n\tilde{\lambda}_i\left.\omega_U^i\right|_{\gamma(t)}\\
            \dot{\gamma}(t)=\sum \limits_{i=1}^n\dot{\tilde{\gamma}}(t)\mathcal{A}_{i}^U(\gamma(t)).
        \end{cases}
    \end{align*}
    The horizontal part is given by
    \begin{align}
        \dot{\tilde{\gamma}}^i=\langle\omega_U^i,\dot{\gamma}\rangle=\tilde{\lambda}^i\quad i=1,\ldots,n.\label{eq:HorizontalPart}
    \end{align}
    This follows from the computation of the Hamiltonian via the Fenchel transform.
    
    The vertical part is given by
    \begin{align*}
        \dot{\tilde{\lambda}}_i&\overset{\eqref{eq:Liouville}}{=}\left\{\mathfrak{h},\tilde{\lambda}_i\right\}=\left\{\frac{1}{2}\sum_{j=1}^n \tilde{\lambda}_j\tilde{\lambda}^j,\tilde{\lambda}_i\right\}=\left\{\frac{1}{2}\sum_{j=1}^n \alpha^j(\cdot)\tilde{\lambda}_j^2,\tilde{\lambda}_i\right\}\\
    &=\sum_{j=1}^n\tilde{\lambda}_j\alpha^j(\cdot)\left\{\tilde{\lambda}_j,\tilde{\lambda}_i\right\}+\frac{1}{2}\sum_{j=1}^n\tilde{\lambda}_j^2\left\{\alpha^j(\cdot),\tilde{\lambda}_i\right\}\\
        &=\sum_{j=1}^n\sum_{k=1}^n\tilde{c}_{ji}^k\alpha^j(\cdot)\tilde{\lambda}_j\tilde{\lambda}_k=\sum_{j,k=1}^n\tilde{c}_{ji}^k\tilde{\lambda}_k\tilde{\lambda}^j.
    \end{align*}
    In the above derivation, we have used that \[\left\{\tilde{\lambda}_i,\tilde{\lambda}_j\right\}=\sum_{k=1}^n\tilde{c}_{ij}^k(\cdot)\tilde{\lambda}_k.\]
Additionally, it is important to note that
\begin{align*}
    \{\alpha^j(\cdot),\tilde{\lambda}_i\}&=\sum_{k=1}^n\mathcal{A}_k^U\alpha^j(\cdot)\frac{\partial\tilde{\lambda}_i}{\partial\tilde{\lambda}_k}-\frac{\partial\alpha^j(\cdot)}{\partial\lambda_k}\mathcal{A}_k^U\tilde{\lambda}_i\\
    &=\mathcal{A}_{i}^U\alpha^j(\cdot),\numberthis\label{eq:SolvedPoissonBrackets}
\end{align*}since $\mathcal{A}_k^U\tilde{\lambda}_i=0$. Consequently, we find
\begin{align*}
    \sum_{j=1}^n\frac{1}{2}\tilde{\lambda}_j^2\left\{\alpha^j(\cdot),\tilde{\lambda}_i\right\}&\overset{\eqref{eq:SolvedPoissonBrackets}}{=}\sum\limits_{j=1}^n\frac{1}{2}\tilde{\lambda}_j^2\mathcal{A}_{i}^U \alpha^j(\cdot)\\
    &\overset{\eqref{eq:Hamiltonian}}{=}\mathcal{A}_{i}^U\mathfrak{h}(\lambda)=0,
\end{align*} since the Hamiltonian is constant.
    So in total 
    we have
    \begin{align*}
        \begin{cases}
            \dot{\tilde{\gamma}}^i=\tilde{\lambda}^i(t)\\
            \dot{\tilde{\lambda}}_i(t)=\sum\limits_{j=1}^n\sum\limits_{k=1}^n\tilde{c}_{ji}^k(\gamma(t))\; \tilde{\lambda}_k(t)\tilde{\lambda}^j(t).
        \end{cases}
    \end{align*}
    This is equivalent to
    \begin{align*}
        &\begin{cases}
            \dot{\tilde{\gamma}}^i(t)=\left(\left(\left.\mathcal{G}^U\right|_{\gamma(t)}\right)^{-1}\tilde{\lambda}(t)\right)^i\\
            \dot{\tilde{\lambda}}_i(t)-\sum \limits_{j,k=1}^n\tilde{c}_{ji}^k(\gamma(t))\tilde{\lambda}_k(t)\tilde{\lambda}^j(t)=0
        \end{cases}\\
        \Leftrightarrow&\begin{cases}
            \dot{\gamma}(t)=\left(\left.\mathcal{G}^U\right|_{\gamma(t)}\right)^{-1}\lambda(t)\\
            \left(\nabla^{U}\right)_{\dot{\gamma}}^*\lambda=0
        \end{cases}\\
        \Leftrightarrow&\begin{cases}
            \lambda(t)=\left(\left.\mathcal{G}^U\right|_{\gamma(t)}\right)\dot{\gamma}(t)\\
            \left(\nabla^{U}\right)_{\dot{\gamma}}^*\lambda=0.
        \end{cases}
    \end{align*}
\noindent
\begin{remark}[Justification of Eq.~\eqref{eq:HorizontalPart}]
    We have Lagrangian
    \begin{align*}
        \mathcal{L}(\gamma,\dot{\gamma})=\frac{1}{2}\sum \limits_{i=1}^n \alpha_i(\cdot)\left(\dot{\tilde{\gamma}}^i\right)^2.
    \end{align*}
    From this expression, we can determine the Hamiltonian $\mathfrak{h}$:
    \begin{align*}
        \mathfrak{h}(\gamma,\lambda)&=\sup_{\dot{\gamma}\in T_\gamma(G)}\left\{\langle\lambda,\dot{\gamma}\rangle-\mathcal{L}(\gamma,\dot{\gamma})\right\}\\
        &=\max_{(v^1,\ldots,v^n)}\left\{\sum_{i=1}^n\tilde{\lambda}_iv^i-\frac{1}{2}\alpha_i(\cdot) \left(v^i\right)^2\right\}.
    \end{align*}
    The maximum is calculated by taking the derivative w.r.t. $v$:
    \begin{align*}
        \begin{pmatrix}
            \tilde{\lambda}_1-\alpha_{1}(\cdot) v_{\max}^i\\
            \vdots\\
            \tilde{\lambda}_n-\alpha_n(\cdot) v_{\max}^i\\
        \end{pmatrix}=\begin{pmatrix}
            0\\\vdots\\0
        \end{pmatrix}\Leftrightarrow\begin{pmatrix}
            \tilde{\lambda}_1\\\vdots\\
            \tilde{\lambda}_n
        \end{pmatrix}=\begin{pmatrix}
            v_1^{\max}\\
            \vdots\\
            v_n^{\max}
        \end{pmatrix},
    \end{align*} where \begin{align*}
        \tilde{\lambda}_i=\alpha_i(\cdot)\tilde{\lambda}^i \text{ and } v_i^{\max}=\alpha_i(\cdot) v_{\max}^i.
    \end{align*}
    Consequently, we find $\tilde{\lambda}^i=v_{\max}^i$ for $i=1,\ldots,n$, and the Hamiltonian is given by
    \begin{align*}
        \mathfrak{h}&=\sum \limits_{i=1}^n\tilde{\lambda}_i\tilde{\lambda}^i-\frac{1}{2}\alpha_i(\cdot)\left(\tilde{\lambda}^i\right)^2\\&=\sum_{i=1}^n\left(\tilde{\lambda}_i\tilde{\lambda}^i-\frac{1}{2}\tilde{\lambda}_i\tilde{\lambda}^i\right)=\frac{1}{2}\sum_{i=1}^n\tilde{\lambda}_i\tilde{\lambda}^i.
    \end{align*}
    Hence, we also have found that $\tilde{\lambda}_i=\dot{\tilde{\gamma}}_i$ and $\tilde{\lambda}^i=\dot{\tilde{\gamma}}^i$, which we aimed to show. Note that reformulation in a coordinate-free matter yields  \[
    \forall_{i \in \{1,\ldots,n\}}\dot{\tilde{\gamma}}^i =\tilde{\lambda}^i \Leftrightarrow 
    \left.\lambda\right|_{\gamma(t)} = \left.\mathcal{G}_{U}^{-1}\right|_{\gamma(t)} \dot{\gamma}(t), 
    \]
    for all $t \in [0,W(\g)]$. 
\end{remark}
\begin{remark}
In Theorem~\ref{th:th} we give a backtracking formulation (where geodesics go `down-hill' to the origin via steepest descent) where we rescaled time back to the interval $t \in [0,1]$ (as this is more practical). Similar to \cite[Thm.4, Eq.28]{duits2018optimal} this is done as follows: in the ODE backtracking for geodesics (\ref{BT}) we included an extra negative scaling factor $-W(\g)$ in comparison to all the canonical ODEs above.     
\end{remark}

{\color{blue}
Thirdly, we address the symmetry (\ref{symmetry}) of the data-driven Cartan connection.} \\[6pt]    
By construction of (\ref{Finslerfin}) and (\ref{eq:connection}), we have the correct symmetry (\ref{symmetry}).
Indeed from (\ref{Finslerfin}), it follows that
\[
\begin{array}{l}
\left.\mathcal{A}_{i}^{\mathcal{L}_{\q}U}\right|_{\g\p}=
(L_\g)_* \left.\mathcal{A}_{i}^{U}\right|_{\p} \textrm{ and } \left.
\omega^{i}_{\mathcal{L}_\q U} \right|_{\g \p} = (L_{\g})^* \left.\omega^{i}_U \right|_{\p}
\end{array}
\]
and via 
(\ref{eq:connection}) we get $(L_{\g})^*
\left(\nabla^{\mathcal{L}_{\g}U}\right)^*=
(\nabla^U)^*$,
where we use the fact that $\mathcal{G}^U$ is diagonal w.r.t. basis $\{\mathcal{A}_i^U\}$ and where we respectively applied the 
pushforward of a vector field, the pullback of a co-vector field, and the pullback of a dual connection.
\begin{remark}\label{remark:Pullback}
    In general the pullback $\Phi^* \nabla^*$ of a dual connection $\nabla^*$ on manifold $G$ under a smooth mapping $\Phi: G \to G$ is given by $(\Phi^* \nabla^*_X)(\Phi^* \lambda)=\Phi^*(\nabla^*_{\Phi_* X}\; \lambda)$ with $\lambda \in T^*(G)$ and $X \in T(G)$. 
\end{remark}
{\color{blue} Finally, we address the integration of geodesics and their symmetry.} \\[6pt]
Eq.~(\ref{BT}) follows by (\ref{relevant}). Here $\lambda={\rm d}W$ implies we must indeed set the following momentum components:
\[
\tilde{\lambda}^k= \mathcal{A}_{k}^U W(\gamma(\cdot))\ , \textrm{ for all }k=1,\ldots,n.
\]
Furthermore in (\ref{relevant}) we invert the data-driven left-invariant metric tensor field $\mathcal{G}^U$ (recall Eq.~(\ref{eq:DDMTFDiagonalization})) and express the velocities as \mbox{$\dot{\tilde{\gamma}}^{k}= g_{kk}^{-1} \lambda_k = g^{kk} \lambda_k$}. Then via Remark~\ref{rem:proofeq29toEq31} we obtain that the geodesics indeed follow by integration of the vector field $v(W)$ on $G$. Clearly, this vector field is data-driven left-invariant (as all the vector fields $\mathcal{A}_{i}^U$ are) and thereby one has:
\[
\frac{d}{dt}\left( \gamma_{\q \g, \q \g_0}^{\mathcal{L}_{g}U} \right)(t)=
(L_{\q})_* \frac{d}{dt} \gamma_{\g, \g_0}^{U}(t),
\]
for all $\q,\g,\g_0 \in G$, $t \in [0,1]$, from which the symmetry (\ref{symgeods}) follows by integration. 
\qed

\section{The Used Metric Tensor Field is indeed a Data-Driven Left Invariant Metric Tensor Field}\label{app:GisDDLIV}

We first rely on a convenient standard formula of the Hessian of smooth function on a manifold relative to a connection on that manifold in Lemma~\ref{def:Hessian}. Then we provide an alternative formulation of such a Hessian in Lemma~\ref{lemma:AlternativeFormulationHessian} (via the notion of parallel transport). 

Finally, we prove that $\mathcal{G}^U$, that heavily relies (\ref{Finslerfin}) on a Hessian $HU$ of a sufficiently smooth orientation score $U:SE(2)\to \R$,  is indeed a data-driven left-invariant vector field in Lemma~\ref{lemma:DDLIF}.
\begin{lemma}\label{def:Hessian}
The Hessian $HU=\nabla^* {\rm d}U$ of a smooth function $U: M \to \R$ relative to connection $\nabla$ on manifold $M$ satisfies
\begin{equation} \label{eq:convenient}
HU(X,Y)= X(YU)- \nabla_X Y U.
\end{equation}
\end{lemma}
\begin{proof}
One can easily see that
\begin{align*}
    HU(X,Y)&=\nabla^*\mathrm{d}U(X,Y)=\left\langle\nabla_X^* \mathrm{d}U,Y\right\rangle\\
    &\overset{(\ref{dualCON})}{=}X\left\langle\mathrm{d}U,Y\right\rangle-\left\langle\mathrm{d}U,\nabla_X Y\right\rangle\\&=X(YU)-(\nabla_X Y)U.\qquad\qquad\qquad\qquad\qquad\qed
\end{align*}
\end{proof}
\begin{remark}[Alternative Formulation Hessian]
    Let $M$ be a smooth manifold with connection $\nabla$. Let $\p\in M$ and $X_\p,Y_\p\in T_\p(M)$, i.e. two tangent vectors not necessarily associated to a vector field. Let $f\in C^\infty(M,\R)$.
    
    Let $\mathcal{X}:[-\delta,\delta]\to M$, with $\delta>0$, such that
    \begin{align*}
        \begin{cases}
            \mathcal{X}(0)=\p\\
            \dot{\mathcal{X}}(0)=X_\p\\
            \nabla_{\dot{\mathcal{X}}(t)}\dot{\mathcal{X}}=0\qquad\forall t\in[-\delta,\delta],
        \end{cases}
    \end{align*}
    i.e. $\mathcal{X}$ is the unique autoparallel curve through $\p$ with tangent vector $X_\p$.
    For all $s,t\in[-\delta,\delta]$ let $P_{s,t}^\mathcal{X}: T_{\mathcal{X}(s)}M\to T_{\mathcal{X}(t)}M$ be the parallel transport operator along the curve $\mathcal{X}$, which is uniquely defined by the following properties
    \begin{enumerate}
        \item $P_{t,t}^\mathcal{X}=id\quad\forall t\in[-\delta,\delta]$,
        \item $P_{t_2,t_3}^\mathcal{X}\circ P_{t_1,t_2}^\mathcal{X}=P_{t_1,t_3}^\mathcal{X}, \quad\forall t_1,t_2,t_3\in[-\delta,\delta]$,
        \item smooth with respect to $\mathcal{X}$, $t$ and $s$.
    \end{enumerate}
    Then $t\mapsto P_{0,t}^\mathcal{X} Y_\p\in T_{\gamma(t)}M$ gives a smooth vector field along the curve $\mathcal{X}$ that is unique parallel transport of $Y_\p$ along that curve, i.e. with the property $\nabla_{\dot{\mathcal{X}}(t)}\left(P_{0,\cdot}^\mathcal{X} Y_\p\right)=0$.
\end{remark}
\begin{lemma}\label{lemma:AlternativeFormulationHessian}
  We can now define the Hessian of a (sufficiently) smooth function $f:M \to \R$ also as follows
    \begin{align*}
        \left.Hf\right|_\p(X_\p,Y_\p):&=\partial_t\left((P_{0,t}^\mathcal{X} Y_\p)f\right)(0)\\
        :&=\lim_{t\downarrow 0}\frac{(P_{0,t}^\mathcal{X} Y_\p)f-Y_\p f}{t}.\numberthis\label{eq:newHessian}
    \end{align*}
\end{lemma}
\begin{proof}
    If $X$ and $Y$ are smooth vector fields then
    \begin{align*}
        \left.(\nabla_X Y)\right|_\p=\nabla_{X_\p} Y=\lim_{t\to 0}\frac{P_{t,0}^\mathcal{X}Y_{\mathcal{X}(t)}-Y_\p}{t}.
    \end{align*}
    Note that $\lim_{t\to0} P_{t,0}^\mathcal{X}=id$, so
    \begin{align*}
        \lim_{t\to 0}\frac{Y_{\mathcal{X}(t)}-P_{0,t}^\mathcal{X}Y_\p}{t}&=\lim_{t\to 0}P_{t,0}^\mathcal{X}\frac{Y_{\mathcal{X}(t)}-P_{0,t}^\mathcal{X}Y_\p}{t}\\
        &=\lim_{t\to 0}\frac{P_{t,0}^\mathcal{X}Y_{\mathcal{X}(t)}-Y_\p}{t}=\nabla_{X_\p} Y.\numberthis\label{eq:TransportHessian}
    \end{align*}
    Now, we have
    \begin{align*}
        Hf(X_\p,Y_\p)&\overset{\eqref{eq:newHessian}}{=}\lim_{t\to0}\frac{P_{0,t}^\mathcal{X} Y_\p f-Y_\p f}{t}\\
        &\overset{\eqref{eq:TransportHessian}}{=}\lim_{t\to 0}\frac{(Yf)(\mathcal{X}(t))-(Yf)(\mathcal{X}(0))}{t}-\nabla_{X_\p}Yf\\
        &=X_\p(Yf)-\nabla_{X_\p}Yf,
    \end{align*}
    which is the same as Eq. \eqref{eq:convenient}.\qed
\end{proof}
\begin{lemma}\label{lemma:DDLIF}
    The metric tensor field $\mathcal{G}^U$ introduced in Eq. \eqref{Finslerfin} is a data-driven left invariant metric tensor field.
\end{lemma}
\begin{proof}
    We recall that the dual norm used in the definition of the data-driven metric tensor field $\mathcal{G}^U$ is given by
    \begin{align*}
        \|HU|_\p(\dot{\p},\cdot)\|_*=\sup_{\substack{Y\in T_\p(\mathbb{M}_2)\\\|Y\|=1}}|HU|_\p(\dot{\p},Y)|.
    \end{align*}
    In order to prove that $\mathcal{G}^U$ is a data-driven left invariant metric tensor field, we need the following identities:
    \begin{align}
        \left(\left(L_\g\right)_* Y\right)_{\g\p}\left(\mathcal{L}_\g U\right)=Y_\p (U),\label{eq:defPushForward}\\
        \nabla_{\left(L_\g\right)_* X} \left(L_\g\right)_* Y=\left(L_\g\right)_*\nabla_X Y,\label{eq:EquivCartanConnection}
    \end{align}
    where \eqref{eq:defPushForward} is the definition of the pushforward, and \eqref{eq:EquivCartanConnection} is the equivariance of the Cartan connection $\nabla=\nabla^{U=1}$. In addition, it is important that
    \begin{align}
        Y_\p\mapsto\left(L_\g\right)_* Y_\p\text{ is an isometry }T_\p(\mathbb{M}_2)\to T_{\g\p}(\mathbb{M}_2),\label{eq:Isometry}
    \end{align}
    so $\|Y\|=\|\tilde{Y}\|$, where $Y\in T_\p(\mathbb{M}_2)$ and $\tilde{Y}\in T_{\g\p}(\mathbb{M}_2)$.
     
    We check the data-driven left invariance for the separate terms of the metric tensor field, starting with $\mathcal{G}_\p(\dot{\p},\dot{\p})$:
    \begin{align*}
         \mathcal{G}_{\g\p}&\left(\left(L_\g\right)_*\dot{\p},\left(L_\g\right)_*\dot{\p}\right)\\
         &=\sum_{i,j=1}^n g_{ij}\left.\omega^i\right|_{\g\p}\left(\left(L_\g\right)_*\dot{\p}\right)\left.\omega^j\right|_{\g\p}\left(\left(L_\g\right)_*\dot{\p}\right)\\
         &=\sum_{i,j=1}^n g_{ij}\left.\omega^i\right|_\p\left(\dot{\p}\right)\left.\omega^j\right|_\p\left(\dot{\p}\right)=\mathcal{G}_\p\left(\dot{\p},\dot{\p}\right).
    \end{align*}
    Then, we study the data-driven left invariance of the term
    \begin{align*}
        \left\|\left.HU\right|_\p\left(\dot{\p},\cdot\right)\right\|_*^2&=\sup_{\substack{Y\in T_\p(\mathbb{M}_2)\\ \|Y\|=1}}\left|\left.HU\right|_\p\left(\dot{\p},Y\right)\right|^2\\
        &=\sup_{\substack{Y\in T_\p(\mathbb{M}_2)\\ \|Y\|=1}}\left|\dot{\p}(YU)-\mathrm{d}U\nabla_{\dot{\p}}Y\right|^2,
    \end{align*}
    which is satisfied since (set $\tilde{Y}=(L_\g)_* Y$)
    \begin{align*}
        &\left\|\left.H\left(\mathcal{L}_\g U\right)\right|_{\g\p}\left(\left(L_\g\right)_*\dot{\p},\cdot\right)\right\|_*^2\\
        &\qquad=\sup_{\substack{\tilde{Y}\in T_{\g\p}(\mathbb{M}_2)\\ \|\tilde{Y}\|=1}}\left|\left(L_\g\right)_*\dot{\p}(\tilde{Y}\mathcal{L}_\g U)\right.\left.-\mathrm{d}\mathcal{L}_\g U\nabla_{(L_\g)_*\dot{\p}}\tilde{Y}\right|^2\\
        &\qquad=\sup_{\substack{Y\in T_\p(\mathbb{M}_2)\\\|Y\|=1}}\left|\dot{\p}YU-\mathrm{d}U\nabla_{\dot{\p}}Y\right|^2=\left\|HU|_\p(\dot{\p},\cdot)\right\|_*^2.
    \end{align*}
    Note that since $\|HU|_\p(\dot{\p},\cdot)\|_*^2$ satisfies the data-driven left invariant property, we also have
    \begin{align*}
        \frac{\|H(\mathcal{L}_\g U)|_{\g\p}((L_\g)_*\dot{\p},\cdot)\|_*^2}{\max\limits_{\substack{\dot{\q}\in T_\q(\mathbb{M}_2)\\ \|\dot{\q}\|=1}}\|H(\mathcal{L}_\g U)|_{\g\p}((L_\g)_*\dot{\q},\cdot)\|_*^2}=\frac{\|HU|_{\p}(\dot{\p},\cdot)\|_*^2}{\max\limits_{\substack{\dot{\q}\in T_\p(\mathbb{M}_2)\\ \|\dot{\q}\|=1}}\|HU|_{\p}(\dot{\q},\cdot)\|_*^2}.
    \end{align*}
    Thereby, $\mathcal{G}^U$ is data-driven left-invariant:
    \begin{align*}
        \mathcal{G}_{\g\p}^{\mathcal{L}_\g U}&((L_\g)_*\dot{\p},(L_\g)_*\dot{\p})=\mathcal{G}_{\g\p}((L_\g)_*\dot{\p},(L_\g)_*\dot{\p})\\
        &\quad+\lambda\frac{\|H(\mathcal{L}_\g U)|_{\g\p}((L_\g)_*\dot{\p},\cdot)\|_*^2}{\max\limits_{\substack{\dot{\q}\in T_\p(\mathbb{M}_2)\\ \|\dot{\q}\|=1}}\|H(\mathcal{L}_\g U)|_{\g\p}((L_\g)_*\dot{\q},\cdot)\|_*^2}\\
        &=\mathcal{G}_\p(\dot{\p},\dot{\p})+\lambda\frac{\|HU|_\p(\dot{\p},\cdot)\|_*^2}{\max\limits_{\substack{\dot{\q}\in T_\p(\mathbb{M}_2)\\\|\dot{\q}\|=1}}\|HU|_\p(\dot{\q},
        \cdot)\|_*^2}=\mathcal{G}_\p^U(\dot{\p},\dot{\p}).\hfill\qed
    \end{align*}
\end{proof}
\section{Cost Function $C$: A Multi-scale Crossing-Preserving Vesselness Filter Variant}\label{app:CostFunction}
The differentiable cost function $C:\mathbb{M}_2\to [\delta,1]$ is an important tool used to encode where vascular structures are located. 
Costs are high (``1'') outside the blood vessels, and low ($\delta$ given in experiments) at the center of the blood vessels, stimulating the geodesic to move along the vascular structure. Many approaches to automatically calculate the vessel locations have been proposed over the years \cite{hannink2014crossing,Zhang2016,bekkers2015pde}. In the calculation of the tracking results, we use a cost function inspired by \cite{hannink2014crossing}. The precise relation between the cost and the vesselness filter follows later in \eqref{eq:CostFunction}. The vesselness expression \mbox{$\mathcal{V}^{SE(2)}(U_f^a):SE(2)\to\mathbb{R}^+$} is, as in \cite{hannink2014crossing,Frangi2000Multiscale}
\begin{align}
    \mathcal{V}^{SE(2)}(U_f^a)=\begin{cases}
        0 &\text{if }\mathcal{Q}\leq 0,\\
        \exp\left(\frac{-\mathcal{R}^2}{2\sigma_1^2}\right)\left(1-\exp\left(\frac{-\mathcal{S}^2}{2\sigma_2^2}\right)\right)&\text{if }\mathcal{Q}>0,
    \end{cases}\label{eq:multiscaleVesselness}
\end{align}
where $U_f^a(\x,\theta)$, $a>0$ fixed, is a single layer of a multilayer wavelet transform. In all experiments, we set $\sigma_1=0.5$ and $\sigma_2=0.5\left\|\mathcal{S}\right\|_\infty$. In \eqref{eq:multiscaleVesselness}, the anisotropy measure $\mathcal{R}$, structure measure $\mathcal{S}$ and convexity criterion $\mathcal{Q}$ are given by
\begin{align*}
    &\mathcal{R}\!=\!\left|\frac{\left(\!\mathcal{A}_{1}^2 U_f^a\!\right)^{\!s,\beta,\sigma_{s,Ext},\sigma_{a,Ext}}}{\left(\!\mathcal{A}_{2}^2U_f^a\!\right)^{\!s,\beta,\sigma_{s,Ext},\sigma_{a,Ext}}}\right|\!,
    \mathcal{Q}\!=\!\left(\!\mathcal{A}_{2}^2U_f^a\!\right)^{\!s,\beta,\sigma_{s,Ext},\sigma_{a,Ext}}\!\!,\\
    &\mathcal{S}\!=\!\sqrt{\!\left(\!\left(\!\mathcal{A}_{1}^2U_f^a\!\right)^{\!s,\beta,\sigma_{s,Ext},\sigma_{a,Ext}}\!\right)^{\!\!2}\!\!+\!\left(\!\left(\!\mathcal{A}_{2}^2U_f^a\!\right)^{\!s,\beta,\sigma_{s,Ext},\sigma_{a,Ext}}\!\right)^{\!\!2}}\!,\\
\end{align*}
where $\mathcal{A}_{i}^2 U_f^a:=\mathcal{A}_i\mathcal{A}_i U_f^a$, and where the superscripts $^{s,\beta}$ denote Gaussian derivatives at spatial scale $s=0.5\sigma_s^2$ and angular scale $0.5\beta^2$, where $\beta=0.75$, and where the superscripts $^{\sigma_{s,Ext},\sigma_{a,Ext}}$ denote external regularization with spatial scale $\sigma_{s,Ext}=\sigma_s=a$ and angular scale $\sigma_{a,Ext}$.

Here, we implement the dual norm $\|HU|_\p(\dot{\p},\cdot)\|_*^2$ is computed by \eqref{eq:dualnorm} using Gaussian derivatives with scales $\sigma_{s,Ext}$ and $\sigma_{a,Ext}$.

\begin{figure}
    \centering
    \includegraphics[width=0.48\textwidth]{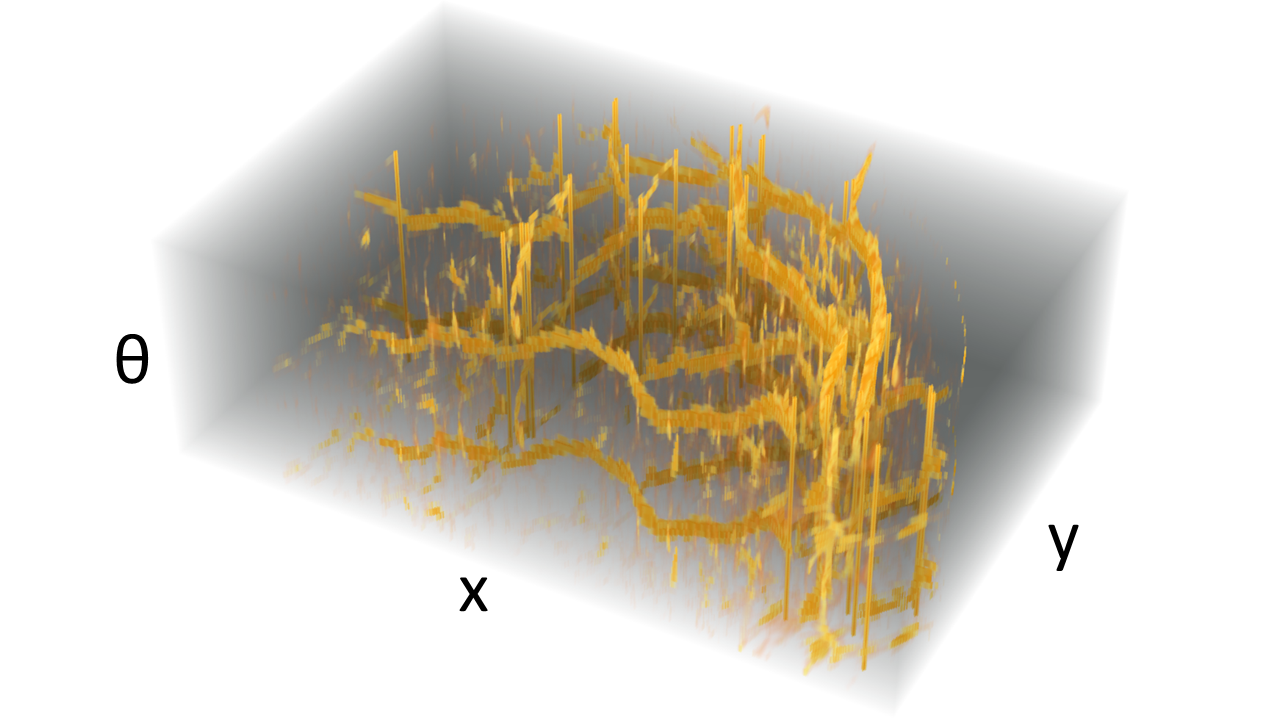}
    \caption{3D visualisation of a cost function $(x,y,\theta) \mapsto C(x,y,\theta)$ calculated with the introduced multi-scale crossing-preserving vesselness variant, of a retinal image $f$, with $\sigma_s\in\{1,2\}$ and $\sigma_{a,Ext}=0$.}
    \label{fig:3DcalculatedCostFunction}
\end{figure}

\begin{figure}
\centering
    \begin{subfigure}[b]{0.2\textwidth}
    \centering
    \includegraphics[width=\textwidth]{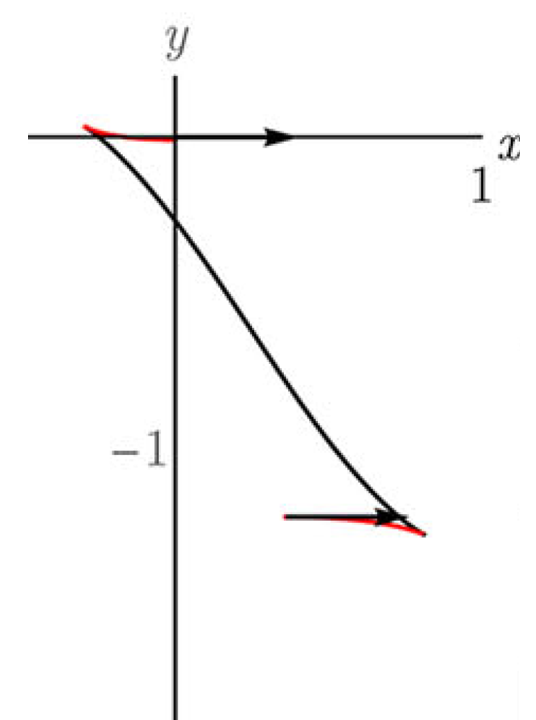}
        \caption{Cusp.}
        \label{fig:cusp}
    \end{subfigure}
    \begin{subfigure}[b]{0.2\textwidth}
    \centering
    \includegraphics[width=\textwidth]{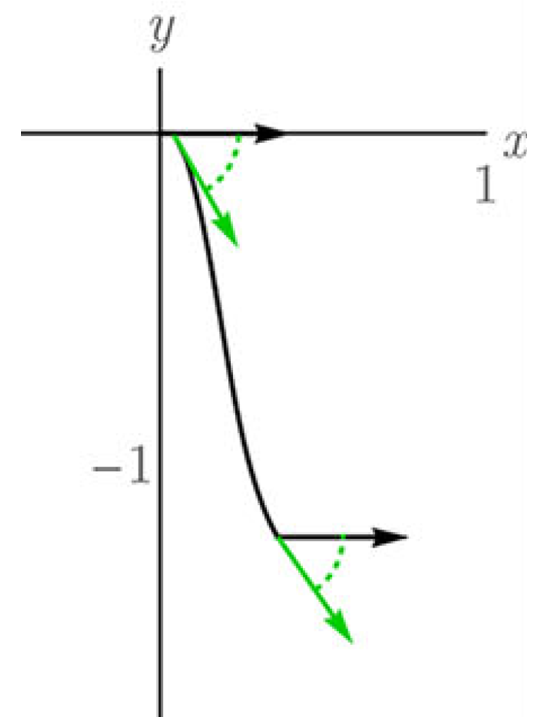}
        \caption{In-place rotation.}
        \label{fig:InPlaceRotation}
    \end{subfigure}
    \caption{Visualization of a cusp in a spatially projected geodesic in the Riemannian manifold $(\mathbb{M}_2,\mathcal{G}^{U=1})$ (left) and an in-place rotation in spatially projected Finslerian geodesic $(\mathbb{M}_2,\mathcal{F}^{U=1})$ (right).}
    \label{fig:cuspVsInPlaceRotation}
\end{figure}

Last, we apply erosion with scale $s_e$ on $\mathcal{V}^{SE(2)}$, denoting the result by $\mathcal{V}_{s_e}^{SE(2)}$. Then, the cost function \mbox{$C:SE(2)\to\mathbb{R}^+$} is defined as (similar to \cite{bekkers2015pde})
\begin{align*}
    &\left(\mathcal{V}^{SE(2)}\left(U_f\right)\right)(\x,\theta):=\mu_{\infty}^{-1}\sum_{l=1}^{N_s}\left(\mathcal{V}_{s_e}^{SE(2)}\left(U_f^{a_l}\right)\right)(\x,\theta),\\
    &C(\x,\theta):=\left(1+\lambda\left(\left(\mathcal{V}^{SE(2)}\left(U_f\right)\right)(\x,\theta)\right)^p\right)^{-1}\numberthis\label{eq:CostFunction}
\end{align*}
where $N_s$ denotes the number of scales, and $a_l$ denotes the different scales that are considered. The scaling parameter $\mu_\infty$ is defined as $\mu_\infty:=\left\|\sum_{l=1}^{N_s}\mathcal{V}_{s_2}^{SE(2)}\left(U_f^{a_l}\right)\right\|_\infty$. In all experiments, the values of parameters $\sigma_{a,Ext}$, $\lambda\geq 0$ and $p>0$ are chosen to be $0$, $1000$ and $2$ respectively. In Fig.~\ref{fig:3DcalculatedCostFunction}, a cost function constructed by the above formulation is visualized. Here, the vertical structures at bifurcations allow for in-place rotations, as depicted in Fig.~\ref{fig:cuspVsInPlaceRotation}.

\section{Adaptation to Asymmetric Data-Driven Finsler Functions} \label{app:AdaptationAsymmetricDataDrivenFinslerFunctions}
A generalization of \cite[Thm.~1,2,4]{duits2018optimal} to go from the symmetric model $(\mathbb{M}_2,\mathcal{G}^U)$ to the asymmetric model $(\mathbb{M}_2,\mathcal{F}^U)$, means in practice that we have to adapt the Eikonal equation \eqref{eq:EikonalPDE} to
\begin{align*}
    \left(\alpha_1^U(\cdot)\right)^{-1}\left((\mathcal{A}_1^U W)_+\right)^2+\sum_{j=2}^3 \left(\alpha_j^U(\cdot)\right)^{-1}\left(\mathcal{A}_j^U W\right)^2=1,
\end{align*}
where $(x)_+=\max\{x,0\}$, 
and backtracking \eqref{BTsimple} to
\begin{align*}
    \dot{\tilde{\gamma}}^1&=\frac{1}{W(g)}|\alpha_1^U|^{-1}\left((\mathcal{A}_1^U W)_+\right)\left(\tilde{\gamma}\right)\\
    \dot{\tilde{\gamma}}^k&=\frac{1}{W(g)}|\alpha_k^U|^{-1}\left(\mathcal{A}_k^U W\right)\left(\tilde{\gamma}\right), k=2,3.
\end{align*}
This adapted model does work reasonably well in practice. However, cusps may still occur in projected geodesics of  this adopted model since the required (`no reverse gear') condition
\begin{align} \label{reversegearcond}
    \dot{\gamma}^1=\dot{x}\cos\theta+\dot{y}\sin\theta=\dot{\x}\cdot\n(\theta)\geq 0
\end{align}
differs from the actually applied condition $\dot{\tilde{\gamma}}^1\geq 0$.
If the angle between $\left.\mathcal{A}_1^U\right|_{\tilde{\gamma}}$ and $\left.\mathcal{A}_1\right|_{\tilde{\gamma}}$ is not too large (which is often the case when geodesics pass locations with low cost), projected geodesics usually do not exhibit cusps.

A much less obvious and more precise solution -\emph{that does exclude cusps altogether}- in spatially projected data-driven geodesics, is given in Lemma~\ref{corol:rsf_eps_coefs} in the Numerical Section with computational scheme \eqref{eq:asym_quad_consistency} and backtracking \eqref{eq:BacktrackingSteepestDescent}. 

Essentially, in this approach (used in our experiments!) one takes the positive part of $\langle M^{-1}(\omega_U^1-\alpha \omega_U^2)|_{\tilde{\gamma}},\dot{\tilde{\gamma}}\rangle$ to ensure (\ref{reversegearcond}), rather than taking the positive part of $\dot{\tilde{\gamma}}^1=\langle \omega_U^1|_{\tilde{\gamma}},\dot{\tilde{\gamma}}\rangle$.

If $U$ is constant, then $\alpha=0$, $\omega^{1}_U=\omega^1$ and then we have $\dot{\tilde{\gamma}}^1=\dot{\gamma}^1\geq 0$, thereby if $U$ is constant we indeed end up in the standard Reeds-Shepp car model without reverse gear \cite{duits2018optimal}.

\section{Proof of Lemma~\ref{corol:rsf_eps_coefs}} \label{app:proofcor1}
Only in this proof we briefly write $\omega^i:=\omega^i_U$ for $i=1,2,3$. \\ 
Let $P$ be a $3\times k$ matrix of rank $k$ ($1\leq k\leq3$). Then, by the Woodbury formula, one has
\begin{align*}
    (\text{Id}_3+\epsilon^{-2}P P^\top)^{-1}&=\text{Id}_3-\epsilon^{-2} P(\textrm{Id}_k+\epsilon^{-2}P^\top P)^{-1}P^\top\\
    &=\textrm{Id}_3-P(P^\top P)^{-1}P^\top+\cO(\epsilon^2),
\end{align*}
which is up to $\cO(\epsilon^2)$ error equal to  the orthogonal projection $\Id_3-P(P^TP)P^T$ onto $\Span(P)^\perp$. Calculating the expression for $D_\epsilon$ by Woodbury formula and Taylor expansion gives:
\begin{align*}
    D_\epsilon:=(M^0+\epsilon \tilde{P}\tilde{P}^\top)^{-1}=(M^0)^{-\frac{1}{2}}\left(\textrm{Id}_3+\epsilon^{-2} PP^\top\right)^{-1}(M^0)^{-\frac{1}{2}},
\end{align*}
where $P:=(M^0)^{-1/2}\tilde{P}$ and  $\tilde{P}=\begin{bmatrix}\omega^1 &\omega^2\end{bmatrix}$. Noticing $\text{Span}\{P\}^\perp=\text{Span}\{\tilde{\omega}\}$, with $\tilde{\omega}:=(\sqrt{M^0}^{-1}\omega^1)\times(\sqrt{M^0}^{-1}\omega^2)$, one has
\begin{align*}
    D_\epsilon=(M^0)^{-1/2}\frac{\tilde{\omega}\tilde{\omega}^\top}{\tilde{\omega}^\top M^0\tilde{\omega}}(M^0)^{-1/2}+\cO(\epsilon^2).
\end{align*}
Simplification yields $D=\frac{\omega\omega^\top}{\omega^\top M^0\omega}+\cO(\epsilon^2),$ with $\omega:=\omega^1\times\omega^2$ as stated in \eqref{eq:Depsilon}.

Likewise, via Lemma~\ref{lem:asym_quad_dual}, the Woodbury formula, and applying a Taylor expansion, we obtain that 
\begin{align*}
    M_\epsilon^{-1} &= (M^0+\epsilon^{-2} \omega^2 (\omega^2)^\top)^{-1} \\
    & =(M^0)^{-1/2}(\textrm{Id}_3-\tilde{\omega}^2((\tilde{\omega}^2)^\top\tilde{\omega}^2)^{-1}(\tilde{\omega}^2)^\top)(M^0)^{-1/2}\\
    &\qquad+\cO(\epsilon^2),
\end{align*}
where $\tilde{\omega}^2:=(M^0)^{-1/2}\omega^2$ and $$(\textrm{Id}_3-(M^0)^{-1/2}\omega^2((\omega^2)^\top(M^0)^{-1}\omega^2)^{-1}(\omega^2)^\top (M^0)^{-1/2})$$ is up to $\cO(\epsilon^2)$ error the orthogonal projection onto\\ $\Span\{(M^0)^{-1/2}\omega^2\}^\perp$. The further simplification of
\begin{equation*}
    \eta_\epsilon=\frac{\epsilon^{-1}M_\epsilon^{-1}\omega^1}{\sqrt{1+\epsilon^{-2}(\omega^1)^\top M_\epsilon^{-1}\omega^1}}
\end{equation*}
boils down to \eqref{eq:etaepsilon}.\qed
\end{appendices}


%
%
\section*{Declarations}
\subsection*{Competing Interests}
R.~Duits is a member of the editorial board of JMIV.
\subsection*{Authors' Contributions}
N.J.~van~den~Berg developed the used model, wrote the first draft of the manuscript, did all the experiments, wrote the code for the experimental section supported by B.M.N.~Smets, J.-M.~Mirebeau and R.~Duits (available via GitHub \cite{githubNicky}), and created all figures. She also contributed to important parts of the proof of Theorem~\ref{th:th} (App.~\ref{app:th}). 
B.M.N.~Smets proposed the data-driven metric tensor model (beneficial over previous choices of locally adaptive frames in the context of image enhancement \cite{Duitsbookchapter}) on which the paper is built. He developed Lemma~\ref{lemma:AlternativeFormulationHessian}, and contributed to Lemma~\ref{lemma:ConnectionInComponents}. 
G.~Pai acted as co-supervisor of the project, providing relevant input and suggestions to experiments, helped with careful exposition of the results, and structuring of the manuscript. 
J.-M.~Mirebeau developed the new anisotropic fast-marching algorithm and theoretical analysis (Lemma~\ref{corol:rsf_eps_coefs} and asymptotic expansion in \eqref{eq:scheme_consistent_0}), which is explained in 
Section~\ref{sec:ExtensionAFM}, and provided relevant input on Theorem~\ref{th:th}. 
R.~Duits is the supervisor of the project, proposed and build up the underlying theory (Section~\ref{ch:model}~and~\ref{sec:ThShortVsStraight}, Theorem~\ref{th:th} \& Lemma \ref{lemma:DDLIF}), and polished the manuscript on mathematical details.

All authors collaborated closely, wrote parts of the manuscript and reviewed the manuscript.

Author contributions per section (in order of appearance):
Sec.~\ref{sec:Introduction} - N.J.~van~den~Berg \& G.~Pai; 
Sec.~\ref{sec:background} - N.J.~van~den~Berg \& R.~Duits;
Sec.~\ref{ch:model} - N.J.~van~den~Berg, B.M.N.~Smets \& R.~Duits;
Sec.~\ref{sec:ThShortVsStraight} - N.J.~van~den~Berg \& R.~Duits;
Sec.~\ref{sec:ExtensionAFM} - J.-M.~Mirebeau \& R.~Duits;
Sec.~\ref{sec:experiments} - N.J.~van~den~Berg \& G.~Pai;
App.~\ref{app:lemma} - N.J.~van~den~Berg \& R.~Duits;
App.~\ref{app:th} - N.J.~van~den~Berg \& R.~Duits;
App.~\ref{app:GisDDLIV} - N.J.~van~den~Berg, B.M.N.~Smets \& R.~Duits;
App.~\ref{app:CostFunction} - N.J.~van~den~Berg \& R.~Duits;
App.~\ref{app:proofcor1} - N.J.~van~den~Berg \& J.-M.~Mirebeau.


\subsection*{Funding}
We gratefully acknowledge the Dutch Foundation of Science NWO for its financial support by Talent Programme VICI 2020 Exact Sciences (Duits, Geometric learning for Image Analysis, VI.C. 202-031).
\subsection*{Availability of Data and Materials}
The STAR-dataset \cite{Zhang2016,AbbasiSureshjani} is, together with the \textit{Mathematica} notebooks, available via \cite{githubNicky}.

\bibliographystyle{spmpsci}      
\bibliography{references}

%
%



\end{document}